\input amstex







\font\rm=cmr10 \rm

\font\bf=cmb10
\font\Rm=cmr9 at 11pt
\rm
\font\it=cmsl9 at 10pt
 at 7pt

\font\Rrm=cmr17 at 16pt
   \font\Rm=cmr12 at 11.5pt

\long\def\Pf{\par\noindent {\it Proof.} }
\def\({\left(}
\def\){\right)}
\def\st{such that }
\def\qed{\hfill$\bullet$\vskip 4pt}

\def\brcs#1{\left\{ #1\right\}}

\def\iso{\cong}
\def\wrt{with respect to }
\def\:{\,:}

\def\ker{\text{ker\,}}

\def\C{\text{\bf C}}

\def\I{\text{I\,}}

\def\R{\text{\bf R}}
\def\N{\text{\bf N}}
\def\Z{\text{\bf Z}}
\def\Q{\text{\bf Q}}

\def\Arrow #1;#2.{#1\:#2 \to }

\def\Set#1#2{\brcs{#1 \left|\vphantom{#1 #2} \right.#2}}



\def\Rrr#1,#2{{\Cal J}_{#1,#2}}
\def\slfrac#1#2{{\raise -.07 ex\hbox{$^{#1}$}}\!/\raise .35 ex \hbox{${}_{#2}$}}
\def\ssf #1/#2{\slfrac {#1}{#2}}

\def\pd #1,#2.{\frac {\partial #1}{\partial #2}}

   \long\def\Lem
#1.#2\par{\vskip4pt{\baselineskip=13pt\font\it=cmsl12 at
11.5pt\Rm
   \noindent {\rm \uppercase{#1}} #2\vskip3pt

   }} 

\long\def\Proclaim #1.#2 \endproclaim{\vskip4pt{\baselineskip=13pt\font\it=cmsl12 at
11.5pt\Rm
   \noindent {\rm \uppercase{#1}} #2\vskip3pt

   }} 

\long\def\remark #1\endremark{\vskip 2pt \noindent {\it Remark\/} #1\par}

\long\def\Sectionhead #1.#2:\par #3{\vskip 4pt \noindent {\bf #1 #2}vskip 2pt\noindent\nospace #3}

\long\def\Title #1\par {\noindent{\Rrm #1}\vskip 9pt}

 \long\def\SubT #1.{\noindent {\it #1\/} } 
 
 \long\def\SecT
#1\par{\vskip 3pt \noindent {\bf #1}\vglue1pt
   \noindent}

\long\def\subtitle #1.{\vskip 2pt \noindent {\it #1}}

\long\def\Rmk#1\par{\vskip 1pt \noindent {\it
Remark.} #1\vskip2pt}

\long\def\Abstract #1\par{{\leftskip= 3 true cm \rightskip = 3 true cm \font\it=cmsl10 \font\rm=cmr10 \baselineskip = 10pt
\parindent=.35 true cm\rm\noindent 
{\it Abstract} #1\vskip 8pt

}}

\long\def\Author #1 \par{\noindent{\it #1}}

\scrollmode\NoBlackBoxes
\magnification=1100
\long\def\Abstract #1\par%
{\vskip .2 true cm{\leftskip 1 true in \rightskip 1 true in \font\rm=cmr8 \rm
\baselineskip=1pt \font\it=cmsl8 \font\bf=cmb10 at 8pt
\parindent=0em {\bf Abstract} #1

}}
\comment
\font\rm=Times at 10pt

\font\bf=TimesB
\font\Rm=Times at 11pt
\rm
\font\it=TimesI at 10pt
\endcomment

\long\def\Pf{\par\noindent {\it Proof.} }
\def\({\left(}
\def\){\right)}
\def\st{such that }
\def\qed{\hfill$\bullet$\vskip 4pt}

\def\brcs#1{\left\{ #1\right\}}
\def\Set#1#2{\brcs{#1 \left|\vphantom{#1 #2} \right.#2}}

\def\C{\text{\bf C}}

\def\I{\text{I\,}}

\def\iso{\cong}
\def\wrt{with respect to }
\def\:{\,:}
\def\Arrow #1;#2.{#1\:#2 \to }


\def\R{\text{\bf R}}
\def\N{\text{\bf N}}
\def\Z{\text{\bf Z}}
\def\Q{\text{\bf Q}}
 
\def\Rrr#1,#2{{\Cal J}_{#1,#2}}

\def\slfrac#1#2{{\raise -.07 ex\hbox{$^{#1}$}}\!/\raise .35 ex \hbox{${}_{#2}$}}
\def\ssf #1/#2{\slfrac {#1}{#2}}

\def\pd #1,#2.{\frac {\partial #1}{\partial #2}}


   \long\def\Title #1\par {\noindent{\Rrm #1}\vskip 9pt}
 \long\def\SubT #1.{\noindent {\it #1\/} }   \long\def\SecT
#1\par{\vskip 3pt \noindent {\bf #1}\vglue1pt
   \noindent}
\long\def\subtitle #1.{\vskip 2pt \noindent {\it #1}}

\long\def\Rmk#1\par{\vskip 1pt \noindent {\it
Remark.} #1\vskip2pt}


\def\One{1}
\def\oneone{\One.1}
\def\onetwo{\One.2}
\def\onethr{1.3}
\def\onefou{1.4}
\def\onefiv{1.5}
\def\onesix{1.6}
\def\onesev{1.7}
\def\oneeig{1.8}
\def\onenin{1.9}

\def\twoone{3.1}
\def\twotwo{3.2}

\def\throne{2.1}

\def\fouone{4.1}
\def\foutwo{4.2}
\def\fouthr{4.3}
\def\foufou{4.4}
\def\foufiv{4.5}

\def\fivone{5.1}
\def\fivtwo{5.2}
\def\fivthr{5.3}
\def\fivfou{5.4}
\def\fivfiv{5.5}
\def\fivsix{5.6}
\def\fivsev{5.7}
\def\fiveig{5.8}
\def\fivnin{5.9}
\def\fivten{5.10}
\def\fivele{5.11}

\def\sixtwo{6.2}
\def\sixthr{6.3}

\def\sevone{7.1}
\def\sevtwo{7.2}
\def\sevthr{7.3}
\def\sevfou{7.4}
\def\sevfiv{7.5}
\def\sevsix{7.6}

\def\seveig{7.8}
\def\sevnin{6.3}
\def\sevten{6.4}
\def\sevele{6.5}


\long\def\Remark{\vskip 2pt \noindent{\it Remark\/ }}
\def\flo #1{\lfloor #1 \rfloor}

\def\rank{\text{{\rm rank}}\,}
\def\Q{\text{\bf Q}}
\def\Inf{\text{Inf\,}}
\def\Ann{\text{Ann\,}}

\let\hat=\widehat

\def\Aff{\text{Aff\,}}
\NoBlackBoxes

\def\One{1}
\def\oneone{\One.1}
\def\onetwo{\One.2}
\def\onethr{1.3}
\def\onefou{1.5}
\def\onefiv{1.6}
\def\onesix{1.7}
\def\onesev{1.9}
\def\oneeig{1.8}
\def\onenin{1.10}
\def\oneten{1.11}
\def\oneele{1.12}
\def\onetwe{1.4}

\def\twoone{2.1}
\def\twotwo{2.2}

\def\throne{3.1}

\def\fouone{4.1}
\def\foutwo{4.2}
\def\fouthr{4.3}
\def\foufou{4.4}
\def\foufiv{4.5}

\def\fivone{5.1}
\def\fivtwo{5.2}
\def\fivthr{5.3}
\def\fivfou{5.4}
\def\fivfiv{5.5}
\def\fivsix{5.6}
\def\fivsev{5.7}
\def\fiveig{5.8}
\def\fivnin{5.9}
\def\fivten{5.10}
\def\fivele{5.11}

\def\sixtwo{6.1}
\def\sixthr{6.2}

\def\sevone{7.3}
\def\sevtwo{B.2}
\def\sevthr{B.1}
\def\sevfou{B.3}
\def\sevfiv{B.2}
\def\sevsix{B.5}

\def\seveig{7.13}
\def\sevnin{B.4}
\def\sevten{7.12}
\def\sevele{7.2}
\def\sevtwe{7.4}
\def\sevthi{B.6}

\def\appone{A.1}
\def\apptwo{A.2}

\def\appbtwo{7.7}
\def\appbthr{7.1}
\def\appbfou{7.2}
\def\appbfiv{7.5}
\def\appbsix{7.6}
\def\appbsev{7.8}
\def\appbeig{7.9}
\def\appbnin{7.10}
\def\appbten{7.12}
\def\appbele{7.11}
\def\appbtwe{7.0}

\def\exone{1}
\def\extwo{2}
\def\exthr{3}
\def\exfou{4}
\def\exfiv{5}
\def\exsix{7}
\def\exsev{8}
\def\exeig{9}
\def\exnin{10}
\def\exten{11}

\def\extwe{6}
\def\exthi{13}

\Title Measures on Cantor sets: the good, the ugly,  the bad%
\plainfootnote{\Rm*}{\rm This is the correct translation of the title of the Italian spaghetti western ({\it Il buono, il brutto, il cattivo,} D: Sergio Leone), which was changed for US audiences. The English language title, {\it The good, the bad, and the ugly,} is clich\' ed now---over 175 articles in engineering and mathematics use it in their title---but we could not find any titles using the original order.}

\Abstract We translate   Akin's notion of {\it good\/} (and related concepts) from measures on Cantor sets to traces on dimension groups, and particularly for invariant measures of minimal homeomorphisms (and their corresponding simple dimension groups), this yields characterizations and examples, which translate back to the original context. Good traces on a simple dimension group are characterized by their kernel having  dense image in their annihilating set of affine functions on the trace space; this makes it possible to construct many examples with seemingly paradoxical properties. {\par}
In order to study the related property of {\it refinability,} we consider goodness for sets of measures (traces on dimension groups), and obtain partial characterizations in terms of (special) convex subsets  of Choquet simplices. {\par}These notions  also very closely related to unperforation of quotients of dimension groups by convex subgroups (that are not order ideals), and we give partial characterizations. Numerous examples illustrate the results.

\Author Sergey Bezuglyi \& David Handelman %
\plainfootnote{$^{1}$}{\rm Supported in part by an NSERC Discovery Grant.} 

\SecT Introduction

Let $X$ be a non-atomic separable zero dimensional compact set; we refer to it as a Cantor set. Let $\mu$ be a probability measure on $X$. Motivated by the problem of homeomorphism of measures,   Akin  initiated a systematic study  of  Borel probability measures on a Cantor set [Ak1, Ak2]. He showed that the  {\it clopen value set\/} $S(\mu)$  (the set of values of $\mu$ on clopen sets) plays
a significant role in classification of measures up to a homeomorphism. It is a countable dense subset of the unit interval, and this set provides an invariant for topologically equivalent measures, although not complete. However, for the class of {\it good\/} measures, $S(\mu)$ {\it is\/} a complete invariant. By definition, a  (full non-atomic) probability measure $\mu$ is good if whenever $U$, $V$ are clopen sets with $\mu(U) < \mu(V)$, there exists a clopen subset $W$ of $V$ such that $\mu(W) = \mu(U)$. It turns out that such measures are exactly the invariant measures which can arise in uniquely ergodic minimal homeomorphisms of Cantor sets [Ak2, GW].

A series of papers in the last decade, by Akin, Austin, Dougherty, Mauldin, Yingst ([ADMY, Au, DMY, Y], have  focused on various properties of Bernoulli and related measures.  In particular, the notion of {(\it weak\/{\rm)} refinability\/} (weaker than good) was formulated and studied.
More recently, [BK, K] have exhaustively studied invariant finite and infinite measures of aperiodic substitution dynamical systems.

Akin [A2] introduced the notion of a good measure, yielding criteria for two measures on $X$ to be homeomorphic (in the sense that there there is a self-homeomorphism of $X$ sending one measure to the other). Specifically if $\mu$ and $\mu'$ are good probability measures, then they are homeomorphic if and only if their sets of values on clopen sets are equal. Glasner and Weiss [GW] showed that if there is a minimal homeomorphism $T$, and $A$ and $B$ are clopen sets \st  $\mu(A) < \mu (B)$ for all invariant measures, then there exists an element, $\gamma$, of the topological full group  of $(X,T)$ \st$\gamma(A) \subset B$. A particular consequence is that if $T$ is uniquely ergodic, then the measure is good.

Akin  posed the question, if $\mu$ is a measure \st the group of homeomorphisms of $X$ that preserve $\mu$ acts minimally, then is  $\mu$ good? If this were true, then all invariant measures of minimal homeomorphisms would be good. Early in this paper, we show that extremal (also known as {\it ergodic\/ {\rm or} pure\/}) invariant measures of a minimal homeomorphism need not be  good, and are rarely so.  For example, we can arrange examples wherein the extremal measures are not good, or the set of good pure measures and the set of not good pure measures are simultaneously dense in the set of pure (that is, ergodic) measures. We also introduce ugly and bad measures (formulated in terms of traces on dimension groups).

The language of dimension groups permits us to characterize  good measures among the measures (now called  traces) in terms of the ordered group  and Choquet space structures, and this makes  construction of examples relatively easy. For example, if $G$ is a simple dimension group with order unit $u$ and $\tau$ is a pure trace, then the trace $\tau$ is good (that is, when we translate back to invariant measures on a Cantor set with a minimal homeomorphism, the corresponding invariant measure) if and only if
the image of $\ker \tau$ (a subgroup of $G$) has dense range in its the vector space of affine functions that vanish at $\tau$.

Based on work of [GPS1] and before that of Vershik [V], we can go back and forth between invariant measures of minimal homeomorphisms of Cantor sets and trace of simple dimension groups. In fact, this translation can be extended to some non-minimal homeomorphism, via either Vershik's adic map, or via pre-ordered K-theory of crossed product C*-algebras. In the former case, Medeyents [Me] has shown that any aperiodic homeomorphism admits a representation as a Vershik map on (generally non-simple) Bratteli diagrams.  Dimension groups and Choquet theory have a large repetoire of results and examples, and it seems convenient to translate the original problems of measures to this context.

When $(X,T)$ is a Cantor set with a self-homeomorphism, we can form the usual $H^1(X,T)$ (Cech cohomology), which can be identified with K$_0 \( C(X,\C) \times_T \Z\)$, the pre-ordered Grothendieck group of the C*-algebra crossed product. This can also be identified [BH, P; going back to [GPS] in the minimal case] with the quotient group $G \equiv G(X,T)= C(X,\Z)/(\I-T)C(X,\Z)$ (the kernel representing the coboundaries) with the quotient pre-ordering (that is, a coset is in the positive cone if it has a representative $f \in C(X,\Z)$ that is nonnegative). In many cases (minimal: [GPS]; shifts of finite type: [BH2]), this pre-ordering is actually a partial ordering (that is, $G^+ \cap -G^+ = {\pmb 0}$) and much more. We abbreviate the statement $g-h \in G^+$ to $h \leq g$ when the ordering is unambigously given. Moreover, the invariant probability measures on $(X,T)$ are in natural (affine) bijection with the normalized traces on $G$.

We define a {\it trace\/} on a partially ordered abelian group $G$ to be a nonzero positive group homomorphism $\Arrow \tau; G. \R$. An element $u$ of $G^+$ is an {\it order unit of $G$\/} if for all $g \in G$, there exists a positive integer $N$ \st $-Nu \leq g \leq Nu$. The set of order units will be denoted $G^{++}$. A trace $\tau$ is {\it normalized \wrt $\tau$\/} if $\tau(u) = 1$. We may abbreviate this to   {\it normalized trace of $(G,u)$\/}. The assignment $\mu \mapsto \tau_{\mu}$ where $\tau_{\mu}([\chi_U]) = \mu (U)$ (where $U$ is a clopen subset of $X$ and $\chi_U$ denotes the indicator function of $U$) actually yields an affine homeomorphism between the Choquet simplices of invariant measures and of normalized traces on the preordered abelian group $(G(X,T), [\chi_X])$ with order unit.

When $T$ is minimal, [GPS] showed that the resulting pre-ordered group is not only an ordered group, but a simple dimension group. A partially ordered abelian group is {\it simple\/} if every nonzero element of $G^+$ is an order unit. A partially ordered abelian group is a {\it dimension group\/} if it satisfies Riesz decomposition (whenever $0 \leq a \leq b + c$ where $b,c \in G^+$, there exist $0 \leq b' \leq b$ and $0 \leq c' \leq b$ \st $a = b' + c'$) and {\it unperforation\/} (if $n$ is a positive integer and $g$ is an element of $G^+$ \st $ng \in G^+$, then $g \in G^+$). Equivalently, by [EHS], a dimension group is a direct limit of simplicial groups $\Z^{n(\alpha)}$ ($\alpha$ varying over a directed set) with positive homomorphisms between them. For countable dimension groups, the index set in the limit of simplicial groups can be taken to be the positive integers.

Not only did [GPS] show that when $T$ is minimal, $G(X,T)$ is a simple dimension group, but they also showed that given a countable simple dimension with order unit, $(G,u)$, there exists a minimal $T$ on $X$ \st $(G(X,T), [\chi_X] ) $ is order isomorphic to $(G,u)$, and moreover, $(G(X,T), [\chi_X] ) $ is a complete invariant for strong orbit equivalence. 

Simple dimension groups are completely understood (if we think that we understand Choquet simplices and countable dense subgroups of Banach spaces). Specifically, if $(G,u)$ is a partially ordered abelian group, we may construct a compact convex set $S(G,u)$ inside $\R^G$ 
 consisting of the normalized traces. Then $(G,u)$ admits an affine representation, that is, the map $\Arrow \hat{}; G. \text{Aff\,} S(G,u)$, where $\widehat g (\tau) = \tau(g)$ (so that $\hat u $ is the constant function $\pmb 1$). For details on affine representations, Choquet theory, etc,  see [Al], [As], and [G, sections 5,8,9,11]. As in those references, the extremal boundary of $S(G,u)$ consisting of the extreme (or pure) points will be denoted $\partial_e S(G,u)$.

The map $\ \hat {}\ $ is order preserving (that is $g \in G^+$ entails $\hat g \geq 0$ as a function). If $G$ is unperforated, then $\hat g > \delta \pmb 1$ for some $\delta > 0$ entails that $g$ is an order unit [GH1] (the converse is trivial), and in fact it is sufficient that $\tau(g) > 0$ for all pure $\tau$. When $G$ is a dimension group, $S(G,u)$ is a Choquet simplex [G, 10.17].  Then the representation theorem for simple dimension groups asserts that if $G$ is a noncyclic simple dimension group, then the image of $G$, $\hat G$, is dense in $\text{Aff\,}S(G,u)$ and $G^+ \setminus \brcs{0}$ consists of $g$ \st $\hat g $ is bounded below away from zero. Moreover, if $G$ is any torsion free abelian group and there is a group homomorphism $\Arrow f; G.\text{Aff\,}K $ where $K$ is a Choquet simplex, and $f(G)$ is dense in $\text{Aff\,}K$, then equipped with the partial ordering $G^+ := \brcs{0} \cup f^{-1}((\text{Aff\,}K)^{++})$ (equivalently, nonzero $g$ belongs to $G^+$ if and only if $f(g)(k) > 0$ for all $k \in K$), $G$ is a simple dimension group, and if $f(u) = \pmb 1$ for some $u$ in $G$, then $S(G,u)$ is affinely homeomorphic to $K$.

For non-simple dimension groups, there are partial characterizations, but these are not usually so complete. Moreover, in the non-minimal case, $G(X,T)$ need not be partially ordered, but even if it is, it need not satisfy the Riesz decomposition property (for example, all nontrivial shifts of finite type), and so need not be a dimension group.

We noted that there is a direct translation between (pure) invariant probability measures on $(X,T)$ and (pure) normalized traces on $(G(X,T), [\chi_X])$. For this article, we have to translate goodness and related properties of invariant measures to properties of traces on $G(X,T)$.

Let  $\mu$ be a probability measure on $X$. Define $S(\mu) = \Set{ \mu(U)}{U \text{ clopen in } X}$. So $S(\mu)$ is a countable subset of the reals. Akin [Ak2] defined $\mu$ to be {\it good\/} if for all $\alpha \in S(\mu)$  and $V$ a clopen in $X$ \st $\mu(V) > \alpha$, there exists a clopen set $U \subseteq V$ \st $\mu(U) = \alpha$. We can of course translate this directly to traces on $G(X,T)$ (assuming $\mu$ is $T$-invariant), but the direct translation is awkward to use, and instead, we obtain an equivalent formulation which is much easier to deal with. That is part of the reason for the battle of the definitions (section 1).

For a partially ordered abelian group with order unit $(G,u)$ and $b$ in $G^+$, define the {\it interval generated by $b$\/} to be $[0,b]:= \Set{a \in G^+}{a \leq b}$. The naive direct translation of goodness for a trace $\tau$ of a dimension group with order unit $(G,u)$ would be, for all $b \in [0,u]$, $\tau([0,b]) = \tau([0,u]) \cap [0,\tau(b)]$ (the last is the usual interval in $\R$).  Our first order of business is to show this is equivalent to   more usable criteria, for example, for all $b \in G^+  $,
$\tau([0,b]) = \tau(G) \cap [0,\tau(b)]$, and it is equivalent to the condition replacing $\tau(G)$ by $\tau(G^+)$.  Necessary and sufficient that  $\tau$ be good is that there exist a sequence of elements $b_i \in G^+$ \st $\hat b_i \to 0$ (uniformly), $\tau(b_i) \neq 0$ and $\tau([0,b_i]) = \tau(G^+) \cap [0,\tau(b_i)]$ (Corollary \oneten).

With $u = [\chi_X]$, and $\mu$ translating to the trace $\tau$, $S(\mu)$ translates to $\tau([0,u])$.
These definitions, however, while useful, are still not easy to work with. There is a complete characterization of traces that are good, specifically, the image of $\ker \tau$ is dense (in the affine representation) in the space of affine functions that vanish at $\tau$ (a well known object in Choquet theory). It is easy to use this characterization to construct examples and non-examples.

In addition to numerous classes of examples, we also deal with dimension group extensions  (by order ideals); these correspond to closed invariant subsets of $X$; so we can ask, when goodness is preserved by the natural extension of $\mu$? This is answered in Proposition \foutwo, extending a result from stationary systems. 

If $G$ and $H$ are dimension groups, there is a natural construction of a dimension group $G\otimes H$, and pure traces of the latter are all of the form $\sigma \otimes \tau$ where $\sigma$ and $\tau$ are pure traces of $G$ and $H$ respectively. If $\sigma$ and $\tau$ are good, then so is $\sigma\otimes \tau$, but the latter can be good without either being good, and we investigate classes of examples, leading to a result with an algebraic geometrical flavour, Proposition \fivten. 

In section 7, we consider goodness en masse, that is, if $U$ is a collection of traces \st whenever $a$ and $b$ are in $G^+$ and $\tau(b)-\tau(a)$ is bounded below away from zero as $\tau$ varies over $U$. In contrast to the situation with a single trace, a characterization of good sets of traces is more complicated, having to deal with   special subsimplices of the trace space.

Refinability, a weakening of goodness, is discussed. We obtain a sufficient condition for a trace to be refinable (without being good) in terms of goodness en masse of a set of traces related to the trace. We use this to show that if $(R,1)$ is a partially ordered ring with $1$ as order unit that is also a simple dimension group (a lot of examples of this type are considered in section 5),  then all pure traces are refinable. This is somewhat surprising, since it is relatively easy to construct simple dimension groups of rational dimension three with exactly pure three traces (hence the trace space is a Bauer simplex), not all of which are refinable. 

In examining refinability, we are lead back to goodness en masse, which in turn leads to a definition of goodness for subsets of Choquet simplices that arise as intersections with closed flats. Very few of these turn out to be good, but among the good ones are faces and singleton sets (not consisting of an extreme point), and their coproducts in the category of Choquet simplices.

In   Appendix A, we outline a construction of a minimal $\Z$ action that corresponds to the tensor product of dimension groups. In Appendix B, motivated by a result in section 7 that the difference between refinability and goodness sometimes reduces to whether the corresponding ordered quotient group (by convex---but not directed---subgroups) is unperforated, we give some necessary and some sufficient geometrical conditions that the quotient be unperforated.

References for Choquet theory include Alfsen [Al], Asimow \& Ellis [AE], and various sections of Goodearl [G, sections 5, 8, 9, 11]. The latter is also good for dimension groups (and other partially ordered groups). A brief introduction to dimension groups is given in [E].

\SecT 1 The definitions shoot it out

\noindent Goodness and related notions deal with lifting properties, both in the
original context of measures on Cantor sets, and traces on dimension groups. In
order to implement the translation, we first note a property of minimal
homeomorphisms, observed in [Pu, Theorem 1.1]. Let $(G,u)$ be the simple dimension
group corresponding to $(X,T)$ where $T$ is a minimal self-homeomorphism
of the Cantor set $X$. Then $G$ can be obtained as $C(X,\Z)/(1-T)C(X,Z)$
with the quotient ordering, and $u$ is the image of the constant function
$\pmb 1$. Putnam showed that if $0 \leq a \leq u$, then there exists a
clopen set $U$ \st the image of the indicator function $\chi_U$ is $a$. The application of flow equivalence to this result, using the fact that every
nonzero element of $G$ is an order unit, yields:
\item{(1)} Suppose $\Arrow B;X.\Z$ is a nonnegative continuous function whose
image in $G$ is denoted $b$. If $a$ is an element of $G$ \st $ 0 \leq a
\leq b$, then there exists a continuous function $\Arrow A;X.\Z$ \st $ 0
\leq A \leq B$ (as functions) and $A \mapsto a$.
 
\noindent This permits the translation of the notion of good (and later on, goodness
en masse) from measures invariant under a minimal homeomorphism to simple
dimension groups.
 
In the non-simple case, the situation is somewhat more complicated. Again,
let $(X,T)$ be a Cantor set with homeomorphism, this time not necessarily
minimal, and let $(G,u)$ be defined as in the minimal case. Kim, Roush, and S
Williams [KRW, Proposition 3.1] showed that if $(X,T)$ is an irreducible shift of finite
type, then single projections lift, as in the minimal case---that is, if
$0 \leq a \leq u$, then there exists clopen $U$ \st $\chi_U \mapsto a$. By
applying inverse limits and flow equivalence, [BH3] showed that this
extends, with  provisos, to all $(X,T)$:
\item{(2)} Suppose $\Arrow B;X.\Z$ is a nonnegative continuous function whose
image in $G$, $b$,  {\it  is an order unit of $G$}. If $a$
is an element of $G$ \st $ 0 \leq a \leq b$, then there exists a
continuous function $\Arrow A;X.\Z$ \st $ 0 \leq A \leq B$ (as functions)
and $A \mapsto a$.
 
\noindent A minor proviso is that the pre-ordering on the quotient need not be a
partial ordering, although it is in most cases of interest; the major
proviso restricts the larger element to be an order unit. This is enough
however, for the notion of order unit good to translate directly in the
non-minimal case.
 
We are restricting ourselves throughout to the case(s) that the  $G$
obtained from $(X,T)$ is a dimension group (or more generally, we are
considering goodness and related notions for dimension groups). If $(X,T)$
is a primitive shift of finite type, the corresponding $G$ is not a
dimension group, and in fact, the full lifting result (that is, without
$b$ being an order unit) does not hold. However, it is plausible that when
$G = C(X,\Z)/(1-T)C(X,\Z)$ is a dimension group, then the full lifting
result does hold. If this were so, the translation of goodness from traces
on dimension groups to invariant measures of $(X,T)$ would be complete.
 
Since most of our attention is devoted to simple dimension groups---all of which
come out of minimal homeomorphisms---we don't regard the apparently
incomplete translation in the non-minimal case to be a serious problem,
especially since order unit goodness does translate in general.
 
Let $G$ be a dimension group, and let $\tau$ be a pure trace thereon. We say $\tau$ is {\it good\/} if for all $a,b \in G^+\setminus\brcs{0}$ \st $\tau(a)  < \tau (b)$, there exists $c$ in $G^+$ \st $c \leq b$ and $\tau (c) = \tau (a)$. 

Let $b$ be a positive element of $G$. Denote by $[0,b]$, the {\it interval generated by $b$,} $\Set{c \in G^+}{c \leq b}$. Then an alternative formulation of the definition of a good trace is simply that for all $b$ in $G^+$, $\tau([0,b]) = \tau(G^+) \cap [0,\tau(b)]$ (the second interval is the usual one in $\R$). This implies a stronger property, with the plus superscript deleted (Lemma \oneone\  below). The direct translation of the definition of good measure would be, for  fixed order unit $u$, and any $b \leq u$, that  $\tau([0,b]) = \tau([0,u]) \cap [0,\tau(b)]$. We show in Lemmas \oneone\ and \onethr\ that this is equivalent to our definition of good. 

Let $\tau$ be a trace of $G$. We say an element $b$ of $G^+$ is {\it weakly $\tau$-good\/} or {\it weakly good \wrt $\tau$\/} if for every $a$ in $G^+$ \st $0 < \tau(a) < \tau(b)$, there exists $c$ in $[0,b]$ \st $\tau(c) = \tau(a)$. The element $b$ in $G^+$ is {\it $\tau$-good\/} or {\it good  \wrt $\tau$} if the $a$ can vary over $
G$, that is, $\tau([0,b]) = \tau(G) \cap [0,\tau(b)]$. Obviously, $0$ is always $\tau$-good, and $\tau$ is good if and only if every element of $G^+$ is weakly $\tau$-good. 
The next result says that this is equivalent to all elements of the positive cone being  good \wrt $\tau$. 

\Lem Lemma \oneone. 
(a) If $\tau$ is a normalized trace of $(G,u)$ \st for all $b$ in $[0,u]$ and $a \in [0,u]$ \st $\tau(a) < \tau(b)$, there exists $c \in [0,b]$ \st $\tau(c) = \tau(a)$, then $\tau$ is good. 
\item{(b)} If $\tau$ is a good  trace of $G$, then $\tau([0,b]) = \tau(G) \cap [0,\tau(b)]$ for all $b$ in $G^+$; that is, every element of $G^+$ is weakly good implies every element of $G^+$ is good \wrt $\tau$.

\Pf (a) Select $a \in G^+ $ \st $\tau(a) < \tau(b)$. There exists an integer $n$ \st $a \leq nu$ (from the definition of order unit). By Riesz decomposition, there exist $a_i \in [0,u]$ ($i = 1,2, \dots, n$) \st $a = \sum a_i$. Since $a_n \leq a$, we have $\tau(a_n) < \tau(b)$, so there exists $c_n \in [0,b]$ with $\tau(c_n) = \tau(a_n)$. Then $\tau\(\sum_{i< n} a_i\) < \tau(b- c_n)$, and we proceed by induction on $n$. 

\noindent (b) Given $a$ in $G$ and $b$ in $G^+$ \st $0 < \tau(a) < \tau (b)$ (if either equality holds, the result is immediate), there exist positive elements $c$ and $d$ \st $a = c-d$ (since $G$ is directed, part of the definition of partially ordered abelian group). Since $\tau(c) -\tau(d)= \tau(a) > 0$, we can apply weak goodness to  $d$ and $c$, yielding an element $e$ in $[0,c]$ \st $\tau(e) = \tau(d)$. Then $c-e$ is a positive element of $G$ and $\tau(c-e) = \tau(a)$, so we can apply weak goodness to the pair $c-e$ and $b$.\qed

If $\tau$ is a  trace on $(G,u)$, we say that $\tau$ is {\it reasonable\/} if for all $b$ in $G^+$, $\tau([0,b]) $ is dense in $\tau(G) \cap [0,\tau(b)]$. Most traces that we are likely to encounter are reasonable, but some are not: let $(H,v)$ be any nonzero dimension group with order unit, take $G= \Z \oplus H$ with the coordinatewise ordering, and let $\sigma$ be a normalized trace of $(H,v)$. Then $\tau :(m,h)\mapsto m + \sigma (h)$ is not reasonable if $|\sigma([0,v])| \geq 3$ (take $b = (1,0)$). If $\sigma$ has dense range in $\R$, then $\tau$ does as well, but $\tau([0,b])$ is discrete. More generally, if $G$ has two simple quotients (by order ideals) at least one of which is cyclic, then $G$ admits a trace that is  not reasonable.  By Lemma~\oneone, a good trace is reasonable. 

A trace is {\it discrete\/} [G] if its image is cyclic. A dimension group $(G,u)$ with order unit has no discrete (pure) traces if and only if its image in $\text{Aff\,}S(G,u)$ is dense, \wrt the supremum norm [GH1; Theorem 4.8], and then is called {\it approximately divisible.}

We also recall the {\it purity criterion\/} for traces on dimension groups with order unit, $(G,u)$, given in [GH1; Theorem 3.1]:
\itemitem{} the trace $\tau$ is pure if and only if for all $\epsilon > 0$ and all $a$ and $b$ in $G^+$, there exists $c \in [0,a] \cap [0,b]$ (that is, $c \leq a$; $c \leq b$; and $c \in G^+$) \st $\tau(c) > \min \brcs{\tau(a), \tau(b)}-\epsilon$.

\Lem Lemma \onetwo. Let $(G,u)$ be a dimension group with order unit. 
\item{(a)} Suppose  all order ideals admitting an order unit are approximately divisible. Then all traces of all order ideals with order unit (including $(G,u)$ itself) are reasonable and have dense range in $\R$.
{\par}
\item{(b)} If $\tau$ is a pure trace of $G$, then $\tau $ is reasonable.

\Pf (a) It suffices to do this when $b$ is an order unit, since we can take
the order ideal generated by $b$. Then $G$ has dense image in
$\text{Aff\,} S(G,b)$ (so $\hat b = \pmb 1$, the constant function). Hence
for each $n$, there exists $b_n$ \st $n^{-1} \pmb 1 > \hat b_n > (n+1)^{-1}\pmb
1$ (approximate the constant affine function with value $(1/2) (n^{-1} +
(n+1)^{-1})$ with an error of less than $1/2n(n+1)$). Then $kb_n$ belongs to
$[0,b]$ for $1 \leq k \leq n$, and the value of $kb_n$ at every normalized
trace belongs to the interval $(k/(n+1), k/n)$, and thus the set of values
of $\brcs{kb_n}_{1 \leq k \leq n; n = 1,2, \dots}$ is dense in $[0,1] = [0,
\tau(b)]$.

\noindent (b)  For fixed $b$ in $G^+$, suppose $c$ is an element of $G$ \st $ 0 < \tau(c) < \tau(b)$. Then we can write $c = g-h$ where $g$ and $h$ are positive elements. Apply the purity criterion to the pair $g,h$; given $\epsilon$, we may find positive $f \leq g,h$ \st $\tau(f) > \min\brcs{\tau(g), \tau(h)}- \epsilon = \tau(h)-\epsilon$. Then $c = (g-f)-(h-f)$, and $\tau(h-f) < \epsilon$. If $\epsilon < \tau(b-c)/2$, then $\tau(g-f) < \tau(b) $. Since $g-f \geq 0$, we may apply the purity criterion again, obtaining $0 \leq e \leq g-f, b$ \st $\tau(e) > \tau(g-f)-\epsilon$, and thus $|\tau (e)-\tau(c)| < |\tau(e-g + f)| + |\tau(g-f -c)| < 2\epsilon$, and $e \in [0,b]$.
\qed

Suppose $H$ is a subset of $G^+$ \st every element of $G^+$ is a sum of elements of $H$; we call $H$ a {\it generating subset\/} of $G$. Examples include $H = [0,u]$ if $u$ is an order unit of $G$, and $H = \Set{h_1 \otimes h_2}{h_i \in H_i}$ if $G = G_1 \otimes G_2$ (the ordered tensor product; this is a dimension group if both $G_i$ are, and will be discussed later), and  $H_i$ is a generating subset of $G_i$.

\comment
\Lem Proposition \onetwo. [GH1] (Purity criterion) Let $(G,u)$ be a dimension group with order unit, and let $\tau$ be a trace on $G$. Then $\tau$ is pure if and only if for all $g,h$ in $G^+$ and all $\epsilon > 0$, there exists $e$ in $G^+$ \st $e \leq g,h$ and $\tau(e) \geq \min\brcs{\tau(g),\tau (h)}-\epsilon$.
\endcomment

\Lem Lemma \onethr. Let $G$ be a dimension group with a reasonable trace $\tau$ with $\tau(G)$ dense in $\R$. Let $v_1$ and $v_2$ be nonzero elements of $G^+$. If both are  (weakly) $\tau$-good then $v_1 + v_2$ is as well. 

\Pf 
\comment
We can deal with both cases simultaneously. Suppose $\tau(g) < \tau(v_1 + v_2)$. By reasonableness, for $\epsilon > 0$, we may find $w_i \in [0,v_i]$ \st $\tau(w_i) \in (\tau(v_i)-\epsilon, \tau(v_i)-2\epsilon)$ and so that $\tau(g) < \tau(w_1 + w_2)$ (restrict to $\epsilon < \tau(v_1 + v_2-g)$). If $\tau(g) \leq \max\brcs{\tau(v_i)}$, the result is immediate. Since $\tau([0,w_1]) + \tau([w_2]) $ is dense in $[0,\tau(w_1)] + [0,\tau(w_2)] = [0,\tau(w_1 + w_2)]$,  we may find $\alpha_i \in \tau([0,w_i])$ \st $|\alpha_1 + \alpha_2-\tau(g)|\leq \epsilon/4$
\endcomment
First suppose that both $v_i$ are  $\tau$-good; we wish to show $v_1 + v_2$ is as well. Select $f$ in $G$ \st $0 < \tau(f) < \tau(v_1 + v_2)$ (if either inequality is replaced by equality, there is nothing to do). Select $\epsilon < \tau(v_1 + v_2-f)/2$. If $\tau(f) \leq \max\brcs{\tau(v_i)}$, there is an obvious reduction, so we may assume that $\tau(f) > \max\brcs{\tau(v_i)}$. Since $\tau(f) \in  [0,\tau(v_1) + \tau(v_2)] =   [0,\tau(v_1)] +  [0, \tau(v_2)]$, there exist $s_i$ in $[0,\tau(v_i)]$ \st $\tau(f) = s_1 + s_2$. Then the larger one of them, say $s_1$ (after relabelling) must exceed $\tau(f)/2$, and the other one must be nonzero. By reasonableness of $\tau$, given $\epsilon > 0$, there exists $h $ in $[0,v_1]$ \st $|s_1-\tau(h) |< \epsilon$, where we choose $\epsilon < \min\brcs{\tau(v_1 + v_2-f)/2, \tau(f)/4, s_2}$. Then $|\tau(f-h)-s_2| < \epsilon$, so that $\tau(f -h) > 0$, and since $s_2 \leq s_1$, we have $\tau(f -h) < \tau(v_1)$. By goodness applied respectively to $h$ and $f -h$, there exist $h_i \in [0,v_i]$ \st $\tau(h) = \tau(h_1)$ and $\tau(f-h) = \tau(h_2)$. Hence $\tau(f) = \tau(h_1 + h_2) \in [0,v_1 + v_2]$. 
\comment
Now assume that $\tau$ is pure and $v_i$ are weakly good \wrt $\tau$. Select $f$ in $G^+$ \st $0 < \tau(f) < \tau(v_1+v_2)$. Select $\epsilon < \tau(v_1 + v_2-f)/2$. Again, we may assume $\tau(f) > \max\brcs{\tau(v_i)}$. By the purity criterion, there exists $0 \leq f_1 \leq f,v_1$ \st $\tau(f_1) > \min\brcs{\tau(f),\tau(v_1)}-\epsilon = \tau(v_1)-\epsilon$. In particular, $f':= f-f_1$ belongs to $G^+$, and again applying the purity criterion, there exists $0 \leq f_2 \leq f',v_2$ \st $\tau(f_2) > \min\brcs{\tau(f'), \tau(v_2)}-\epsilon$. 

If $\tau(f') \leq \tau(v_2)$, we may apply weak goodness of $v_2$ and obtain $f'' \in [0,v_2]$ \st $\tau(f'') = \tau(f')$. Then $f'' + f_1 \in [0,v_1 + v_2]$, and $\tau(f'' + f_1) = \tau(f' + f_1) = \tau(f)$ as desired. 

We conclude by showing $\tau(f') > \tau(v_2)$ is impossible. If $\tau(f') > \tau(v_2)$, then 
$$
 0 \leq \tau(f-f_1-f'') < \tau(f)-\tau(v_1)-\tau (v_2) + 2\epsilon < 0,
$$
a contradiction. Hence $v_1 + v_2$ is weakly $\tau$-good.
\endcomment

Now suppose the $v_i$ are $\tau$-weakly good. Suppose $a \in G^+$ and $0 < \tau(a) < \tau(v_1 + v_2)$. Set $\alpha_i = \tau(v_i) \tau(a)/\tau(v_1 + v_2)$, so that $0 < \alpha_i < \tau(v_i)$ and $\tau(a) = \alpha_1 + \alpha_2$. Choose $\epsilon < \min\brcs{\tau(v_i)-\alpha_i}$. By reasonableness of $\tau$, we may find $d_1 \in [0,a]$ \st $\alpha_1-\epsilon < \tau(d_1) < \alpha_1 < \tau(h_1)$. Then $\tau(a-d_1) < \alpha_2 + \epsilon < \tau(v_2)$. By weak goodness of the $v_i$,  there exist $c_1 \in [0,v_1]$ \st $\tau(c_1) = \tau(d_1)$ and $c_2 \in [0,v_2]$ \st $\tau(a-d_1) = \tau(c_2)$. Then $c_1 + c_2 \in [0,v_1] + [0,v_2] \subseteq [0,v_1+v_2]$ and $\tau(c_1 + c_2) = \tau(a)$. 
\qed

\Lem Corollary \onetwe. Suppose that $G$ is a dimension group, and $\tau$ is a reasonable trace with dense range in $\R$. If $G^+$ admits a generating set consisting of weakly $\tau$-good elements, then $\tau$ is good.

This will be superseded (in the simple case) by Corollary \oneten,  where it is shown that sufficient for goodness is that  there exist a sequence of weakly $\tau$-good elements $v_i$ \st $\hat v_i \to 0$ (that is, \wrt the supremum norm on $S(G,u)$). 

Suppose that $b$ is an element of the positive cone of the dimension group $G$. Recall that  $[0,v]$ denotes $\Set{g \in G^+}{g \leq v}$. Then $\tau$ is {\it group-like\/} \wrt $v$ if whenever $a$ belongs to $G$ and $0 \leq \tau (a) \leq \tau (v)$, there exists $c$ in $[0,v]$ \st $\tau(c) = \tau (a)$, that is, $\tau([0,v])
= \tau(G) \cap [0,\tau(v)]$; we say $\tau$ is {\it weakly group-like \wrt $v$} if $a$ is restricted to $G^+$---that is, $\tau([0,v])
= \tau(G^+) \cap [0,\tau(v)]$.

\noindent{\it Examples.} There exist simple dimension groups (which are free abelian groups of rank 3 or 4) with two pure traces \st neither one is good, or exactly one is good, or both are good.  There also exist simple dimension groups with two pure traces \st neither trace  is good but become good on tensoring with the rationals, and also that don't become good on so tensoring. Similarly, there exist simple dimension groups for which being group-like depends on the choice of order unit.

Many of  the examples are of the form $G \subset \R^2$, where $G$ is a dense subgroup of rank three, equipped with the strict ordering inherited from $\R^2$, that is, $G^+ \setminus\brcs{0} = (\R^2)^{++} \cap G$. These are dimension groups, and the two projection maps (onto the first and second coordinates, respectively), $\tau_1$, $\tau_2$, are the only pure traces. If $G$ is the abelian group spanned by the three vectors, $v_1, v_2, v_3$ of $\R^2$, then $G$ is dense if and only if the set of three $2\times 2$ determinants, $\brcs{v_1 \wedge v_2, v_1 \wedge v_3,v_2 \wedge v_3}$ is linearly independent over the rationals.

\Lem Example \exone. Let $\brcs{1,\alpha, \beta}$ be a set of real numbers that is linearly independent over the rationals and \st $ 0 < \alpha < 1 < \beta$;  set $G = \left\langle(1,0), (0,1),(\alpha,\beta) \right\rangle$. Linear independence of $\brcs{1,\alpha, \beta}$ guarantees $G$ is dense in $\R^2$, and thus is a simple dimension group \wrt the strict order. Let $\tau = \tau_2$ be the projection onto the second coordinate.

\noindent (a) {\it $\tau$ is not good.} Set $a = (1,1)$ and $b = (\alpha,\beta)$; we show the existence of $c$ in $G^+$ \st $\tau(a) = \tau(c)$ and $c \leq b$ is impossible. As a group homomorphism $\Arrow \tau; G. \R$, its kernel is spanned by $(1,0)$. Hence $c -a = (m,0)$ for some integer $m$, that is, $c = (m+1,1)$; since $c$ is in $G^+\setminus {0}$, we must have $m \geq 0$, whence $\tau_1 (c) \geq 1$. But $\tau_1(b) =  \alpha < 1 \leq \tau_1 (c)$, contradicting $c \leq b$.

\noindent (b) {\it The extension of $\tau$ to $G \otimes \Q$ is good.} Write $a = (a_1,a_2)$ and $b = (b_1,b_2)$ in $(G\otimes \Q)^+$ with  $0 < a_2 < b_2$ (so $\tau(a) < \tau(b)$). Since $a$ and $b$ are nonzero positive elements, $a_1, b_1 > 0$. Since the positive rationals are dense in the positive reals, there exists a positive rational number $q$ \st $0 < a_1-q < b_1$. Set $c = a-q\cdot(1,0)= (a_1-q, a_2)$; then $c$ belongs to $G \otimes \Q$, is strictly positive (hence is in the positive cone), $\tau(c) = a_2 = \tau(a)$, and $c $ is strictly less than $b$, so $b-c \in (G\otimes\Q)^+$, verifying goodness.

\noindent (c) {\it With respect to the order unit $u = (1,1)$, $\tau $ is group-like.}
Suppose $a = (r,s)$ is an element of $G^+$ \st $\tau (a) \leq 1$, that is, $s\leq 1$. If $s = 1$, set $b = u$ and we are done. We may thus suppose $0 < s < 1$. Decompose $r = \flo{r} + \brcs{r}$ into its integer and fractional parts. If $\brcs{r}$ is not zero (that is, $r$ is not an integer), then set $b = a-\flo{r}(1,0) = (\brcs{r}, s)$, so that $b$ belongs to $[0,u]$ and $\tau (b) = s$.

Otherwise, $r$ is an integer. Then we can write $a = (r,s) = m(1,0) + n(0,1) + p\cdot (\alpha,\beta)$ where $m,n,p$ are integers. Evaluating at the first coordinate yields $r = m + p\alpha$; since $\brcs{1, \alpha}$ is linearly independent over $\Q$ and $r$ is an integer, we deduce $r = m$ and $p = 0$. Hence $a = (m,n)$ which is impossible, since the second coordinate of $a$ is strictly between $0$ and $1$.

\noindent (d) {\it With respect to the order unit $u = (\alpha,\beta)$, $\tau$ is not group-like.} Set $a = (1,1)$; thus  $\tau(a) = 1 < \beta = \tau(u)$. If $b$ were an element of $[0,u]$ \st $\tau (b) = 1$, then $a-b \in \ker \tau = \Z(1,0)$; then $b = (m,1)$ for some integer $m$; as $b\geq 0$ (but $b \neq 0$), we must have $m \geq 1$. On the other hand, $b \leq (\alpha,\beta)$ implies $m \leq \alpha < 1$, a contradiction.

\Lem Example \extwo. Let $\alpha_i, \beta_i$ ($i=1,2$) be four positive real numbers \st each set  $\brcs{1,\beta_1, \beta_2}$ and $\brcs{\beta_1- \alpha_1, \beta_2-\alpha_2, \alpha_1 \beta_2-\beta_2 \alpha_1}$ is linearly independent over the rationals (this will occur, for example, if $\brcs{1,\alpha_i,\beta_i}$ is algebraically independent over the rationals). Assume in addition that $\alpha_i < \beta_i$  
and $\alpha_1 + \alpha_2 < 1 < \beta_1 + \beta_2$.  Set $G = \langle(1,1), (\alpha_1, \beta_1), (\alpha_2, \beta_2) \rangle$

Linear independence of the second set guarantees $G$ is dense in $\R^2$, so equipped with the strict ordering, $G$ is a simple dimension group. The important property here is that $\ker \tau$ is trivial, as follows from linear independence of $\brcs{1,\beta_1,\beta_2}$.

\noindent{\it For no order unit $u$ is $\tau = \tau_2$ group-like.} This is a special case of the following easy result, which is preliminary to our later characterization of good pure traces on simple dimension groups. The group of infinitesimals is denoted $\Inf (G)$, and is equal to the intersection of the kernels of all traces (when $G$ has an order unit).

\Lem Proposition \onefou. Suppose that $(G,u)$ is a  dimension group, and $\tau $ is a  trace of $G$ \st $\ker \tau = \Inf(G)$. If $G$ has more than one trace, then no positive element $b$ \st $\tau(b) \neq 0$ is weakly $\tau$-good (in particular, every element  is not $\tau$-good).

\Pf 
\comment
Let $\tau'$ be some other pure trace, and let $u$ be an order unit. We may assume that both $\tau$ and $\tau'$ are normalized \wrt $u$.  Represent $G$ as a group of affine continuous functions on its (normalized) trace space, $S(G,u)$. Since $G$ is a simple dimension group, the image of the positive cone is dense in the positive cone of $\text{Aff\,}S(G,u)$. We may thus find an element $a$ in $G^+$ \st $\hat a(\tau) < 1/2$ but $\hat a(\tau') > 2$; that is, $\tau(a) < 1/2$ and $\tau'(a) > 2$. Select $b$ in the positive cone \st $\tau(b) > 1/2$ but $\tau'(b) < 2$ (e.g., $b = u$). If $c$ is any element of $G$ \st $\tau(a) = \tau (c)$, then from $\ker \tau = 0$, we have $a = c$. Since $\tau'(b)$ is thus greater than $2$, $c \leq b$ is impossible.\qed
\endcomment
Without loss of generality, $\Inf (G)  = 0$ and $\tau(u) = 1$. Suppose there exists another normalized trace unequal to $\tau$, $\sigma$. If there exist $a$, $b$ in $G^+$ \st $\tau(a) < \tau(b)$ but $\sigma(a) > \sigma (b)$, and $0 \leq c \leq b$ \st $\tau(c) = \tau (a)$, then    $a = c$, and thus $\sigma(c)> \sigma(b)$, a contradiction. Hence if $b$ is weakly $\tau$-good, then $\tau (a) < \tau (b)$ implies $\sigma (a) \leq \sigma (b)$ for all other traces $\sigma$. Multiply $a$ and $b$ by a big enough integer $K$ \st   $\tau (Ka)<  \tau((K-1)b)$; then $\sigma (Ka) \leq  \sigma (K-1)b$ for all traces unequal to $\tau$, and thus $\tau(a) < \tau(b)$ entails $\sigma(a) < \sigma(b)$ for all traces, and thus $b-a$ is an order unit. However, since the elements of $G^+$ separate  traces, for any pair of distinct normalized traces, there exist $a$ in $G^+$ \st $\tau(a) < \tau(b)$ and $\sigma (a) > \sigma(b)$. \qed

\Lem Lemma \onefiv. If $G$ is a simple dimension group,  $\tau$ is pure, $G$ has more than one trace, and $\ker \tau/\Inf(G)$ is cyclic, then $\tau$ fails to be good.

\Pf Since $G/\Inf(G)$ is a simple dimension group with the same traces and values,  we may assume that $\Inf (G) = 0$. If $\ker \tau = \brcs{0}$, we are in the case of Proposition\,\onefou\ above (good implies group-like \wrt every order unit). Otherwise, we note that $\ker \tau$ is obviously torsion-free, so must be isomorphic to  $\Z$, i.e., there exists $g$ in $G$ \st $\ker \tau = \Z g$. Since $\Inf (G)= 0$, there exists a pure trace $\tau_1$ \st $\tau_1 (g) \neq 0$; by replacing $g$ by $-g$ if necessary, we may assume $\tau_1(g) > 0$. Select any nonzero positive element  $u$ \st $\tau_1(u) < \tau_1(g)$ (such exist from elementary results); then $u$ is an order unit for $G$ from simplicity. Renormalize the traces at $u$ (so $\tau_1 (g) > \tau_1 (u) = \tau(u) = 1$), and represent $G$ in the group of affine functions on its normalized trace space.

By density of the positive cone in the positive cone of the latter, there exists $a$ in $G^+$ \st $\tau(a) < 1/2$ and $1 < \tau_1 (a) < \tau_1(g)$. If $b$ were any element of $G^+ $ \st $\tau (a) = \tau(b)$, then $b = a + mg$ for some integer $m$. Applying $\tau_1$, we have $\tau_1(b) = \tau_1 (a) + m\tau_1 (g)$. If $m$ were negative, we would obtain $\tau_1(b) < 0$, a contradiction. If $m \geq 0$, then $\tau_1(b) \geq \tau_1(a) > 1$; hence $b$
 cannot be in $[0,u]$. Hence $\tau$ is not group-like \wrt some order unit, and therefore $\tau$ is not good.
\qed

Both of these easy results will be superseded, but they are worth mentioning, as the criteria are usually easy to verify. 
\comment
\Lem Conjecture. If $G$ is a simple dimension group and $\tau$ is a pure trace, then $\tau$ is good iff $\tau$ is group-like for all order units $u$ of $G$.

Let $u$ be an order unit of $G$, and let $K= S(G,u)$ be the usual normalized trace space, a Choquet simplex. Form $\text{Aff\,} (K)$, equipped with the usual (not the strict) order. If $\tau$ is a pure, normalized trace of $(G,u)$, then $\brcs{\tau}$ is a face of $K$. Let $I(\tau)$ denote the set of $\alpha$ in $\text{Aff\,} (K)$  \st $\alpha (\tau) = 0$. Then $I(\tau)$ is a maximal order ideal of $\text{Aff\,} (K)$, and every maximal order ideal is of this form (for a pure trace $\tau$). 

\Lem Conjecture. If $G$ is a simple dimension group and $\tau $ is a pure trace, normalized \wrt an order unit $u$, then $\tau$ is good if and only if the image of $\ker \tau$ in $\text{Aff\,} (K)$ (under the map, $g \mapsto \hat g$, where $\hat g(k) = k(g)$) is (norm)-dense in $I(\tau)$. 

This is true if $K$ is finite-dimensional (that is, $G$ has only finite many pure traces), and when I have time after doing my NSERC proposal, I will write it up. This extends drastically the earlier proposition and lemma.

\noindent {\bf Right definition (sort of)}. 

\noindent If $v$ is an element of $G^+$,  recall $[0,v]:= \Set{g \in G}{0 \leq g \leq v}$.
Let $G$ be a dimension group. Let $\tau$ be a pure (finite) trace on $G$ and let $u$ be an order unit. Then {\it $\tau$ is weakly group-like \wrt $u$\/} if for all $a$ and $a'$ in $[0,u]$ \st $\tau (a+a') < \tau(u)$, there exists $c$ in $[0,u]$ \st $\tau(c) = \tau (a+a')$.

If we replace $\brcs{a,a'}$ by an arbitrary finite set, $\brcs{a_i} \subset [0,u]$ \st $\sum \tau(a_i) < \tau (u)$, then $\tau$ being weakly group-like implies there exists $c$ in $[0,u]$ \st $\tau (c) = \sum \tau(a_i)$ by a simple induction argument.

Recall the definition of a trace being good; a pure trace $\tau$ on $G$ is good if for all $a$ and $b$ in $G^+$ \st $\tau (a) < \tau(b)$, there exists $0 \leq c \leq b$ \st $\tau(c) = \tau (a)$. The first two results are practically tautological. The converse requires simplicity.

\Lem Lemma. If $\tau$ is a good pure trace on $G$, then $\tau$ is weakly group-like \wrt every order unit $u$ of $G$.

\Pf Choose order unit $u$ and set $b = u$; if $\tau (a+a') < \tau (b)$, there exists $ 0 \leq c \leq b$ (so $c$ belongs to $[0,u])$ \st $\tau(c) = \tau(a+a')$, and we are done.\qed

\Lem Proposition. Suppose that $\tau$ is a pure trace of $G$ that is weakly group-like \wrt every order unit $u$. If $G$ is simple, then $\tau$ is good.

\Pf Since $G$ is simple, every nonzero positive element is an order unit. Given $a$ and $b$ in $G^+$ with $\tau (a) < \tau(b)$, set $u = b$. Since $G$ is a dimension group (thus satisfies Riesz decomposition), we may decompoose $a = \sum a_i$ with $a_i $ in $[0,u]$. Then $\sum \tau(a_i) < \tau (u)$, so by the iterated form of the definition of weakly group-like, there exists $c$ in $[0, u]$ \st $\tau(c) = \sum\tau(a_i)$. In particular, $0 \leq c \leq b$.
\qed
\endcomment

We say a dimension group $G$ with an order unit  is {\it approximately divisible\/} if for all pure traces $\sigma$, $\sigma (G)$ is non-discrete in $\R$. This is equivalent to the natural map $G\to \text{Aff} (S(G,u))$ having norm-dense range, and also to the image of $G^+$ being dense in the positive cone of $\text{Aff} (S(G,u)) $ [GH1; Theorem 4.8]; if $G = \text{K}_0 (A)$ for some unital AF algebra, approximate divisibility of $G$ is equivalent to $A$ having no finite dimensional representations.

Fix an order unit $u$ of $G$ and a pure trace $\tau$, which we may now assume is normalized \wrt $u$. Form the affine representation $G\to \text{Aff} (S(G,u))$, $g \mapsto \hat g$ where $\hat g (\sigma) = \sigma (g)$ for normalized traces $\sigma$. With respect to the usual ordering, $\text{Aff} (S(G,u))$ is itself a dimension group, and $\tau $ induces a map $\Arrow {\tilde \tau}; \text{Aff} (S(G,u)).\R$ with kernel denoted $\brcs{\tau}^{\perp} = \Set{f \in \Aff S(G,u)}{f(\tau) = 0}$. Since $\tau$ is pure, $\brcs{\tau}^{\perp}$ is an order ideal of $ \text{Aff} (S(G,u))$, in fact, a maximal order ideal. 

If $\tau$ is not pure, it still makes sense to consider $\brcs{\tau}^{\vdash} : = \Set{h\in \text{Aff\,}S(G,u)}{h(\tau) = 0}$ (note the slight difference in notation: when $\tau$ is pure, we use the perp symbol ${}^{\perp}$, but if $\tau$ is not, we use ${}^{\vdash}$; the former will signal that it is an order ideal in $\Aff S(G,u)$), but the latter is only  a closed convex codimension one subspace of $\text{Aff\,}S(G,u)$, not (necessarily) an order ideal. Obviously $\hat{\ }$ sends $\ker \tau$ to $\brcs{\tau}^{\vdash}$. 

If we restrict the definition of good by requiring that   $b$ be an order unit, then we call the corresponding property (of $\tau$) {\it order unit good.} Obviously if $G$ is simple and $\tau$ is order unit good, then it is good (since all nonzero elements of $G^+$ are order units).

\Lem Proposition \onesix. Let $G$ be an approximately divisible dimension group with order unit, and let $\tau$ be a 
trace.
Then  $\tau$ is order unit good if and only if the image of  $\ker \tau$ is norm dense in $\brcs{\tau}^{\vdash}$.

\comment
Suppose there exists a pure trace $\sigma$ (unequal to $\tau$), \st $\sigma (\ker \tau)$ is cyclic (possibly even zero); it follows there exists a pure trace with this property, which we now assume. Then we show that $\tau$ is not good. Fix an order unit $u$. There exists unique nonnegative real number $\alpha$ \st $\ker \sigma = \alpha \Z$ ($\alpha = 0$ is possible). If $\alpha > 0$ use  density to find two order units (elements $g$ \st that $\hat g$
is strictly positive), $a$ and $b$ in $G$ with the following properties:
$$\eqalign{
\hat b(\sigma) &< \hat a (\sigma) < \alpha\cr
\hat a(\tau ) & < \hat b(\tau). \cr
}$$
If $\tau (c) = \tau (a)$, then $a-c $ belongs to $\ker \tau$, so $\sigma (a-c) = m\alpha$ for some  integer $m$, that is, $\sigma (c) = \sigma(a)-m\alpha$. If $c \leq b$, then $\sigma (c ) \leq \sigma(b)$, so $\sigma (a) \leq \sigma(b) + m\alpha$; the  first displayed inequality forces $m \geq 1$, but this means that $\sigma(c) < 0$, contradicting $0 \leq c \leq b$.

If $\alpha = 0$, simply find order units $b$ and $a$ in $G$ \st $\hat b(\sigma) < \hat a (\sigma)$ and $\hat a( \tau) < \hat b(\tau)$. If $\tau (c) = \tau (a)$, then $\sigma(c) = \sigma (a)$, so that $0 \leq c \leq b$ (which entails $\sigma (c) \leq \sigma(b)$) is impossible.
\qed

Note that in the latter case, we barely need anything---we do not require density, simply that the traces can be distinguished by positive elements.

This might already prove the first conjecture, as this implies that for a good trace $\tau$,  we must have $\sigma (\ker \tau)$ is dense in $\R$ for all traces $\sigma$. This might be enough to show density (sufficient to deal with differences of traces, that is, $\sigma_1-\sigma_2$; if image is discrete, obtain weirdo thing; if image is always dense, then density is immediate.
\endcomment

\Pf Fix an order unit, $u$, \wrt which the trace space, $S(G,u)$, and the corresponding affine representation, $\hat{\ },$ are determined. We may assume that $\tau(u) = 1$. 

Suppose that $\tau$ is order unit good.  Select $h$ in $\brcs{\tau}^{\vdash}$. From approximate divisibility of $G$, given $\epsilon> 0$, there exists $g$ in $G$ \st $\|\hat g-h\| < \epsilon$, so that $|\tau(g)| < \epsilon$. Again by density of $\hat G$, there exists $g_0$ approximating to within $\epsilon$ the constant function with value $2\epsilon$, so that $\|\hat  g_0 \| < 3\epsilon$, and thus $g' := g_0 + g$ satisfies $\| \hat g'-h\| < 4\epsilon$ and $\tau(g') >0$. 

Again by density of $G$ in the affine representation, there exists $b$ in $G$ \st $4\epsilon <\hat b(\sigma)< 5\epsilon$ for all $\sigma$ in $S(G,u)$
(for example, a sufficiently fine approximant of the constant affine function with value $9\epsilon/2$). Automatically, $b$ is an order unit, thus positive, and $\tau(g') < \tau(b)$. By order unit goodness, there exists $0 \leq c \leq b$ \st $\tau(c) = \tau(g')$; in particular, $\| \hat c\| \leq \|\hat b\| < 5\epsilon$. Thus $y:=g'-c $ belongs to $\ker\tau$, and $\|h-\hat y \| \leq \|h-\hat g' \| + \|\hat c \| < 9 \epsilon$.

\comment
Recall that if $\phi$ is a continuous linear functional on $\text{Aff\,}S(G,u)$, we can decompose $\phi = \phi_0 + \lambda \tau$ where $\lambda$ is a real number, and $\phi_0$ is in the linear span of the complementary face to $\brcs{\tau}$. To show that the image of $\ker \tau$ is dense in $M(\tau)$, it suffices to show that for every nonzero continuous linear functional $\phi$ on $\text{Aff}(S(G,u))$ for which the corresponding $\lambda$ is zero (that is, $\phi$ is in the span of the complementary face to $\brcs{\tau}$), the range of $\phi$ on $\ker \tau$ is dense in $\R$. For such $\phi$, we may write $\phi = \sigma_1-\sigma_2$ where $\sigma_i$ are unnormalized (and probably not pure) traces, each of which is disjoint from $\tau$ (that is, for all $\epsilon> 0$, there exists $a$ in $[0,u]$ \st $\sigma_i(a) < \epsilon $ and $\tau(a) > 1-\epsilon$). In what follows, we exploit the fact that $\brcs{\tau}$ is a compact (face) in $S(G,u)$, using standard (?) separation theorems [Alfsen \& Schultz, and the later one]. 

Suppose $\phi$ restricted to the image of $\ker \tau$ is zero. Select an order unit $a$ \st $|\sigma_1(a)-\sigma_2(a)| > 1$. Select $b$ an order unit \st $\tau(b) > \tau(a)$, but $(\sigma_1+\sigma_2)(b) < 1/2$ (this is possible, since $\brcs{\tau}$ is a compact face disjoint from $\brcs{\sigma_1/\sigma_1(u), \sigma_2/\sigma_2(u)}$.  If $\tau(a) = \tau (c)$, then $(\sigma_1- \sigma_2)(c ) = (\sigma_1- \sigma_2)(a)$, so $|(\sigma_1- \sigma_2)(c)| > 1$; but $c \leq b$, so $|(\sigma_1- \sigma_2)(c) | \leq (\sigma_1+\sigma_2)(b) < 1/2$, a contradiction.

Now suppose $\phi$ restricted to the image of $\ker \tau$ has discrete range, that is, values in $\Z \alpha$ for some positive real $\alpha$. Select an order unit $b$ \st $(\sigma_1  + \sigma_2) (b) < \alpha/3$, and then try to find an order unit $a$ \st  $\slfrac{\alpha}3 < |\sigma_1 (a)-\sigma_2 (a)| < \slfrac{2\alpha}3$, and in addition $\tau(a) < \tau (b)$. If $\tau(c) = \tau (a)$, then $(\sigma_1-\sigma_2 )(c) =
(\sigma_1-\sigma_2 )(a) + m\alpha$, for some integer $m$.

The possible values of $c$ under $\sigma_1-\sigma_2$ lie in the middle thirds intervals, $m\alpha + \alpha(\slfrac13,\slfrac23)$. However, if $0 \leq c \leq b$, then $|(\sigma_1 -\sigma_2)(c)| \leq ( \sigma_1 + \sigma_2) (b)$, so $(\sigma_1 -\sigma_2)(c)$ lies in the first or third interval, that is,  $|(\sigma_1 -\sigma_2)(c) | < \slfrac{\alpha}3$, a contradiction.
\endcomment 

Now suppose that the image of $\ker \tau$ is dense in $\brcs{\tau}^{\vdash}$.  Fix an order unit $u$ (for constructing the affine representation), let $b$ be an order unit and let $a$ be an element of $G^+$ \st $0 <\tau(a) < \tau(b)$. Set $h = \tau(a) \hat b/\tau (b)$, so that $h$ is strictly positive and $(h-\hat a)(\tau) = 0$, hence $h - \hat a \in \brcs{\tau}^{\vdash}$. Select $\epsilon$ so that
 $$
0 < \epsilon <\(\inf_{\sigma \in S(G,u)}\brcs{\sigma(h)}\)\cdot \min\brcs{\frac{\tau(b)- \tau(a)}{\tau(b)}, 1 };
$$
 we may find $g$ in $\ker \tau$ \st $\|\hat g-h + \hat a\| < \epsilon$. Now consider $c = g+a$. Its value at $\tau$ is $\tau(a)$, and for every normalized trace $\sigma$, we have $|\sigma (g+a)-\sigma (h)| < \epsilon$, so $\sigma(g+a) > \sigma(h)-\epsilon > 0$; thus $g+a$ is an order unit, so in particular, $c$ is positive. Next $\sigma (b)-\sigma (g+a) =(\hat b-h)(\sigma) + (h-\hat g-\hat a)(\sigma)$, and this is at least $(\tau (b)(1- \tau(a)/\tau(b))-\epsilon > 0$. Hence $\hat b-\hat g -\hat a$ is strictly positive, and thus $b \geq g+a = c$. \qed

\Lem Corollary \oneeig. If $G$ is a simple dimension group other than $\Z$ and $\tau$ is a  trace, then $\tau$ is good if and only if for any (hence for every) order unit, the image of $\ker \tau$ is norm dense in $\brcs{\tau}^{\vdash}$. 

This   even applies when $\tau$ is the unique (up to scalar multiple) trace, since in that case $\brcs{\tau}^{\perp} = \brcs{0}$. It is a trivial application of unperforation that uniqueness of the trace implies goodness. A minor modification of the arguments above yield  semi-characterizations of weakly good elements (\onenin\ and \oneten). In examples (to follow), we usually deal with pure traces, because it is easier to calculate $\ker \tau$ and $\brcs{\tau}^{\perp}$. 

An order preserving group homomorphism between dimension groups with order unit  $\Arrow \phi; (G,u).(H,v)$ (\st $\phi(u) = v$) is called an {\it extended homomorphism\/} (or {\it extended\/}) if for all pure normalized traces $\tau$ of $(G,u)$, there exists a unique (normalized) trace $\sigma$ of $(H,v)$ \st $\tau = \sigma \circ \phi$. (Uniqueness of $\sigma$ guarantees it is pure, as the pre-image of an extreme point is a closed face.)

This entails that $\phi^{-1}(H^{++}) = G^{++}$ (that is, for $g$ in $G$, if $\phi(g)$ is an order unit of $H$, then $g$ is an order unit of $G$), but does not imply $\phi$ is an order-embedding---for example, if $G$ is simple and $\Inf(G) \neq 0$, then the quotient $ H = G/\Inf(H)$ is a simple dimension group (for general dimension groups, it is {\it not\/} true that if $G$ is a dimension group, then $G/\Inf(G)$ is as well), and the quotient map is extended. (However, if $G$ and $H$ are simple dimension groups and $\Inf (G) = 0$, then an extended map from $G$ to $H$ {\it is\/} an order embedding.)

Examples arise when $G$ is a simple and noncyclic dimension group, and $H$ is a dense subgroup of $\text{Aff\,}(S(G,u))$ containing $\hat G$ (and $H$ is equipped with the strict ordering); then $\phi$ is just the map \ $\hat{ }$ \  followed by the inclusion. When $G$ and $H$ are simple, modulo infinitesimals and order automorphisms, this is the only type of extended map between them. 

In general, extended maps induce a natural affine homeomorphism $S(G,u) \leftarrow S(H,v)$
given by $\sigma \mapsto \sigma \circ \phi$. This preserves the extremal boundaries (the sets of extreme points), and induces an isometric isomorphism between their affine function spaces.  The following corollary is somewhat unexpected---for example, if $G$ is  a dense subgroup of $\text{Aff\,}(K)$ where $K$ is a Choquet simplex, then we can take $H$ to be any intermediate group, $G \subseteq H \subseteq \text{Aff\,}(K)$ (both $G$ and $H$ are equipped with the strict ordering, so they are simple dimension groups), and there is no obvious reason why goodness of $\sigma|G$ should imply that of $\sigma|H$ when $\sigma$ is an extreme point of $K$. However, this is a consequence of the next result, since $G \subseteq H$ is an extended map.

\Lem Corollary \onesev. Let $\Arrow \phi;(G,u).(H,v)$ be  an extended homomorphism of   dimension groups. Suppose $\tau$ is a (normalized) trace of $(G,u)$ that is order unit good, and $\sigma$ is the (normalized) trace of $(H,v)$ \st $\tau = \sigma \circ \phi$. Then $\sigma$ is order unit good. 

\Pf It suffices to show $\ker \sigma$ has norm-dense image in $\brcs{\sigma}^{\vdash}$. However, we note that $\phi(\ker\tau) \subseteq \ker \sigma$, and from Proposition \onesix, $\ker \tau$ has dense range in $\brcs{\tau}^{\vdash}$. It is trivial that $\phi$ induces an isometric  isomorphism $\brcs{\tau}^{\vdash} \to \brcs{\sigma}^{\vdash}$, so that $\phi(\ker \tau)$ has dense image in $\brcs{\sigma}^{\vdash}$, and thus so does $\ker \sigma$. \qed

Obviously we cannot say anything when $\tau$ is not good (e.g., if $H$ is a real vector space, then all traces of $H$, even impure ones, are good). We will see in Example \exsix(b)  that extended maps do not generally preserve goodness; of course this phenomenon can only occur in the  non-simple case.

Here is a situation where there appears to be a difference between pure and impure traces, although not very much.

\Lem Lemma \onenin. Let $(G,u)$ be an approximately divisible dimension group with order unit, let $\tau$ be a  trace thereon, and let $v$ in $G^+$ be weakly $\tau$-good \st $\tau(v) > 0$.
\item{(i)} If $\tau$ is pure, then  for all $\epsilon$, for all $h$ in $\(\brcs{\tau}^{\perp}\)^+$, there exists $g$ (depending on $h$) in $\ker \tau$ \st $\| \hat g-h\| < \epsilon + \|\hat v\|$ (the norm is computed \wrt the order unit $u$). 
\item{(ii)}
If $v$ is $\tau$-good, then the same conclusion holds for all $h$ in $\brcs{\tau}^{\vdash}$. 

\Remark The difference is that a weaker assumption is made on $v$ in the impure case than in the pure case; we require purity with the weak assumption in order to guarantee positive elements of $\Aff S(G,u)$ can be approximated by images of elements of $G^+$.
\Pf (i) There exists $g_0$ in $G^+$ \st $\| \hat g_0-h\| < \epsilon$, by approximate divisibility.
Without loss of generality, $\epsilon < \tau (v)$, and thus $ 0 \leq \tau(g_0) <\epsilon<  \tau(v)$.  By weak goodness of $v$, there exists $g_1$ in $[0,v]$ \st $\tau(g_0) = \tau(g_1)$. Set $g = g_0-g_1$, so that $g \in \ker \tau$, and $\| \hat g-h\| \leq \|\hat g_0-h \| + \|\hat  g_1\| \leq \epsilon + \| \hat v \|$. 

\noindent (ii)
By approximate divisibility, there exists $g_1$ in $G$ \st $\| \hat g_1-h\| < \epsilon/3$.
Without loss of generality, $\epsilon < \tau (v)$, and thus $ |\tau(g_1)| <\epsilon/3<  \tau(v)$. There exists an order unit $w$ in $G^+$ \st $\epsilon/3 < \hat w < 2\epsilon/3$, so that 
with $g_0 = g_1 + w$, we have $0\leq \tau(g_0) < \epsilon$. By  goodness of $v$, there exists $g_2$ in $[0,v]$ \st $\tau(g_0) = \tau(g_2)$. Set $g = g_0-g_2$, so that $g \in \ker \tau$, and $\| \hat g-h\| \leq \|\hat g_0-h | + \|\hat  g_2\| \leq \epsilon + \| \hat v \|$. 
\qed

\Lem Corollary \oneten. (Order unit goodness criterion) Suppose  $(G,u)$ is an approximately divisible dimension group with  trace $\tau$. Either of the following conditions is sufficient  for $\tau$ to be order unit good.
\item{(i)} $\tau$ is  pure and there exists a sequence of weakly $\tau$-good positive elements of $G$, $\brcs{v_i}$, \st $\| \hat v_i\| \to 0$ and $\tau(v_i) > 0$ for almost all $i$.
\item{(ii)} there exists a sequence of  $\tau$-good positive elements of $G$, $\brcs{v_i}$, \st $\| \hat v_i\| \to 0$ and $\tau(v_i) > 0$ for almost all $i$.

Both conditions are of course necessary. This suggests a measurement of how far away a pure trace is from being order unit good; if $u$ is an order unit, use the norm obtained from the unital affine map $G \to \text{Aff\,} S(G,u)$, and define the element of $[0,\infty]$ via 
$$
\Cal G (\tau)\equiv {\Cal G_u} (\tau) = \inf \Set{\|v\|}{v \in G^+ \setminus (\ker \tau)\  \& \text{ $v$ is  good \wrt $\tau$}}.
$$
This depends on the choice of $u$, but not at the endpoints. Then $\tau$ is order unit good if and only if $\Cal G(\tau) = 0$, and no positive element outside $\ker \tau$ is weakly $\tau$-good if and only if $\Cal G(\tau) = \infty$. The latter applies when $\tau$ is bad. As a function on $\partial_e S(G,u)$ (the extremal boundary of the trace space), $\Cal G$ is mysterious (unless $G$ is a rational vector space, in which case it is $\brcs{0,\infty}$-valued). If $\tau$ is pure, we can define a similar function replacing {\it good\/} by {\it weakly good.}

\Lem Corollary \oneele. Let $(G,u)$ be a dimension group and $\tau$ a pure trace thereon. If there exists a weakly $\tau$-good $v$ in $G^+$ \st $\tau(v) \neq 0$, then the trace on $G \otimes \Q$ (the divisible hull of $G$) given by $\tau \otimes 1$, is order unit good. 

\Pf Since $G\otimes \Q$ is approximately divisible, it suffices to show each $v_n:= v \otimes \slfrac 1n$ is weakly $\tau$-good. However, it is  routine  to show that the interval $[0,v_n]$ in $G\otimes \Q$ is just $\slfrac 1n\cdot [0,v_1]$, and  $[0,v_1] = [0,v] \otimes ([0,1]\cap \Q)$. \qed

Here, $\Q$ can be replaced by any noncyclic subgroup of $\Q$. More criteria for goodness will be obtained in the next section. The converse of Corollary \oneele,   whether $\tau \otimes 1$ being good as a pure trace of $G \otimes \Q$ implies there exists a faithful weakly $\tau$-good element in $G$, seems unlikely, but as yet we have found no counter-examples. 

\SecT 2 A fistful of traces

 With the last few results in mind, we say a  trace $\tau$ on a dimension group is {\it ugly,} if the dimension group has more than one trace, the range of $ \ker \tau$ in $\brcs{\tau}^{\vdash}$ is discrete, and the trace $\tau \otimes 1$ on $ G\otimes \Q$ is good. In the case that $G$ is simple, the last condition is equivalent to the divisible hull of $\ker \tau$ (that is, its tensor product with the rationals) having dense range in $\brcs{\tau}^{\vdash}$. By [L2], a countable abelian group with a discrete norm is free; this applies to the image of $\ker \tau$ when the simple dimension group is countable (as here) and the pure trace $\tau$ is ugly. Thus  $\ker \tau$ splits as an abelian group direct sum of a free group and $\Inf(G)$.

A dimension group is {\it critical\/} (see [H1] for properties and examples) if it is simple, has $d$ pure traces, and is free of rank $d+1$. If $G$ is critical, then any subgroup of lesser rank has {\it discrete\/} image in $\R^d$. Hence none of its pure traces are good. 

For example, if $G$ is a critical group of the form $ G = \langle e_1, \dots, e_n; \sum_j \alpha_j e_j\rangle$ where $e_i$ are the standard basis elements in $\R^n$ and $\brcs{1, \alpha_1, \dots , \alpha_n}$ is linearly independent over the rationals, then all pure traces (the coordinate maps) are ugly. It is easy to verify that the element $u = \sum e_i = (1,1,\dots,1)$ is good \wrt every pure trace. More generally, if $v = (v_i) \in G$ with $v_i \geq 1$ for all $i$, then $v$ is good \wrt every pure trace, and thus the function $\Arrow \Cal G_u ;\partial_e S(G,u) = \brcs{1,2,\dots,n}. [0,\infty]$ defined earlier is everywhere at most one, and is probably the constant $1$. In this case, we have positive elements that are uniformly good, that is, simultaneously good \wrt every pure trace.

For other types of critical groups of rank $n+1$, we can arrange examples to admit no or only some ugly traces. 

We say a  trace is {\it bad\/} if $\ker \tau = \Inf (G)$ and the dimension group has more than one trace. This is preserved by tensoring with the rationals, so bad traces remain bad (not only is the trace not good, it's {\it ba-ad,} Proposition \onefou). 
 
Here is a different kind of  critical group with all of its pure traces bad. Let $f $ be an irreducible monic polynomial over $\Q$ with $d+1$ distinct real roots, $\brcs{r_i}_{i=1}^{d+1}$. Discard the last one and form the elements $v_j = (r_1^j,r_2^j, r_3^j, \dots, r_{d}^{j})$ of $\R^d$ for $j = 0,1, \dots, d$, and set $G = \langle v_j\rangle$. A simple Vandermonde determinant argument (executed in [EHS]), shows that $G$ is a dense subgroup of $\R^d$, so that equipped with the strict ordering, $G$ is a simple dimension group and a critical group. The pure traces are the coordinate functions, and since for each $i$, $\brcs{r_i^0, r_i^1, \dots, r_i^d}$ is linearly independent over $\Q$, it follows that each pure trace has zero kernel, that is, each pure trace is bad.

Bad traces are related to the notion of saturated minimal homeomorphisms, introduced [BK2], which refer to the set of all invariant measures, not just a single one.

We now construct two examples each with exactly three pure traces, such that one is good, one is ugly, and one is bad. Both are free abelian groups, the first of rank five (smallest possible for this arrangement); the second is of infinite rank, because it has an additional property: the value groups of the pure traces---$\tau(G)$---are identical. If $G$ has a bad trace and a not bad trace whose images in $\R$ are equal, then $G$ is isomorphic (as a group) to a proper quotient of itself, and being torsion-free, cannot be of finite rank.

\Lem Example \exthr. {\par}
\item{(a)} A simple dimension group free of rank five with exactly three pure traces: one is good, one is ugly, and one is bad.
\item{(b)} A simple dimension group $G$ with exactly three pure traces, $\tau_i$, \st one is good, one is ugly, and one is bad, {\it and\/} the value groups $\tau_i(G)$ are equal to each other.

\Pf (a) Let $\brcs{a,b,c,d,e,f,g}$ be an  {\it algebraically\/} independent set of seven real numbers. Define the following five elements of $\R^3$.
$$
\matrix
v_1 & = & (0 & b & g)\cr
v_2 & = & (0 & 1 & f)\cr
v_3 & = & (0 & d & e)\cr
v_4 & = & (a & b & c)\cr
v_5 & = & (1 & 1 & 1).\cr
\endmatrix
$$
Note that $1$ and $b$ are repeated in the second column.  Let $G$ be the subgroup of $\R^3$ generated by $\brcs{v_i}$. We first claim that $\langle v_5, v_4, v_3, v_2 \rangle$ is a critical subgroup, i.e., is dense: to check this, we note that the set of four $3\times 3$ determinants is linearly independent over $\Q$. Thus $G$ itself is dense, so with the strict ordering, is a simple dimension group.

Let $\tau_j$ ($j=1,2,3$) be the three coordinate maps. These are the pure traces; we will show that $\tau_1$ is good, $\tau_2$ is ugly, and $\tau_3$ is bad.

Since $\brcs{1,a}$ is linearly independent over $\Q$,  $\ker \tau_1$ is spanned by $\brcs{v_1, v_2, v_3}$; we claim $\ker \tau_1$ has dense image in $\R^2$ (after deleting the first coordinate). This corresponds to showing $\langle (b,g), (1,f), (d,e)\rangle$ is dense in $\R^2$; the three determinants are $\brcs{bf-g, be-gd, e-df}$ which is clearly linearly independent over $\Q$ (as the original set of seven elements is algebraically independent). Hence $\ker \tau_1$ has dense range in $\tau_1^{\perp}$, so $\tau_1$ is good.

As $\brcs{1,b,d}$ is linearly independent over $\Q$, $\ker \tau_2 = \langle v_1-v_4, v_2- v_1\rangle$. After deleting the middle zeros, we must show  the subgroup
$\langle (-a, g-c), (-1,f-1)\rangle$ is discrete and spans $\R^2$, i.e., the set $\brcs{(-a, g-c), (-1,f-1)}$  is a basis of $\R^2$. It is  linearly independent (as the determinant is $a(1-f)+ c-g$, obviously not zero), thus a basis. Hence $\tau_2$ is ugly.

Finally, $\ker \tau_3$ is zero, since $\brcs{1,c,e,f,g}$ is linearly independent over the rationals. So $ \tau_3$ is bad.

Since critical groups cannot have a good trace,  this example is the smallest possible (at least in terms of rank) among free abelian simple dimension groups.

\noindent (b) Let $G$ be the free abelian group on the countably infinite set $\brcs{x_i}_{i \in \N}$. Let a bijection $\N \to \N \times\brcs{1,2}$ be denoted $i \mapsto (s(i), l(i))$, and indicate a bijection $\N \times \brcs{1} \to \N \times \brcs{3,4}$ by $(i,1) \mapsto (t(i), m(i))$. Let $\brcs{\alpha_i}_{i \in \N}$ and $\brcs{\beta_i}_{i \in \N}$ be infinite  subsets of $\R$  whose union is algebraically independent over $\Q$ (this is probably overkill; it is likely we could get by with linear independence over $\Q$). 

Now we define three group homomorphisms $\Arrow \tau_j; G.\R$, which will turn out to be the desired pure traces. They are of course uniquely determined by their values on $\brcs{x_i}$.
$$\eqalign{
\tau_1 (x_i) & = \cases 0 & \text{if $l(i) = 2$}\\
\alpha_{t(i)} & \text{if $l(i) = 1$ and $m(i) = 3$} \\
\beta_{t(i)} & \text{if $l(i) = 1$ and $m(i) = 4$} \\
\endcases \cr
\tau_2 (x_i) & = \cases 0 & \text{if $i = \min\Set{j}{l(j)=1}$}\\
 0 & \text{if $i = \min\Set{j}{l(j)=2}$}\\
\alpha_{s(i')} & \text{where $i'$ is the immediate predecessor to $i$ in $l^{-1}(1)$} \\
\beta_{s(i')} & \text{where $i'$ is the immediate predecessor to $i$ in $l^{-1}(2)$} \\
\endcases \cr
\tau_3 (x_i) & = \cases \alpha_{s(i)} & \text{if $l(i) = 1$}\\
\beta_{s(i)} & \text{if $l(i) = 2$} \\
\endcases \cr
}$$

We first observe that $\tau_j (G) = \sum \alpha_i\Z + \sum \beta_i \Z$, so all three group homomorphisms have the same value group. Now consider the map $\Arrow T:= (\tau_1, \tau_2, \tau_3); G. \R^3$. Then  $\ker \tau_3 = 0$ is immediate from the construction (this requires only linear independence over $\Q$), and we will show that $\ker \tau_2$ has discrete range of rank two in $\R \times \brcs{0} \times \R$ and $\ker \tau_1$ has dense image in $\brcs{0} \times \R^2$. Finally, we show that the image of $G$
 is dense in $\R^3$, so that equipped with the strict ordering from the map $T$, $G$ becomes a simple dimension group. It is then routine that the pure traces are exactly the $\brcs{\tau_j}$.

Let $i_1 = \min l^{-1}(1)$ and $i_2 = \min l^{-1}(2)$; note that $i_1 \neq i_2$ because the relevant map is a bijection. Then linear independence of $\brcs{\alpha_i, \beta_j}$ yields $\ker \tau_2 = \langle x_{i_1}, x_{i_2} \rangle$. The respective images under $T$ are $(\gamma, 0, \alpha_{s(i_1)})$ and $(0,0,\beta_{s(i_2)})$, where $\gamma$ is one of $\brcs{\alpha_{t(i_1)}, \beta_{t(i_1)}}$---all that matters at the moment is that $\gamma \neq 0$. Then $\brcs{T(x_1), T(x_2)}$ is linearly independent over $\R$, so that it generates a discrete subgroup of $\R^3$. On the other hand, the subgroup it generates sits inside $\R \times \brcs{0} \times \R$, hence its real span is all of this, and thus its rational span is dense in it.

Now $\ker \tau_1$ contains $\langle x_j \mid l(j) = 2\rangle$ (in fact, the latter is all of $\ker \tau_1$, but we do not need this). If $j \neq i_2$, then for $j \in l^{-1}(2)$, $T(x_j) = (0, \beta_{s(j')}, \beta_{s(j)})$, and $T(x_{i_2}) = (0,0,\beta_{s(i_2)})$. Computing  a few size two determinants, we see easily that the span of $\ker \tau$ is dense in $\brcs{0} \times \R^2$ (this is where we use algebraic independence, rather than just linear independence, over the rationals, but we could get around this by specifying some of the values). 

From the preceding, the closure of  $T(G)$  contains  $\brcs{0} \times \R^2$, and it also has lots of elements whose first coordinates generate a dense subgroup of $\R$. Since the closure is a group, it follows that  $T(G)$ is dense in $\R^3$.

Thus we can define the ordering on $G$ via $G^+ \setminus \brcs{0} = T^{-1}((\R^3)^{++})$, that is, a nonzero element of $G$ is in the positive cone if and only if its image under $T$ is a strictly positive triple. Density of $T(G)$ ensures that $G$ is a simple dimension group, and it is clear that $\tau_i$ are all traces. Moreover, since the ordering is determined by the $\tau_i$ and they form a linearly independent set, it is trivial that they are  pure, and the only pure traces of $G$. 

Now knowing that $\tau_i$ are pure traces, what we showed above is that they are respectively bad ($\tau_3$), ugly ($\tau_2$), and good ($\tau_1$).
\qed 

A criterion which can be useful (see the argument for  tensor products, Proposition \fivtwo) is the following. 

\Lem Lemma \twoone. Let $G$ be a approximately divisible dimension group  and $\tau$ a pure normalized trace \wrt a fixed order unit. Suppose  in the corresponding affine representation, we have:{\par}
\item{}Given $\epsilon' > 0$,  there exists $0 < \epsilon \leq \epsilon'$ for which, given $r \in \tau(G^+) \cap (0,\epsilon)$, there exists $g$ in $G$ \st $\tau(g) = r$ and $\| \hat g\| < f(\epsilon)$ where $\lim_{t \downarrow 0} f(t) = 0$. {\par}
\noindent Then the image of $\ker \tau$ is dense in $\brcs{\tau}^{\perp}$, and thus $\tau$ is order unit good.

\Pf It suffices to approximate nonnegative elements of $\brcs{\tau}^{\perp}$ by images of elements of $\ker \tau$, since $\brcs{\tau}^{\perp}$ is an order ideal, and is thus the abelian group generated by its positive cone. Since $\hat G^+$ is dense in the set of nonnegative affine functions, for $h$ in $\(\brcs{\tau}^{\perp}\)^+$, given one of the $\epsilon$\,s that arise from some $\epsilon'$, there exists $g_1$ in $G^+$ \st $\|\hat g_1-h\| < \epsilon$. Then $\tau (g_1) = (\hat g_1-h)(\tau) < \epsilon$. Hence we may find $g_2$ in $G$ \st $\tau(g_2) = r$ and $\| \hat g_2\| < f(\epsilon)$. Then $w:= g_1-g_2$ is in $\ker \tau$, and $\| h-\hat w\| \leq \| h-\hat g_1\| + \| \hat g_2 \| < \epsilon + f(\epsilon)$. \qed

Note that the $g$ in the  displayed hypothesis  need not be positive, and moreover, if $f(\epsilon)$ is chosen to be something like $10\epsilon$, the hypothesis is a weak form of goodness. Here is a slightly different version of the same idea, almost tautological. 
\Lem Proposition \twotwo. (Density criterion) Suppose $(G,u)$ is an approximately divisible dimension group with order unit, and $\tau$ is a  pure trace. Sufficient density of the image of $\ker \tau$ in $\brcs{\tau}^{\perp}$ is:
{\par}
\item{} For any sequence  $\brcs{g_i}_{i \in \N}$  of elements of the positive cone \st $\tau(g_i) \to 0$ and $\brcs{\hat g_i}$ converges uniformly in $\text{Aff\,} S(G,u)$,   there exists an infinite subsequence $S \subseteq \N$ together with elements of $G$,  $\brcs{g_s'}_{s\in S}$,  \st  $\tau (g_s')= \tau (g_s)$ for all $s$ in $S$ and $\hat g_s' \to 0$ uniformly.

\Remark In the statement, $\tau$ need not be normalized \wrt $u$. The conditions are independent of the choice of order unit, since the  norms induced by different order units are mutually equivalent (but changing the order unit can change the isomorphism class of the Choquet simplex $S(G,u)$---fortunately, this is not relevant here).

\Pf For $h$ in $\(\brcs{\tau}^{\perp}\)^+$, density of $\hat G^+$ in $\text{Aff\,} S(G,u)$ (a consequence of  approximate divisibility) allows us to find a sequence  $\brcs{g_i}$ elements of $G^+$ \st $\hat g_i \to h$ uniformly; obviously $\tau(g_i) \to 0$, so we may apply the condition.  Setting $w_s = g_s-g_s'$, we see that the $w_s $ belong to $\ker \tau$ and $\| h- \hat w_s\| \leq \| h- \hat g_s\| + \| \hat g_s'\|$, so $\hat w_s \to h$ uniformly.
\qed

\SecT 3 A few traces more

Recall that a Choquet simplex $K$ is a {\it Bauer simplex,} if its set of extreme points, $\partial_e K$, is closed (hence compact).

\Lem Proposition \throne. Suppose that $A$ is a partially ordered ring with $1$ as order unit that is also an approximately divisible dimension group. Then{\par}
\item{(i)} A normalized trace of $(A,1)$ is pure if and only if it is multiplicative;
\item{(ii)} $S(A,1)$ is a Bauer simplex and $\text{Aff\,} S(A,1) $ can be identified with $C(X,\R)$ (with the sup norm), where $X = \partial_e S(A,1)$ consists of the pure traces equipped with the point-open (weak) topology;
\item{(iii)} a pure trace  $\tau$ is order unit good if and only if the ideal of $A$, $\ker \tau$, is not killed by any pure trace other than $\tau$. 

\Pf The only part of this that is not  either well-known or trivial is the sufficient condition for goodness in part (iii), that is, if $ \ker \tau$ does not vanish at any point of $X$ other than $\tau$, then the image of $\ker \tau$ is dense in $\Ann(\tau) = \Set{f \in C(X,\R)}{f(\tau) = 0}$. By (i), $\ker \tau$ is an ideal of $A$, therefore its closure in $C(X,\R)$ is an ideal thereof. A closed ideal in this algebra is necessarily of the form 
$\Ann(Y)$ for some closed subset $Y$; but the hypothesis forces $Y = \brcs{\tau}$, so $\ker \tau$ is dense in $\Ann(\tau)$.\qed

Examples (that we will discuss in more detail later) are rational polynomial rings  $A = \Q[x]$ or $\Q[x,y]$ which are embedded in $C(I,\R)$, $C(I\times J,\R)$ respectively, where $I$ and $J$ are compact subsets of the reals; the images are dense, so by imposing the strict ordering, $A$ becomes a simple dimension group, and it is easy to check that in all these cases, $A$ is a partially ordered ring in which $1$ is an order unit, and moreover, $\partial_e S(A,1)$ can be identified with respectively $I$ and $I \times J$. 

\Lem Example \exfou. Partially ordered fields with disconnected pure trace space. All pure traces are bad, except when there is just one.

\noindent This class of examples is from [H6]. Let $K$ be a subfield of the reals, equipped with the sum of squares ordering (that is, an element is in the positive cone if and only if it can be expressed as a sum of squares of elements of $K$). If $K$ is algebraic over $\Q$, then $K$ is a simple dimension group with $1$ as an order unit, and the pure traces given by the real embeddings of $K$. Thus the pure traces, being multiplicative, must have zero kernel, so {\it all\/} pure traces are bad. In particular, if $K$ is the field of all real numbers algebraic over $\Q$, the  pure trace space is an atomless totally disconnected compact set (a Cantor set). 

This is an example of an {\it extremely simple dimension group\/} (a simple dimension group for which all pure traces have zero kernel). It is known that  the pure trace spaces of these can contain no connected subset; however, this is of less significance than would be expected, since the trace space need not be a Bauer simplex. This class is contained in the class of archimedean simple dimension groups (equivalent to, $\tau (g) \geq 0$ for all pure traces $\tau$ implies $g \in G^+$), discussed more fully in [H6]. For these, the trace space can be any metrizable Choquet simplex. \qed

\vskip 3pt\noindent{\it Amusing remarks.} Suppose $G$ is a simple dimension group with two pure traces, $\tau$ and $\sigma$; for goodness, we may assume that $\Inf (G) = \brcs{0}$. 

\noindent (0) {\it General case with two pure traces.} The pure trace $\tau$ is good if and only if $\ker \tau$ is not cyclic (that is, $\tau$ is good if and only if $\ker \tau$ is either of rank more than one, or $\ker \tau$ is rank one and not finitely generated): since $\ker \tau \cap \ker \sigma = \Inf (G) = \brcs{0}$, $\sigma$ is one to one on $\ker \tau$, and will have dense range (as a subgroup of the reals) if and only if $\ker \tau$ is not cyclic.
Since $\brcs{\tau}^{\perp}$ is just a copy of the reals (---the closure of $\sigma(G)$, as   $\hat{\ }$ is given by $(\sigma,\tau)$), we are done.  It follows that for simple dimension groups with exactly two pure traces, every trace is one of good, ugly, or bad. 

\noindent (a) {\it $G$ is not free, but $\rank G = 2$.}  If $\ker \phi = 0$, $\phi$ is bad and thus obviously not good. However, if $\ker \tau \neq 0$, then $G$ is a group extension, $0 \to L = \ker \tau \to G \to \tau(G) \to 0$, and both $L$ and $\tau(G)$ are rank one. Even though $G$ is dense in $\R^2$, examples exist for which $L$ is discrete and nonzero, and in this case, the trace is ugly. 

Here is a special class of examples. Let $(a(i))$, $(b(i))$ be two sequences of positive integers, and set $M_j = \(\smallmatrix a(j) & b(j)\\ b(j) & a(j)\\ \endsmallmatrix \)$. Form $G = \lim \Arrow M_j ; \Z^2.\Z^2$. Then $G$ is simple, since the coefficients are all strictly positive. If we multiply any subcollection of the individual terms by $P = \(\smallmatrix 0 & 1\\ 1 & 0\\ \endsmallmatrix \)$, 
$G$ remains unchanged. Hence we may assume $b(i)\leq a(i)$ for all $i$. Then $G$ has either two pure traces or one pure trace. Set $\alpha_k = \prod_{i\geq k} (1- b(i)/(a(i) + b(i)))$. If $\alpha_k = 0$ for all $k$ (that is, $\sum b(i)/(a(i)+ b(i)) = \infty$), then $G$ has unique trace; if there exists $k$ \st $\alpha \equiv \alpha_k \neq 0$ (that is, $\sum b(i)/(a(i)+ b(i)) $ converges), then $G$ has two pure traces; these are constructed from the rows of the limiting matrix of $\prod_{j \geq k} M_j/(a(j)+ b(j))$. In this case, the automorphism of $G$ induced by $P$ simply switches the traces. Hence they are either simultaneously good or simultaneously bad. 

Now we claim that these traces are good if and only if $\alpha$ is rational, otherwise they are bad. It easily follows from the construction that the rank of the group of values is two if $\alpha$ is irrational, but one if $\alpha$ is rational. Hence $\ker \tau$ is zero if and only if $\alpha$ is irrational, and otherwise it is rank one but not cyclic (this is easy to see). Hence our criterion applies. 
Thus if $\alpha \neq 0$, almost certainly the traces will be bad (any number in the interval $(0,1)$ can be attained by $\alpha$ within this class of constructions).

Traditional examples (discussed in the literature) occur with $a(j) =2^j$ or $a(j) = j^2$ and $b(j) = 1$. 

\vskip 4pt Now  consider the case that $G$ be {\it free\/} of rank $n$ (and still with exactly two pure traces). Then $n \geq 3$. For 
every trace, $\rank \tau (G) \geq 2$, and of course, $\rank \ker \tau = n-\rank \tau(G)$. Hence $\tau$ will be good if and only if $\rank \tau(G) \leq n-2$. It will be ugly if and only if $\rank \ker \tau = 1$, which is easily arranged.

\vskip1pt\noindent (b) $n=3$. Then for every trace $\phi$, $\rank \phi (G) \geq 2$, so  neither trace is good.

\vskip1pt\noindent (c) $n=4$. Then $\phi$ is good if and only if $\rank \phi(G) = 2$. Thus we can arrange that one of the traces be good, the other no good. Let $\brcs{\alpha, \beta_1, \beta_2}$ be a set of real numbers \st each of the two sets, $\brcs{1, \beta_1, \beta_2}$ and  $\brcs{\alpha,\beta_1, \alpha\beta_1}$, is linearly independent over the rationals (this implies $\brcs{1,\alpha}$ is also linearly independent); let $G$ be the subgroup of $\R^2$ generated by $\brcs{(1,1), (\alpha,0),(0,\beta_1), (0,\beta_2)}$. Density in $\R^2$ comes from $\brcs{\alpha,\beta_1, \alpha\beta_1}$ being linearly independent; equip $G$ with the strict ordering, making it into a simple dimension group, with the two pure traces, $\sigma = \tau_1$ and $\tau = \tau_2$, the coordinate projections.

It follows from  linear independence  of $\brcs{1, \beta_1, \beta_2}$ over $\Q$ that $\ker \tau_2$ is spanned by $(\alpha,0)$ , so $\rank \ker \tau_2 = 1$ and thus $\tau_2$ is ugly, therefore not good. On the other hand, $\tau_1(G) $ is spanned by $\brcs{1, \alpha}$, so $\tau_1$ is good. We can obviously modify this example so that both traces are good. \qed

\Lem Example \exfiv. Simple dimension groups with a continuum of pure traces, for which the set of good traces, and the set of not good traces, are both dense, and the former is countable. 

\noindent Let $G = \Q[x]$, and let $I$ be a compact subset of $\R$; the typical (but not the only) examples for $I$ will be intervals, $[a,b]$, in the reals. Put the strict ordering on $G$ resulting from restricting to  $I$, that is, a nonzero polynomial $p$ is in the positive cone if an only if $p$ is strictly positive as a function on $I$. Density of $\Q[x] $ in $\R[x]$ (which by Weierstrass' theorem is dense in $C(I,\R)$) ensures that $G$ is a simple dimension group. The pure traces are given by evaluations at the points of $I$, that is $\tau_{\alpha}:p \mapsto p(\alpha)$ for $\alpha$ in $I$, and the norm on the affine function space translates to the sup norm on functions on $I$. It is easy to check that $\brcs{\tau_{\alpha}}^{\perp}$ is just the ideal consisting of the continuous functions on $I$ that vanish at $\alpha$. Assume $I$ has at least two points.

If $\alpha$ is transcendental, then $\tau_{\alpha}$ has kernel zero, so is bad. If $\alpha$ is rational, then the density criterion of Proposition \twotwo\ (which after an affine transformation, amounts to showing $x\Q[x]$ is dense in the continuous functions on $[0,1]$ that vanish at zero; this follows from $x\R[x]$ being dense) applies, so that $\tau_{\alpha}$ is good.

If $\alpha$ is algebraic but not rational, the situation is more interesting. Let $f_{\alpha}$ be the monic irreducible (over the rationals) polynomial of $\alpha$. Then $\ker \tau_{\alpha}$ is the ideal $f_{\alpha} \Q[x]$. If $f_{\alpha}$ has another root, say $\beta$, in $I$, then $\tau_{\beta}$ vanishes on $\ker \tau_{\alpha}$ (because the kernels are equal), so $\ker\tau_{\alpha}$ cannot be dense in $\brcs{\tau_{\alpha}}^{\perp}$. Thus when $f_{\alpha}$ has another root in $I$, $\tau_{\alpha}$ is not good (but not bad either). 

Finally, suppose that  $f_{\alpha}$ has just the one root in $I$; then $\tau_{\alpha}$ is good. To see this, note that it suffices to show $x-\alpha$ belongs to the closure of $f_{\alpha}\Q[x]$ (for then $(x-\alpha)\R[x]$ will also be contained in the closure, and this is dense in $\brcs{\tau_{\alpha}}^{\perp}$). Consider the function $f_{\alpha}/(x-\alpha)$; this is continuous on $I$ and has no zeros there by hypothesis. Its reciprocal, $h:= (x-\alpha)/f_{\alpha}$, is thus continuous, hence may be approximated by elements of $\Q[x]$, say $g_n \to h$ uniformly on $I$. Then $(g_n f_{\alpha})$ is a sequence of elements in $f_{\alpha} \Q[x]$ that converges to $x- \alpha$. 

To summarize: $\tau_{\alpha}$ is good if and only if $\alpha$ is algebraic and has no other algebraic conjugates in $I$.
For example, if $I = [a,2]$ for some $a 
\leq \sqrt 2$ and $\alpha = \sqrt 2$, then for our $G$
 (which now depends on $a$), $\tau_{\alpha}$ is good if and only if $-\sqrt 2 < a \leq \sqrt 2$ (if $a > \sqrt 2$, then $\tau_{\sqrt 2}$ is not a trace). The bad traces are exactly the evaluations at transcendentals in $I$; there are no ugly traces, since the dimension group is a rational vector space. 

We can put a finer  dimension group ordering on $\Q[x]$; take the pointwise ordering on $[a,b]$, that is $g \geq 0$ if and only if $g$ is nonnegative on the interval. This is a dimension group, with the same pure traces, but of course is not simple. In fact, if $\alpha$ is rational, then $\ker \tau_{\alpha}$ is a maximal order ideal. The $\brcs{\tau}^{\perp}$ are exactly the same, and the same characterization of good traces applies.

If $[a,b]$ contains no integers, then we can probably do the same with $G = \Z[x]$, but this situation is more complicated. \qed

In contrast, there is an explicit (and also an inexplicit) construction [H6]  of an archimedean simple dimension group (necessarily not an ordered ring) whose pure trace space is $[0,1]$, \st all evaluations at algebraic  points have zero kernel (and therefore are bad). Point evaluations at transcendentals may have kernels, but it is not clear whether any of them are good.

It is a plausible conjecture that  if $\tau$ is a pure trace on a simple dimension group $G$ with $\rank \tau (G) = 1$, then $\tau $ is good. This is true if there exists a subgroup $H$ together with a rank one noncyclic subgroup $U$ of $\Q$ \st $H \otimes U \subseteq G \subseteq H \otimes \Q$---if $\phi$ is a continuous linear functional on $\text{Aff\,}(S(G,u))$, then either $\phi(\ker \tau) = 0$ or $\phi(\ker \tau)$ contains $U\alpha$ for some nonzero $\alpha$; the latter entails $\phi(\ker\tau) $ is dense, the former that $\phi|G $ is a scalar multiple of $\tau$ (since $\phi$ induces a homomorphism to the reals on the rank one group $\tau(G)$; any two nonzero such must be scalar multiples of each other), which is excluded. Hence the image of $\ker \tau$ is dense in $\brcs{\tau}^{\perp}$. 

It is also true if the map $\Arrow  \tau; G.U \subseteq \Q \alpha$  (where $U$ is rank one, but not cyclic) is almost split, that is, if there is a subgroup $J$ of $G$ \st $(\ker \tau) \cap J = \brcs{0}$ and $G/(J + \ker \tau)$ is finite. Finite index means the problem reduces to $G_0 = J \oplus  \ker \tau$, and since $\tau$ is essentially the projection onto $J$ followed by multiplication by a real number, it easily follows (from the splitting) that $\hat G$ dense in the affine function space implies $\ker \tau$ has dense range in $\brcs{\tau}^{\perp}$. 

However, the conjecture is not true in general. The following example illustrates a few interesting points. 

\Lem Example \extwe. A simple dimension group with two pure traces, each ugly and having values in $\Z[\frac12]$; moreover, the underlying abelian group is strongly indecomposable of rank two. 

\noindent Let $G_0$ be the free abelian group on a countably infinite set, $\brcs{x_i}_{i=0}^{\infty}$. Let $\Arrow k; \Z^+ . \brcs{0,1}$ be a function \st the collection of initial segmentsÊ ofÊ all $0-1$ strings of the form $k(i)k(i-1)\dots k(0)1$ includes all finite sequences infinitely often, and suppose $k(0)= 0$. Define the rational numbers $\alpha_0 = 1$ and $\alpha_{i+1} = (\alpha_i + k(i))/2$, and now define two real-valued functions $\tau_1$, $\tau_2$ from $G_0$ via
$$
\tau_1 (x_i) = \frac 1{2^i} \qquad \tau_2(x_i) = \alpha_i.
$$
\noindent (i) The range of $\Arrow T= (\tau_1, \tau_2); G_0.\R^2$ is dense; let $J$ denote the closure of the image of $G_0$. We observe that $T(x_i) = (2^{-i}, \alpha_i)$. For any real number $\beta$ in the interval $[0,1]$,Ê there exists an infinite subset $S_{\beta}$ of $\Z^+$ \st $\lim_{s \in S_{\beta}\to \infty} \alpha_s = \beta$. Then $(0, 1/2)$ and $(0, 1/e)$ are limit points of $\brcs{T(x_s)_{S_{1/2}}}$ and $\brcs{T(x_s)_{S_{1/e}}}$ respectively, so that $\brcs{0} \times \(\frac 12 \Z + \frac 1e \Z\)$ is in $J$, whence $\brcs{0} \times \R$ is in $J$. This means that $(2^{-i},0) \in J$, whence the closure of the subgroup it generates contains $\R \times \brcs{0}$. Since $J$ is a group, $J = \R^2$.

\noindent (ii) With the strict ordering from $T$, $G_0$ becomes a simple dimension group, and the only pure traces are $\tau_1$ and $\tau_2$.

\noindent (iii) $\ker \tau_1 = \langle \brcs{2x_i - x_{i-1}}\rangle $. Let $g= \sum_{j=0}^k m_j x_j \to 0$ under $\tau_1$ with $m_k \neq 0$. Write $m_k x_k = \flo {m_k/2}(2x_{k} - x_{k-1})+ \epsilon(k) x_k +Ê \flo {m_k/2} x_{k-1}$, where $\epsilon(k) \in \brcs{0,1}$. This reduces $g$Ê to $\sum \epsilon(i) x_i$. Since this maps under $\tau$ both to zero and to the binary number 
$\sum \epsilon(i)2^{-i}$, it follows that all the $\epsilon(i) = 0$, so that $g$ is in the group y $ \langle \brcs{2x_i - x_{i-1}}\rangle $.

\noindent (iv) The range of $\tau_2$ on $\ker \tau_1$ is cyclic. Since $2\alpha_i - \alpha_{i-1} \in \Z$, we thus have $\tau_2 (2x_i - x_{i-1}) \in \Z$, so $\tau_2 (\ker \tau_2) \subseteq \Z$. 

\noindent (v) Finally, $\tau_1 (G_0) = \Z[1/2^i] \subseteq \Q$, and since $\tau_2(\ker \tau_1)$ is cyclic, $\tau_1$ is in fact ugly, not good. 

This example can be modified to a somewhat more startling one. We check that $\ker \tau_2 = \langle \brcs{2x_{i+1} - x_i}_{k^{-1}(0)} \cup \brcs{2x_{j+1} - x_j - x_0}_{k^{-1}(1)}\rangle$, so that the image of $\ker \tau_2$ is also discrete. Now $\ker \tau_1 \cap \ker \tau_2 = \Inf (G_0)$ (since $\tau_1$ and $\tau_2$ are the only pure traces on $G_0$), and it is easy to check that $\ker \tau_1 \cap \ker \tau_2$ is spanned by $\brcs{2x_{i+1} - x_i}_{k^{-1}(0)} \cup \brcs{2x_{j+1} - x_j -(2x_{l+1} - x_l)}_{j,l \in k^{-1}(1)}$. Let $G$ be $G_0/\Inf (G_0)$, and let $\overline \tau_i$ be the traces corresponding to $\tau_i$. Then $\Arrow \overline T:= (\overline \tau_1, \overline\tau_2); G. \R^2$ is the affine representation mapping, and is a one to one map from $G$ to $\Z[1/2] \oplus \Z[1/2]$ (with the strict ordering on the latter). In particular, $G$ is rank two, and each of $\ker \overline \tau_i$ is cyclic (and not zero). Since $T(G_0) = \overline T(G)$ is dense in $\R^2$, it follows that $G$ is indecomposable, in fact every finite index subgroup is indecomposable as an abelian group. 
\qed

\comment
Possible alternative proof of necessity, part (i) (outline). If density fails, there exists $h$ in $M(\tau)$, together with $\delta > 0$ \st for all $g$ in $\ker \tau$, $\| \hat g-h\| > \delta$ (that is, for all normalized traces $\sigma$, $|\sigma(g)-h(\sigma)| > \delta$. Suppose $\inf_{\sigma \in S(G,u)} h(\sigma) = s$, which we may assume to be negative; then $h-s\hat u$ is nonnegative, hence can be approximated by elements of the image of $G^+$, say $\| h-s\hat u-\hat a\| < \delta/4$ with $a$ in $G^+$. We have $|\tau(a) + s| < \delta/4$ and $\tau(a)> 0$, so that $|s| -\delta/4 <\tau (a)< |s| + \delta/4$. rest??
\endcomment

\SecT  4 Extensions 

This section deals with extensions of traces defined on order ideals of a dimension group $G$ to $G$ itself---whether goodness persists under this process. This generalizes [xxx, xxx].

\Lem Lemma \fouone. Let $G$ be dimension group with order unit $u$, and $J$  an order ideal in $G$. Let $\sigma$ be a trace on $J$. Then $\sigma$ extends to a trace on $G$ if and only if
$$
\sup \Set{\sigma(a) } {a \in J,\ \  0 \leq a \leq u} < \infty.
$$
If $\sigma$ is pure and it extends, then it can be extended to a pure trace.

\Lem Proposition \foutwo. Let $G$ be a dimension group with an order unit, let $J$ be an order ideal of $G$, let $\sigma$ be a pure trace of $J$, and $\tau$  a pure trace of $G$ extending $\sigma$. Then $\tau$ is good if and only if {\par}
\item{(i)} $\sigma$ is good and
\item{(ii)}  $\sigma(J^+) = \tau(G^+)$.

\Pf {\it Necessity.} (i) For $a, b$ in $J^+$ with $\sigma (a) < \sigma(b)$; then $\tau(a)< \tau(b)$, so there exists $0 \leq c \leq b$ in $G^+$ \st $\tau(a) = \tau (c) $; but $0 \leq c \leq b$ and $b$ belongs to the order ideal $J$, so $c$ belongs to $J$.

\noindent (ii) Normalize $\tau$ (and correspondingly $\sigma$) so that $\tau(u) = 1$. Suppose $r \in \tau(G^+) \setminus \brcs{0}$; there exists $a$ in $G^+$ \st  $\tau(a) = r$. Since $\sigma$ is not zero, there exists $b_0 $ in $J^+$ \st $\sigma(b_0) \neq 0 $, hence there exists positive integer $n$ \st $n\sigma(b_0) > r$; set $b = n b_0$ in $J^+$. Since $\tau$ is good, there exists $0\leq c \leq b$ (so $c$ belongs to $J^+$) \st $r= \tau(a) = \tau (c) = \sigma(c)$, hence $r \in \sigma(J^+)$.

\noindent {\it Sufficiency.} Select $a,b$ in $G^+$ with $\tau(b) > \tau (a) > 0$. By (ii), there exists $a',b'$ in $J^+$ \st $\sigma(a') = \tau(a)$ and $\sigma(b') = \tau(b)$. By purity of $\tau$, with $\epsilon < \tau(b)-\tau(a)$, there exists $0 \leq b'' \leq b',b$ \st $\tau(b'') > \tau(b)-\epsilon > \tau(a)$. From  $b'' \leq b'$, $b''$ belongs to $J$. Applying goodness of $\sigma$ to $b'', a'$, there exists $0 \leq c \leq b''$ with $\sigma( c )= \sigma(a')$. Thus $c \leq b$ and $\tau( c ) =\tau(a)$.\qed

Suppose $\tau$ is a pure trace of $G$. If there exists a minimal order ideal $I$ \st $\tau(I) \neq 0$, then $\tau$ is good if and only if $\tau|I$ is good and $\tau(I^+) = \tau(G^+)$. As $I$ is a minimal order ideal and the notion of order ideal is transitive, $I$ is simple as a dimension group. So we have  necessary and sufficient conditions for $\tau$ to be good (assuming $I$ is minimal and $\tau(I) \neq 0$), namely, $\tau(I^+) = \tau(G^+)$ and either $I = \Z$, or $\ker \tau|I$ has dense image in $\brcs{\tau|I}^{\perp}$ (\wrt one or  any order unit of $I$).

Part of the problem in trying to give a general characterization of good traces for nonsimple dimension groups is that the approximately divisible condition does not behave well on taking order ideals, and the corresponding density even less well.
However, we note that the proof of \onesix(ii) shows that if $b$ is an order unit and the density condition holds, then $b$ is $\tau$-good. This has a consequence: if $\tau$ is a pure trace, and for every order ideal with (relative) order unit $(I,b)$ (that is, $I$ is the order ideal generated by $b$), $\ker \tau|I$ has dense image in $\brcs{\tau|I}^{\perp}$ (the order ideal of $\text{Aff\,} S(I,b)$ consisting of the functions vanishing at $\tau|I$), and $\tau(I^+) = \tau(G^+)$, then $\tau$ is good.

In this result, we need only do this for one choice of order unit $b$ for each $I$. Moreover, we can weaken this further by using the generating result. Namely, let $H\subseteq G^+$  be a generating set. Look at all the order ideals with order unit of $G$; we can restrict our interest to the order ideals generated by each element of $H$.

In general, an order ideal $I$ (with order unit) in a  dimension group $G$ need not be approximately divisible when $G$ is. One trivial condition that forces order ideals (with order unit) to be approximatively divisible is that $G$ be divisible (for then the order ideals are themselves divisible), or $2$-divisible, or similar conditions (e.g., $G \iso G_0 \otimes U$ where $U$ is a noncyclic subgroup of the rationals). When every (or sufficiently many) order ideals are approximately divisible, the goodness of the restriction of a trace forces the density property. In this case, we would obtain necessary and sufficient conditions for $\tau$ to be a good trace.

For $h$ a positive element of $G$, let $I(h)$ denote the order ideal of $G$ generated by $h$, that is,
$$
I(h) = \Set{g \in G}{\text{there exists } K \in \N \text{ \st } -K h \leq g \leq Kh}.
$$

\Lem Proposition \fouthr. Let $G$ be a dimension group, and let $\tau$ be a pure trace. Let $H$ be a generating subset of $G^+$.{\par}
\item{(a)} If for all $h$ in $H$ \st $\tau(h) \neq 0$, (i) $\tau(I(h)^+) = \tau(G^+)$ and (ii) the image of $\ker \tau|I(h)$ is dense in the order ideal $\brcs{\tau|I(h)}^{\perp}$ of $\text{Aff\,}S(I(h),h)$, then $\tau$ is good.
\item{(b)} Suppose that all order ideals of $G$ with order units are approximatively divisible and $\tau$ is good. Then for all order ideals $(I,b)$ with order unit \st $\tau|I \neq 0$, it follows that (i) $\tau(I^+) = \tau(G^+)$ and (ii) the image of $\ker \tau|I$ is dense in the order ideal  $\brcs{\tau|I}^{\perp}$ of $\text{Aff\,}S(I,b)$.

These conditions are still unwieldy in general, but work perfectly well if (for example), $G$ has only finitely many order ideals, or if there exists a unique minimal order ideal not killed by $\tau$. 

Let $\tau$ be a trace on the dimension group $G$. Define $\ker^+ \tau$ to be abelian group generated by $ \ker \tau \cap G^+$ (that is,  the set of differences of such elements). It is an order ideal of $G$, and is the largest order ideal contained in $\ker \tau$. We may thus form $\overline G = G/\ker^+ \tau$, the dimension group obtained as the quotient by an order ideal; then $\tau$ induces a trace on $\overline G$ via $\overline \tau (g + \ker^+ \tau) = \tau(g)$. It is routine that $\tau$ is pure (as a trace of $G$) if and only if $\overline \tau$ is pure (as a trace of $\overline G$). It is equally routine (in fact the argument is almost identical and obvious) that the same result holds for goodness.

\Lem Lemma \foufou. Suppose $G$ is a dimension group, $\tau$ is a  trace thereon, and $I$ is an order ideal  of $G$ \st $I \subseteq \ker \tau$. If $H = G/I$ with the quotient ordering, and $\tilde \tau$ denotes the trace induced on $H$ by $\tau$, then $\tau$ is good (as a trace of $G$) if and only if $\tilde \tau$ is good (as a trace of $H$).

\Pf ($\Rightarrow$) Given $A$, $B$ in $H^+$ \st $\tilde \tau (A) < \tilde \tau (B)$, there exist $a$ and $b$ in $G^+$ whose images in $H$ are respectively $A$ and $B$. Apply goodness to $b$ and $a$, obtaining $0 \leq c \leq b$ \st $\tau(c)  = \tau (a) = \tilde \tau(A)$; then $C := c+ I$ is the desired element in $H$.

\noindent ($\Leftarrow$) Suppose $a$ and $b$ are in $G^+$ with $ 0 < \tau (a) < \tau(b)$. Let $A$ and $B$ be their respective images in $H^+$. By goodness of $\tilde \tau$, there exists $C$ in $H^+$ \st $ C \leq B$ and $\tilde \tau(C) = \tilde \tau (A) = \tau(a)$. There exists $c_1$ in $G^+$ \st $C = c_1 + I$. There exists $f$ in $H^+$ \st $c_1 \leq b + f$ (from the quotient ordering). By Riesz decomposition, there exist $0 \leq c \leq b$ and $0 \leq f_1 \leq f$ \st $c_1 = c + f_1$. Obviously, $\tau (c) = \tau(c_1) = \tilde \tau (A) = \tau(a)$, so $c$ is desired element.
\qed

\Lem Corollary \foufiv. If $\tau$ is a  trace \st $\ker^+ \tau = \ker \tau$, then $\tau$ is good. In particular, this applies if $\tau$ is a pure discrete trace.

\Pf If equality holds, then $\overline G$ is order isomorphic (via $\overline \tau$) to a subgroup of the reals, hence is totally ordered, hence has unique trace, which necessarily equals $\overline \tau$, and the latter is thus good. The preceding result then applies. For a pure discrete trace, $\ker^+  \tau = \ker \tau$ [G]. \qed

More important than the corollary (the condition on $\ker \tau$ seldom holds), is that we can combine this with our previous conditions to check whether $\tau$ is good. Given a pure trace $\tau$ on $G$, factor out $\ker^+ \tau$, and try to verify the previous conditions for the presumably less complicated $\overline G$. For example, if the latter has sufficiently many order ideals with order units \st the density condition on the relative kernels of $\overline \tau$ applies, and the range of $\overline \tau$ on the order ideals equals that on the whole group, then $\overline \tau$ (and thus $\tau$) is good.

\Lem Example \exsix. (a) A dimension group with a pure trace that is  order unit good but not good;
(b) an extended homomorphism of dimension groups $G \to H$ together with a pure trace $\tau$ that is good, but its extension to $H$ is not.

\Pf Let $K$ and $K'$ be dense subgroups of the reals, and form an extension (in the category of dimension groups) of $K$ by $K'$, following the prescription in [H] (developed further in [GH3]). Set $G(K,K') = K \oplus K'$ with the following ordering
$$
G(K,K')^+ = \Set{(k,0) \in K \oplus \brcs{0}}{k \in \R^+} \cup
\Set{(k,l) \in K\oplus K')}{k + l,l \in \R^{++}}.
$$
That this is a dimension group with exactly two pure traces, given by $\tau:(k,l)\mapsto k+l$ and $\sigma: (k,l) \mapsto l$ is folklore, but is discussed in  [H]; moreover, $J:= K \oplus \brcs{0}$ is the only proper order ideal of $G(K,K')$, and the quotient by $J$ is $K'$. Since $K$ is totally ordered and simple, $\tau|J$ is good. It is easy to check directly that $G(K,K')$ is approximately divisible (this also follows from [GH1] which asserts that a dimension group with order unit  and no discrete traces is approximately divisible).

Next we look at $\ker \tau = \Set{(k,-k)}{k \in K \cap K'}$. We note that  $\brcs{\tau}^{\perp} = \brcs{(r,-r)\in \R^2} \subset \R^2$. If $K \cap K' = \brcs{0}$, then $\tau$ is bad; if $K \cap K'$ is cyclic, then $\ker \tau$ cannot be dense in $\brcs{\tau}^{\perp}$, so that  $\tau$ is not order unit good, let alone good. However, if $K \cap K'$ is dense in $\R$, then $\sigma(\ker \tau) = K\cap K'$ is dense in the reals, so that $\tau$ is order unit good. Thus we have,
\itemitem{(i)} $\tau$ is order unit good if and only if $K \cap K'$ is dense in $\R$.

Since  $\tau(G) = K'+K$ but $\tau(J) =K $, we also have from Lemma \fouone, that
\itemitem{(ii)} $\tau$ is good if and only if $K' \subseteq K$;

Next, we observe that if $K' \subseteq K''$, then the inclusion $G(K,K') \subseteq G(K,K'')$ is an extended homomorphism, since the traces extend uniquely.

\noindent (a) Suppose that  $K'$ strictly contains $K$; then $\tau$ is not good by (ii), but is order unit good by (i). 

\noindent (b) The inclusion $G:= G(K,K) \subset H:= G(K,K')$ (where $K'$ is as in (a)) is extended; however, $\tau$ restricted to $G(K,K)$ is good by (ii), but by (i), its extension to  $\tau$ as a pure trace on $G(K,K')$ is not good.
\qed 

\SecT 5 Ordered tensor products 

If $G$ and $H$ are partially ordered abelian groups, then we may form their ordered tensor product (over $\Z$; this can be done over other ordered rings as well), $G \otimes H$, whose underlying group is that of their tensor product as abelian groups, with positive cone generated by elements of the form $\Set{g\otimes h}{g \in G^+, \ h \in H^+}$. That this {\it is\/} a partially ordered abelian group is not entirely trivial [GH2], but is trivial if both are dimension groups, and in the latter case, is also a dimension group (also trivial to show). 

Suppose both $G$ and $H$ have order units, $u$ and $v$ respectively. Then $u \otimes v$ is an order unit for $G\otimes H$, and the trace space is easily determined, although a little awkward. What is not awkward, however, is the pure trace space; this is the cartesian product of the pure trace space, that is, $\partial_e S(G \otimes H, u \otimes v) = \partial_e S(G, u ) \times \partial_e S( H,  v)$; the identification is via $(\sigma,\tau) \mapsto \sigma\otimes \tau$ (the latter acts via $\sigma\otimes \tau (g \otimes h) = \sigma(g) \tau(h)$). That is, every pure trace of the tensor product is a tensor product of pure traces, and every such tensor product is pure. Now we can ask when $\sigma \otimes \tau$ is good.

If $G$ and $H$ are simple (respectively, approximately divisible), then so is $G\otimes H$. Thus we can apply (some) results above. For example, suppose  $G$ and $H$ are simple, and that the two traces $\sigma$ and $\tau$ are both good \wrt their respective groups. We show $\sigma \otimes \tau$ to be good.

\Lem Lemma \fivone. Suppose $G$ is an approximately divisible dimension group with order unit. If $\tau$ is an order unit  good  trace, then for all $g$ in $G^+ \setminus \ker \tau$, there exists a positive element $g'$ of $G$ \st $\tau(g') = \tau(g)$ and $\| \hat g' \| \leq 2 \tau (g)$.

\Pf By density of $G$ in its affine representation (\wrt a fixed order unit), there exists $b$ in $G$ \st $\tau (g) \pmb 1 < \hat b  \leq 2 \tau(g) \pmb 1$ (possible since $\tau (g) \neq 0$). Then $\hat b$ being strictly positive implies $b$ is an order unit of $G$, in particular, $b$ belongs to $G^+$, and $\tau(b) >  \tau (g)$. Apply order unit goodness to $b$ and $g$; there exists $g'$ in $G^+$ \st $\tau(g') = \tau(g)$ and $ g' \leq b$. The latter entails $\hat g' \leq \hat b \leq 2 \tau(g) \pmb 1$, so $\| \hat g' \| \leq 2 \tau (g)$.
 \qed

\Lem Proposition \fivtwo. Let $G$ and $H$ be approximately divisible dimension groups with order units and with pure traces $\sigma$ and $\tau$ respectively. If $\tau$ and $\sigma$ are order unit good, then the pure trace $\sigma \otimes \tau$ of $G \otimes H$ is order unit good. In particular, if $G$ and $H$ are simple and $\sigma$ and $\tau$ are both good, then $\sigma \otimes \tau$ is good. 

\Remark If merely $G$ is approximately divisible (and both have an order unit), then $G\otimes H$ is approximately divisible. It is not immediately clear how to show that in this case, $\sigma \otimes \tau$ is good if both $\sigma$ and $\tau$ are. 

\Pf It is immediate that $G \otimes H$ is approximately divisible. We apply the density criterion to show that $\ker \sigma \otimes \tau$ is dense in $\brcs{\sigma \otimes\tau}^{\perp}$, and then \onesix\ applies. When $G$ and $H$ are simple, so is $ G\otimes H$, and thus order unit good is the same as good. 

Fix order units for each of $G$ and $H$, and use the corresponding affine representations; the tensor product of the order units gives an order unit for the tensor product of the groups. Given $g$ in $G^+$ and $h$ in $H^+$, there exist $g'$, $h'$  in $G^+$ and $H^+$ respectively \st $\sigma(g') = \sigma (g)$,  $\tau (h') = \tau (h)$, and $\| \hat g'\|  \leq 2 \sigma (g)$, and $\| \hat h' \| \leq 2 \tau (h)$. Then $g' \otimes h'$ is positive and $\| \widehat {g' \otimes h'}\| \leq 4 \sigma(g) \tau(h)$ \wrt the affine representation for $G \otimes H$; moreover,   $\sigma \otimes \tau (g' \otimes h') = \sigma \otimes \tau (g \otimes h)$.

Applying this to all the elements in a sum $z= \sum g_j \otimes h_j$ (with $g_j \in G^+$ and $h_j \in H^+$), we can find positive elements $g_j'$ and $h_j'$ \st with $z' = \sum g_j' \otimes h_j'$, we have $\sigma \otimes \tau\(z'\) = \sigma\otimes \tau(z)  $ and $\|\hat z' \| < 4 \sigma\otimes \tau (z)$. This is more than enough to verify one of the density criteria (Proposition \twotwo) for $\sigma \otimes \tau$: if $\sigma \otimes \tau( z_i) \to 0 $, then $\| \hat z_i' \| \to 0$.
   \qed

 Here is a rather surprising general result about tensor products of traces.

\Lem Proposition \fivthr. Let $\sigma$, $\tau$ be  normalized pure traces on the simple dimension groups with order unit, $(G,u)$, $(H,v)$, respectively, and suppose moreover that $\tau$ is good, and $\sigma(G^+) \tau (H^+) \subseteq \tau (H^+)$. Then $\sigma \otimes \tau$ is a good pure trace of $G \otimes H$.


\Remark  The condition says that $\tau(H)$ is a unital module over the subring of the reals generated by $\sigma(G)$ (since the traces are normalized, $1$ belongs to both $\sigma(G)$ and $\tau(H)$).

\Pf Suppose $z_i = \sum_j g_{ji} \otimes h_{ji}$ belong to $(G \otimes H)^+$ (with $g_{ji}$ in $G^+$ and $h_{ji}$ in $H^+$) and $\sigma\otimes \tau(z_i) \to 0$. Since $\sigma\otimes \tau (z_i) = \sum_j \sigma (g_{ij}) \tau(h_{ij})$, and the latter belongs to the set of finite sums of elements of $\sigma(G^+) \tau (H^+)$, hence belongs to $\tau (H^+)$, there exist $h_i$ in $H^+$ \st $\tau (h_i) = \sigma\otimes \tau (z_i)$, hence $\tau(h_j) \to 0$. Since $\tau$ good, there exists $h_i'$ in $H^+$ \st $\tau (h_i) = \tau (h_i')$ and $\| \hat h_i'\| < 2 \tau (h_i)$; thus $\hat h_i' \to 0$ uniformly, and it is immediate that if $w_i = u\otimes h_i'$, then  $\hat w_i \to 0$ uniformly \wrt the affine representation of $(G \otimes H, u\otimes v)$. Obviously $\sigma \otimes \tau (w_i) = \tau(h_i') = \tau (h_i) =\sigma\otimes \tau (z_i) $, and the density criterion is verified. Simplicity of the tensor product allows us to apply \oneeig.
\qed

Suppose in the preceding, we change the order unit of $G$, say replacing it $u'$ where $\sigma (u') = \lambda$. Then $\sigma$ is replaced by $\sigma': = \sigma/\sigma(u')$, and of course goodness of $\sigma \otimes \tau$ is equivalent to that of $\sigma' \otimes \tau$ (since goodness is defined without reference to any order unit). If $\sigma(G^+) \tau(H^+)/\lambda \subset \tau(H^+)$, then $\sigma'(G^+)\tau(H^+) \subset \tau(H^+)$, whence (assuming $\tau$ is good), it follows that $\sigma' \otimes \tau$ is good, and thus $\sigma \otimes \tau$  is good.
Thus the condition can be weakened to the following, which simply says that $\tau (H)$ is order isomorphic to a module over the ring generated by $\sigma (G)$.

\Lem Corollary \fivfou. Let $\sigma$, $\tau$ be  normalized pure traces on the simple dimension groups with order unit, $(G,u)$, $(H,v)$, respectively, and suppose moreover that $\tau$ is good, and $\sigma(G^+) \tau (H^+) \subseteq \lambda \tau (H^+)$ for some real $\lambda > 0$ in $\sigma (G^+)$. Then $\sigma \otimes \tau$ is a good pure trace of $G \otimes H$.

There is a necessarily weak result in the opposite direction. If $J$ is a subgroup of $\R$ containing $1$, let $F\langle J\rangle$ denote the smallest subfield of $\R$ containing $J$.

\Lem Lemma \fivfiv. Suppose $\sigma$ and $\tau$ are pure normalized traces on the simple dimension groups $(G,u)$ and $(H,v)$ respectively. Suppose  that $\tau$ is bad (that is, $\ker \tau = \brcs{0}$ and $H$ has more than one trace), and satisfies the following additional property:
{\par}\itemitem {} Whenever $\brcs{\tau(h_i)} $ is linearly independent over $\Q$ (for $h_i \in H$), then $\brcs{\tau(h_i)} $ is linearly independent over $F = F\langle \sigma(G)\rangle $.
{\par} \noindent Then $\sigma \otimes \tau$ is not good.

\Pf We show that $\ker \sigma \otimes \tau =( \ker \sigma ) \otimes H$---for which we may assume  that $H$ is divisible. If   $\tau'$ is a normalized trace on $H$ distinct from $\tau$, then $\sigma \otimes \tau'$ kills $\ker \sigma \otimes \tau$, and thus the latter cannot have dense image in $\brcs{\sigma \otimes \tau}^{\perp}$.

The argument is trivial. Suppose $w = \sum_{i=1}^l g_i\otimes h_i$ is in $\ker \sigma\otimes \tau$, and $l$ is the least number of terms to represent $v$. If $l = 1$, then $0= \sigma\otimes \tau (w) = \sigma(g_1) \tau(h_1)$; if $\tau(h_1) = 0$, then $h _1= 0$ and thus $w = 0$; if $\sigma(g_1) = 0$, then $w \in (\ker \sigma ) \otimes H$.

If $l > 1$, then $0 = \sum_i \sigma(g_i) \tau (h_i)$. None of $\tau(h_i)$ can be zero (else the corresponding $h_i$ would be zero, and we would obtain a shorter representation), and we claim $\brcs{\tau(h_i)}$ is linearly independent over $\Q$.
If not, by relabelling there exist rational numbers $q_i$ \st $\tau(h_l) = \sum_{i=1}^{l-1} q_i \tau(h_i)$. Hence there exists a positive integer $N$ \st $Nq_i$ are all integers and $\tau(N h_l-\sum (Nq_i) h_i) = 0$; as $\tau$ is one to one, this forces $N h_l =\sum_{i=1}^{l-1} Nq_i h_i $.
 Then
$$\eqalign{
Nw &= g_l \otimes Nh_l + \sum_{i=1}^{l-1} g_i \otimes Nh_i\cr
& = \sum_{i=1}^{l-1} (g_{i} + g_l) \otimes ((N+ Nq_i)h_i  ). \cr
}$$
Since $H$ is divisible, we may find $j_i$ in $H$ \st $Nj_i =(N+ Nq_i) h_i $, so that $w =  \sum_{i=1}^{l-1} (g_{i} + g_l) \otimes j_i $, a representation with fewer terms. Hence $\brcs{\tau(h_i)}$ is linearly independent over $\Q$.

By hypothesis, $\brcs{\tau(h_i)}$ is linearly independent over $F$. However, $0 = \sigma \otimes \tau(v) = \sum \sigma(g_i ) \tau(h_i)$.  Since $\sigma(g_i)$ belongs to $F$, this forces all $ \sigma (g_i) = 0$, so that $v \in ( \ker \sigma ) \otimes H$.
\qed

The hypotheses are ridiculously strong. However, it is not clear how to significantly weaken them. For example, if the   linear independence hypothesis is reduced  to $F\langle \sigma(G)\rangle \cap F\langle \sigma(H)\rangle = \Q$, the result fails, as illustrated by the following example.

\Lem Example \exsev. A simple dimension group $(G,u)$ with a bad pure trace $\tau$ together with a subfield $E$ of the reals \st $F\langle\tau(G)\rangle \cap E = \Q$, but the pure trace $1 \otimes \tau$ of $E \otimes G$ is good.

\Pf Pick an irreducible cubic over the rationals with three real roots, $r$, $s$, $t$ for which $\Q[r]$ is not Galois; hence the splitting field must have Galois group $S_3$, and $s \not\in \Q[r]$. Form the subgroup of $\R^2$, $ G = \langle (1,1), (r,s), (r^2,s^2)\rangle$; this is free of rank three, and a dense subgroup of $\R^2$; therefore, with the strict ordering, it is a dimension group (cf, [EHS]; it can alternatively be described as $\Z[r]$ equipped with the  inherited ordering from the field $\Q[r]$ induced by two of the three real embeddings). Since each of $\brcs{1,r,r^2}$ and $\brcs{1,s,s^2}$ is linearly independent over the rationals, both pure traces (arising from the coordinates) are bad. Let $\tau = \tau_2$, the projection on the second coordinate.

Now set $E = \Q[r]$ equipped with the total ordering inherited from $\R$. We note that $r+s+t \in \Q$, so that $\Q[r,s]$ is the splitting field, and its dimension is six. Hence $s$ satisfies a quadratic over $\Q[r]$. This immediately yields a nonzero kernel for $1\otimes \tau_2$ on $E \otimes G$, and since the tensor product is a vector space over $\Q[r]$ (induced by the first factor), the kernel is a rational vector space, and therefore the only other trace, $1\otimes \tau_1$ restricted to the kernel
has dense range in $\R$, and therefore $1 \otimes \tau_2$ is good.

Not relevant to the hypotheses of the lemma, but it is immediate that $1 \otimes \tau_1$ is good, by \fivthr.
\qed

Other combinations (for example, both $\tau$ and $\sigma$ not good, but $\sigma \otimes \tau$ good) are interesting, and anything can happen. As an easy example, suppose 
$G$ is a simple dimension group with two pure traces and no infinitesimals, and one of the traces, $\tau$ has the property that $\ker \tau$ is nonzero and cyclic. Then as we have seen, $\tau$ is not good. However, if $H = \Q$, then the induced trace $ \tau \otimes 1$ has kernel $\Q$, and as there are no infinitesimals, and $ G \otimes \Q$ still has just two pure traces, $ \tau \otimes 1$ is good. In other words, tensoring with the rationals has improved $\tau$ to a good trace. However, if $\ker \tau = 0$, then the kernel of $ \tau \otimes 1$ is still zero, so $ \tau \otimes 1$ is not good. We can arrange an example where $G$ is free of rank three with two pure traces with one trace having nonzero kernel, the other not, so that one of the traces improves to good, the other one does not, on tensoring with the rationals. 

More interestingly, we can easily construct an example of a simple dimension group $G$ \st none of its pure traces are good, but on $G\otimes G$, all of them are good. In this example (about the simplest possible), $G$ is free of rank three with two pure traces, so that $G \otimes G$ has $2 \times 2 = 4$ pure traces. 

\Lem Example \exeig. A simple dimension group $G$ all of whose pure traces are ugly (and thus not good) \st all   pure traces of $G \otimes G$ are good. 

\noindent Let $G = \langle (0,1), (1,0), (\alpha,\beta)\rangle$ inside $\R^2$ with $\brcs{1,\alpha,\beta}$ linearly independent over the rationals; equipped with the strict ordering, this is Example \exone. Then $G$ is a simple dimension group, and its pure traces are the coordinate maps, call them $\tau_1$ and $\tau_2$. We already know that both traces are ugly. The four pure traces on $G \otimes G$ are  $\brcs{\tau_i \otimes \tau_j}_{i,j=1,2}$. 

We will identify $\R^2 \otimes \R^2$ with $\R^4$, via $(a,b) \otimes (c,d) \mapsto (ac, ad, bc, bd) $. It is relatively easy to calculate enough of each $\ker \tau_i \otimes \tau_j$ to show that  they are dense in the copy of $\R^3$ obtained by ignoring the entry corresponding to $(i,j)$. 

With this identification, $G \otimes G$ is the subgroup of $
\R^4$ generated by the following nine elements:
$$
\matrix v_1 & = & (0 & 0 & 0 & 1)\\
v_2 & = & (0 & 0 & 1 & 0)\\
v_3 & = & (0 & 1 & 0 & 0)\\
v_4 & = & (1 & 0 & 0 & 0)\\
\endmatrix
\qquad \qquad
\matrix 
v_5 & = & (0 & 0 & \alpha & \beta)\\
v_6 & = & (\alpha & \beta & 0 & 0)\\
v_7 & = & (0 & \alpha & 0 & \beta)\\
v_8 & = & (\alpha & 0 & \beta & 0)\\
v_9 & = & (\alpha^2 & \alpha\beta & \alpha\beta & \beta^2)\\
\endmatrix
$$
The four pure traces, denoted $\phi_k$ ($k = 1,2,3,4$) are the coordinate projections. Just by looking at the locations of zeros, we see $ \brcs{v_1,v_2, v_3, v_5, v_7} \subseteq \ker \phi_1 $ (of course, $v_8-v_6$ also belongs, but we do not need it). Goodness will occur if this subgroup is dense in the copy of $\R^3$ obtained by deleting the first coordinate. The presence of the real basis ($\brcs{v_1, v_2, v_3}$) makes density  relatively easy to calculate. 

The closure of the group span of $\brcs{v_1, v_2, v_5}$ is $0 \times \R^2$ (we ignore the leftmost coordinate, which is automatically zero), hence $(\alpha,0,0)$ (arising from $v_7$) is in the closure, and together with $v_3$, this yields $\R \times 0 \times 0$ in the closure of $\ker \phi_1$. Since the closure of a group is a group, the closure is all of $\R^3$ and the trace $\phi_1$ is thus good.

For $\phi_2$, we have $\brcs{v_1,v_2, v_4, v_5, v_8} \subseteq \ker \phi_2$. Delete the second coordinate, and again we check density: $\brcs{v_2,v_4,v_8}$ generates a dense subgroup of $0 \times \R^2 $, and using this and $v_8$ (and recalling that the closure of group is a group), we see that $(\alpha,0,0)$ belongs to   the closure of $\ker \phi_2$, and therefore (with $v_1$), $\R \times 0 \times 0$ contained in the closure, and again this means that the closure is all of $\R^3$. 

Flipping the second and third coordinates, then interchanging $\alpha$ with $\beta$ yields the same result for $\phi_3$ (from that for $\phi_2$), and simillarly, flipping the first and fourth coordinates and interchanging $\alpha$ and $\beta$ yields the same for $\phi_4$ from that of $\phi_1$. 

So all four pure traces of $G \otimes G$ are good. We did not even use $v_9$. \qed

Still open: are all pure traces of $G \otimes G \otimes G$ good (same $G$ as in Example \exeig)? This is true if $\alpha$ is a cubic algebraic integer and $\beta = \alpha^2$, as we will see below. In other cases, it is not so clear. 

Another choice for $G$ as free rank three and simple has $G = \langle (1,1), (r,s), (r^2,s^2) \rangle$ where $r$ and $s$ are two of  three real roots of an irreducible cubic (so that $G$ is the ring of algebraic integers, $\Z [\alpha]$, where $\alpha$ satisfies a cubic with all roots real, with the ordering from two of the three real embeddings; this type of example was discussed in [EHS], and later in [G] and [H6]). In this case, there is a symmetry that causes $\ker \phi_1 = \ker \phi_4$ and $\ker \phi_2 = \ker \phi_3$, so that none of the pure traces are good ($\phi_{5-i}$ vanishes on $\ker \phi_i$, so density of the kernel in $\brcs{\phi_i}^{\perp}$ is impossible). 

If $G$ is generated by three elements of $\R^2$ \st the set of six real numbers among its coordinates is algebraically independent, then all four traces of the tensor product have zero kernel, so are all not good.

In Example \exeig, we can easily arrange $\sigma(G^+) \tau(G^+) \not\subseteq \tau(G)$ or $\sigma(G^+) \tau(G^+) \subseteq \tau(G)$, as we wish, by ensuring respectively that $\brcs{1,\alpha,\beta,\alpha\beta }$ be linearly independent, or $\brcs{\alpha^2,\beta^2,\alpha\beta }$ be in the $\Z$-span of $\brcs{1, \alpha, \beta}$; e.g., $\alpha$ a cubic algebraic integer, and $\beta = \alpha^2$. If  the last condition holds, then all pure traces of $G \otimes G \otimes G$
 are good, since the condition in \fivthr\ is now satisfied.

\comment 
And of course, our example with  a pure trace on $G$ not good, but its extension to the rational vector space $G\otimes \Q$ being good, also fails to satisfy the main hypothesis. So   $\sigma \otimes \tau$  may be good for other reasons. In addition if that hypothesis is satisfied, we need that $\tau$ be good (the $\Z[\alpha]$ example). 
\endcomment

Now we consider an interesting example arising from the tensor product of two polynomial rings with strict ordering inherited from points of $\R$. As in, let $I$ and $J$ be compact subsets of $\R$ each having at least two points (typically $I$ and $J$ will be intervals, but need not be). Form the polynomial rings $G = \Q[x]$ and $H = \Q[y]$, and impose on them the strict ordering viewing the elements as functions on $I$ and $J$ respectively. Then each is a simple dimension group (as in  Example \exfiv), and their pure traces are the point evaluations $\brcs{\sigma_{\alpha}}_{\alpha \in I}$, $\brcs{\tau_{\beta}}_{\beta \in J}$. 

\Lem Example \exnin. Good pure traces on $\Q[x,y]$ with the strict ordering from $I \times J$, where $I$ and $J$ are compact subsets of $\R$. 

\noindent The tensor product is naturally order and ring isomorphic to the polynomial ring in two variables, $\Q[x,y]$, equipped with the strict ordering from restricting to $I \times J$. To see this, the natural map $f\otimes g \mapsto fg$ extends to a ring isomorphism, it is positive, and the map induces a bijection on the pure traces. Since both $G\otimes H$  and $\Q[x,y]$ are simple, it follows that the map must be an order isomorphism. So we can identify $G \otimes H$ with $\Q [x,y]$ in the obvious way.

The pure traces of the latter must be multiplicative, and of course are exactly evaluations at  points of $R = I \times J$ ($R$ stands for rectangle, but neither $I$ nor $J$ need be intervals). For $z = (\alpha,\beta)$ in $ I \times J$, let $\phi_{\alpha,\beta}$ be point evaluation at $z$; of course, this corresponds to $\sigma_{\alpha} \otimes \tau_{\beta}$ after our identification. 

We can ask when $\phi_{\alpha,\beta}$ is good. Since we are dealing with Bauer simplices here (taking the constant functions $\pmb 1$ as the order units for the respective dimension groups), the affine function spaces are respectively the continuous function rings  $C(I,\R)$, $C(J, \R)$ and $C(R,\R)$ equipped with the supremum norm, and $\brcs{\sigma_{\alpha}}^{\perp}$ is just the ideal of $C(I,\R)$ consisting of continuous functions that vanish at $\alpha$, which we denote $\Ann (\alpha)$
(and similarly for $\brcs{\tau_{\beta}}^{\perp}$ and $\brcs{\phi_{\alpha,\beta}}^{\perp}$). 

We have already determined when $\sigma_{\alpha}$ (respectively $\tau_{\beta}$) is a good trace: this occurs if and only if $\alpha$ is algebraic and no other roots of its (monic) irreducible polynomial over the rationals, $p_{\alpha}$, belong to $I$. 

From the earlier result, for $(\alpha, \beta)$ in $R$, $\phi_{\alpha,\beta}$ is good if both $\sigma_{\alpha}$ and $\tau_{\beta}$ are good. We want to analyze the general case for points of $R$, but we do not obtain complete results. Nonetheless, interesting phenomena are going on. The rings $\Q[\alpha]$, $\Q[\beta]$, and $\Q[\alpha,\beta]$ that appear are the rational subalgebras of the reals generated by the indicated elements. If both $\alpha$ and $\beta$ are algebraic, all three rings are subfields of $\R$. Of course, $\sigma_{\alpha} (G) = \Q[\alpha]$, $\tau_{\beta} (H) = \Q[\beta]$, and $\phi_{\alpha,\beta} (\Q[x,y]) = \Q[\alpha,\beta]$.

\itemitem{$\bullet$}{\it If $\phi_{\alpha,\beta}$ is good and $R$ is a rectangle (so $I$ and $J$ are intervals), {\bf and\/} $(\alpha,\beta)$ belongs to the interior of $R$, then $\alpha$ and $\beta$   are algebraic over $\Q$.} 

\Pf  If $\brcs{\alpha,\beta}$ is algebraically independent, then $\ker \phi_{\alpha, \beta}$ is $\brcs{0}$, contradicting $\phi_{\alpha,\beta}$ being good. By interchanging $x$ with $y$ if necessary, we may assume that $\alpha$ is transcendental, and $\beta$ is algebraic over the field $\Q(\alpha)$.

Let $f(x,y) = \sum_{j} y^j g_j(x)$ be the polynomial of minimal degree in $y$ with $g_j \in \Q(x)$ \st $f(\alpha,\beta) = 0$. Then $\ker \phi_{\alpha,\beta} \subseteq f\Q(x)[y]$. By multiplying by a suitable polynomial in $x$ (the product of the denominators of the $g_j$), we may assume that $g_j$ themselves are  in $\Q[x]$, and moreover, since $\alpha$ is transcendental, this does not induce any new zeros on $(\alpha,y)$. (We subsequently use the fact that since $\alpha$ is transcendental, none of the denominators can vanish at $\alpha$, hence there exists a neighbourhood of $\alpha$ on which none of the denominators hit zero.)

We see that $f_y:= D_y f = \sum jy^{j-1} q_j(x)$, so by minimality of degree in $y$, $f_y (\alpha,\beta) \neq 0$. By the implicit function theorem [the first time either of us can recall using it!], there exists a continuous (actually analytic, but we don't require that) function, $(x(t), y(t))$ on a neighbourhood of $\beta$ with $x(\beta) = \alpha$ \st $f(x(t),y(t))  $   vanishes thereon (the leading term does not vanish on a neighbourhood of $\alpha$ either, take the intersection). Since $(\alpha,\beta)$ is in the interior of $R$, at least one point on the path other than $(\alpha,\beta)$, call it $(\alpha',\beta')$ will belong to $R$. Then $\ker \phi_{\alpha, \beta}$ will be killed by $\phi_{\alpha',\beta'}$, since $f$ vanishes along the curve. Hence $\ker \phi_{\alpha,\beta}$ is not good.
\qed

This is a simplified version of the higher dimensional case, discussed later. 

What happens on the boundary? As our colleague Damien Roy pointed out, the boundary  of $I \times J$ (if $I$ and $J$ are intervals) can behave miserably. For example, if the upper left corner of $R$ is located at $(e,e)$ (where $e$ is transcendental), then $\ker \phi_{e,e} = (x-y)\Q[x,y]$, and since the line $x=y$ hits $R$ only at $(e,e)$, it follows that the closure of $\ker \phi_{e,e} $ in $C(R,\R)$ is $(e,e)^{\perp}$, hence $\phi_{e,e}$ is good. This is the first example wherein both $\sigma$ and $\tau$ have zero kernel (are bad), but $\sigma \otimes \tau$ is good. This can also occur on an edge, e.g., if $(e,e^2)$ is a point on the upper horizontal edge of $R$, then $\ker \phi_{e,e^2}$Êis $(x^2-y)\Q[x,y]$, and the zero set of $x^2-y$ has intersection with  $R$ consisting only of the point $(e,e^2)$, and thus $\phi_{e,e^2}$ is good. 

In contrast, $\phi_{\alpha,\beta}$ is bad if and only if $\brcs{\alpha,\beta}$ is algebraically independent.

Because of this pathology on the boundary, we  assume from now on that $\alpha$ and $\beta$ are algebraic over $\Q$.

\itemitem{$\bullet$}{\it If $\Q[\alpha] \subseteq \Q[\beta]$ and $\tau_{\beta}$ is good, then $\phi_{\alpha,\beta}$ is good.} This is an immediate consequence of \fivthr. 

In particular, if $\Q[\alpha]= \Q[\beta]$, then sufficient for $\phi_{\alpha,\beta}$ to be good is that either $\sigma_{\alpha}$ or $\tau_{\beta}$ be good. We will improve  on this, but before doing so, obtain an easy result. 

If $K$ and $L$ are subfields of a field $M$ with each containing a subfield $N$, we say they are {\it disjoint\/} (\wrt $N$) if $K \otimes_N L$ is naturally isomorphic to $KL$, the field they generate inside $M$; equivalently, $\dim_N KL = (\dim_N K)(\dim_N L) $ (and there are lots of other equivalent formulations, for example, $K \otimes_N L$ is a field). In our case, we will restrict to $N = \Q$, $K = \Q[\alpha]$, and $L = \Q[\beta]$; then $KL = \Q[\alpha,\beta]$, and disjointness occurs if and only if $p_{\beta}$ (the irreducible polynomial of $\beta$ over $\Q$) is irreducible over $\Q[\alpha]$. This is stronger than merely $K \cap L = \Q$ (in fact, $K$ and $L$ could be unequal but isomorphic, in which case their tensor product will not be a field). Disjointness is very closely related to the property discussed in Proposition  \fivfiv.

\itemitem{$\bullet$} {\it Assume that  $\Q[\alpha]$ and $\Q[\beta]$ are disjoint; then $\phi_{\alpha,\beta}$ is good if and only if both $\sigma_{\alpha}$ and $\tau_{\beta}$
 are good.}

\Pf Let $\Cal I$ be the ideal of $\Q[x,y]$ generated by $\brcs{p_{\alpha}(x), p_{\beta}(y)}$ (that is, the polynomials in $\Q[x]$ and $\Q[y]$ that are the monic irreducible polynomials of $\alpha$, respectively $\beta$, over $\Q$). Then $\Cal I$ is contained in $\ker \phi_{\alpha, \beta}$, and obviously $\Q[x,y]/\Cal I$ has dimension $(\deg p_{\alpha})(\deg p_{\beta})$; however, this is the same as $\dim_{\Q}\Q[\alpha,\beta]$, the image of  $\phi_{\alpha,\beta}$; thus, $\Cal I= \ker \phi_{\alpha,\beta}$.

If $\phi_{\alpha}$ were not a good trace of $G = \Q[x]$, then $p_{\alpha}$ would have another root, $\alpha'$, in $\Cal I$. Then $\phi_{\alpha',\beta}$ kills $\Cal I$, and therefore $\Cal I$ cannot be dense in $(\alpha,\beta)^{\perp}$, so that $\phi_{\alpha} $ is not good. Similarly, $p_{\beta}$ not good implies $\phi_{\alpha,\beta}$ is not good. The converse is a consequence of \fivtwo. \qed

\itemitem{$\bullet$} {\it General criterion.} If $\alpha$ and $\beta$ are both algebraic over $\Q$, and $p(x,y)$ is such that  $p(\alpha,y)$ is the monic irreducible polynomial of $\beta$ over $\Q[\alpha]$, then $\ker \phi_{\alpha,\beta} = (p_{\alpha}(x), p(x,y))$. In particular, $\phi_{\alpha,\beta}$ is good if and only if for every algebraic conjugate $\alpha'$ of $\alpha$ in $I$, and every algebraic conjugate $\beta'$ of $\beta$ over $\Q[\alpha]$ in $J$, $p(\alpha',\beta') = 0$ implies $(\alpha',\beta') = (\alpha,\beta)$.

This is a consequence of a more general result (obtained later) about generating sets for annihilating ideals, but is too difficult to work with, so we consider  a number of special cases. 
\comment
First, a modest generalization of Galois extension of a field. Let $K \subseteq L \subset \R$ be subfields, with $L$ finite dimensional over $K$. We say $L$ over $K$ is {\it $\rho$-Galois\/} if for every homomorphism $\Arrow \psi;L.\R$ that fixes $K$ pointwise, $\sigma(L) = L$; equivalently, $L$ is the unique maximal formally real subfield of the Galois closure of $L$ over $K$.

Non-Galois examples include everything for which there is just one real embedding, e.g., $K = \Q$, $L = \Q[m^{1/n}]$ with $m$ a square-free integer and $n$ odd. Non-examples include any extension field of $K = \Q$ that is not Galois, but for which its Galois closure is embeddable in the reals (e.g., $L$ is a cubic extension with three real roots, for which the Galois group is $S_3$).

\endcomment

Now we give sufficient conditions for  $\phi_{\alpha,\beta}$ to be  good; these are also necessary when $\Q[\alpha] \subseteq \Q[\beta]$ (still assuming they are algebraic). We will see interesting phenomena after we cut through the technicalities.

If $f$ is a polynomial in one variable (either $x$ or $y$), $Z(f)$ will denote the set of its roots in $\C$. Suppose that $K:= \Q[\alpha]\cap \Q[\beta]$ is nontrivial, that is, the intersection is not the rationals. By the primitive element theorem, there exists $\theta$ in $K$ \st $K = \Q[\theta]$. We can express $\theta$ uniquely in each of the forms, $\theta = \sum_{i=0}^s q_i \alpha^i$ and $\theta = \sum_{j=0}^t r_j \beta^j$ where $q_i$ and $r_j$ are all rational, $q_s r_t \neq 0$, $s < \deg p_{\alpha}$, and $t < \deg p_{\beta}$. Set $f \equiv f_{\theta} = \sum q_i x^i$ and $h \equiv h_{\theta} = \sum r_j y^j$. Then $f \in \Q[x]\subset \Q[x,y]$, $h \in \Q[y] \subset \Q[x,y]$, and $f(\alpha) = \theta = h(\beta)$. Let $\Cal I$ denote the ideal of $\Q[x,y]$ generated by $\brcs{p_{\alpha(x)}, p_{\beta}(y), f(x)-h(y)}$ ($\Cal I$ sometimes depends on the choice of $\theta$). Then $\Cal I$ is contained in $\ker \phi_{\alpha,\beta}$, frequently equalling it. 

\Lem Lemma \fivsix. Sufficient for $\Cal I = (p_{\alpha}(x), p_{\beta}(y), f(x)-h(y))$ to be dense in $(\alpha,\beta)^{\perp}$ is  the following set of conditions:
{\par}
\item{(i)} $f^{-1}(f(\alpha)) \cap Z(p_{\alpha}) \cap I = \brcs{\alpha}$;
\item{(ii)} $h^{-1}(h(\beta)) \cap Z(p_{\beta}) \cap J = \brcs{\beta}$;
\item{(iii)} $f(Z(p_{\alpha}) \cap I) \cap h(Z(p_{\beta}) \cap J) = \brcs{\theta} = \brcs{f(\alpha)}$. {\par}
\noindent These guarantee that $\phi_{\alpha,\beta}$ is good. {\par}
\noindent When $\Cal I = \ker \phi_{\alpha,\beta}$, these conditions are also necessary for $\phi_{\alpha,\beta}$ to be good.

\Remark All the conditions refer to the behaviour of finite sets (the zeros of $p_{\alpha} $ in $I$ and the zeros of $p_{\beta}$ in $J$) under $f$ and $h$ respectively. They simplify drastically in special cases. 

\Pf We shall prove density of $\Cal I$ in $(\alpha,\beta)^{\perp}$; this will force density of $\ker \phi_{\alpha,\beta}$, which in turn yields that  $\phi_{\alpha,\beta} $ is good.  Let $\overline{\Cal I}$ denote the closure of $\Cal I$ in $C(R,\R)$. Since $\Cal I$ is an ideal of $\Q[x,y]$ contained in $\ker \phi_{\alpha,\beta}$, $\overline{\Cal I}$ is a closed ideal of $C(R,\R)$, obviously contained in $(\alpha,\beta)^{\perp}$. Hence to show equality, it suffices to show that for all points $v$ of $R \setminus{(\alpha,\beta)}$, evaluation at $v$ does not kill $\Cal I$. 

Write $v = (\gamma,\delta) \in R = I \times J$, and suppose $\phi_{\gamma,\delta} (\Cal I) = \brcs{0}$. Evaluating the generators, $p_{\alpha}(x)$, $p_{\beta}(y)$, and $f(x)-h(y)$ at $v$, we obtain (a) $\gamma \in Z(p_{\alpha}) \cap I$, (b) $\delta \in Z(p_{\beta}) \cap J$, and  (c) $f(\gamma) = h(\delta)$. 

If $\gamma \neq \alpha$ but $\delta = \beta$, we deduce $f(\gamma)= f(\alpha)$, contradicting (i). If $\gamma = \alpha$ and $\delta \neq \beta$, then $h(\delta) = h(\beta)$, contradicting (ii). 

If $\gamma \neq \alpha$ and $\delta \neq \beta$, we note by (i) and (ii), that $f(\gamma) \neq f(\alpha) = h(\beta) \neq h(\delta)$, and $f(\gamma) = h(\delta)$ contradicts (iii). Hence $\Cal I$ is dense in $\phi_{\alpha,\beta}$.

Now suppose $\Cal I = \ker \phi_{\alpha,\beta}$ and at least one of (i-iii) is violated. We will find $v$ in $R\setminus \brcs{(\alpha,\beta)}$ \st evaluation at $v$ kills 
$\Cal I$, hence $\ker \phi_{\alpha,\beta}$ is not dense, so by \onesix, $\phi_{\alpha,\beta}$ is not good. If (i) fails, there exists $\alpha' $ in $Z(p_{\alpha}) \cap I$ \st $\alpha' \neq \alpha$ and $f(\alpha') = f(\alpha) = h(\beta)$. Set $v = (\alpha',\beta) \in I \times J$; obviously, $\phi_{\alpha',\beta}(p_{\alpha}(x)) = \phi_{\alpha',\beta}(p_{\beta}(y)) = 0$, and $\phi_{\alpha',\beta}(f(x) -h(y)) = f(\alpha')-h(\beta) = 0$, so $\phi_{\alpha',\beta}$ kills $\Cal I$.

If (ii) fails, there exists $\beta'$ in $Z(p_{\beta}) \cap J$, unequal to $\beta$, \st $h(\beta') = h(\beta)$; set $v = (\alpha,\beta') \in I \times J$, and as in the previous case, we see that $\phi_{\alpha,\beta'}$ kills $\Cal I$.

If (iii) fails, there exist $\alpha'$ in $Z(p_{\alpha}) \cap I$ and $\beta'$ in $Z(p_{\beta}) \cap J$ \st $f(\alpha') = h(\beta')$. Then $v = (\alpha',\beta')$ is the desired point.
\qed 

The ideal  $\Cal I$ depends on the choice of $\theta$, and it might be useful to examine the larger ideal obtained by considering all choices of $\theta$. However, still left unconsidered are the cases wherein $\Q[\alpha] \cap \Q[\beta] = \Q$ and the two fields are not disjoint. 

\Lem Corollary \fivsev. Suppose $\Q \neq \Q[\alpha] \subseteq \Q[\beta]$; with unique choice of $h$ of degree less than $\deg p_{\beta}$ \st $ \alpha = h(\beta)$, define  $\Cal I$ to be the $\Q[x,y]$-ideal generated by $\brcs{p_{\alpha}(x), p_{\beta}(y), x-h(y)}$. Then  $\Cal I = \ker \phi_{\alpha,\beta}$. Moreover, $\phi_{\alpha,\beta} $ is good if and only if for all $\beta' \neq \beta$ in $Z(p_{\beta}) \cap J$, (ii \& iii) $h(\beta') \not\in  Z(p_{\alpha}) \cap I$.

\Pf We observe that the codimension of $\Cal I$ is at most $\deg p_{\beta}$, but this is the dimension of $\Q[\beta] = \Q[\alpha,\beta] = \phi_{\alpha,\beta}(\Q[x,y])$; since $\Cal I \subseteq \ker \phi_{\alpha,\beta}$, equality holds.

Now condition (i) of the previous result is automatic, and (ii) and (iii) translate to the combined condition hypothesized. 
\qed

\comment

\endcomment

There is one more easy case: if $\alpha$ is rational (that is, $\Q[\alpha] = \Q$); in particular, $\sigma_{\alpha}$ is good. Then  $\ker \phi_{\alpha,\beta}$ is the ideal generated by $\brcs{x-\alpha,p_{\beta}(y)}$. Now it is obvious that $\phi_{\alpha,\beta}$ is good if and only if $\phi_{\beta}$ is. This is also covered by the disjointness result.

\Lem Corollary \fiveig. (a) Suppose that $\alpha = \beta^k$  for some  positive integer $k < \deg p_{\beta}$. Then $\phi_{\alpha,\beta}$ is good if and only if  no $k$th power of a root in $J$ of $p_{\beta}$ other than $\beta$  equals a root of $p_{\alpha}$ in $I$, and if $k$ is even, $-\beta$ is not in $J \cap Z(p_{\beta})$. {\par}
\noindent (b) If $\alpha \in I \cap J$, then $\phi_{\alpha,\alpha}$ is good if and only if $I \cap J \cap Z(p_{\alpha}) =\brcs{\alpha}$. {\par}
\noindent (c) If $\alpha/\beta = q \in \Q\setminus\brcs{0}$, then $\phi_{\alpha,\beta}$ is good if and only if $q^{-1}I \cap J \cap Z(p_{\beta}) = \brcs{\beta}$. 

\Pf We verify conditions (ii \& iii) of  Corollary \fivsev. (a) Here $h(y) = y^k$, so $h(y) = \alpha$ has only nonreal roots (other than $\beta$) if $k$ is odd, otherwise it has $-\beta$ as a zero. (b)
Here $h(y) = y$. (c) Here $h(y) = qy$.\qed

This allows us to construct lots of examples. For instance (always assuming $\alpha$ is algebraic), suppose $\alpha \in I \cap J$, and none of its other algebraic conjugates belong to $I \cap J$, but $Z(p_{\alpha}) \cap I$ and $Z(p_\alpha) \cap J$ both contain at least two elements (one of which is $\alpha$). Then neither  $\sigma_{\alpha}$ nor $\tau_{\alpha}$ is good, but their tensor product, $\phi_{\alpha,\alpha}$, is good. Specifically, we can arrange this by setting $I = [a,b]$, $J = [c,d]$, and $p$ an irreducible rational polynomial with at least three real roots, $\alpha_2 < \alpha < \alpha_1$  with $\alpha_2 < a < \alpha < \alpha_1 < b$ and $c < \alpha_2 < \alpha < d < \alpha_1$ \st all other real roots are either bigger than $b$ or less than $c$. 

For example, the irreducible polynomial of $\alpha=\sqrt 2 + \sqrt 3$ is degree four, and the four roots, all real, are $\pm \sqrt 2 \pm \sqrt 3$. By varying the locations of $a \leq \alpha \leq b$ and $c\leq \alpha \leq d$, various phenomena can be observed.

\Lem Corollary \fivnin. Suppose $\tau_{\beta}$ is good. Then $\phi_{\alpha,\beta}$ is good if $f(\alpha') \neq f(\alpha)$ for all $\alpha' \neq \alpha$ in $Z(p_{\alpha}) \cap I$. 

\Pf We have already seen that $\tau_{\beta}$ is good if and only if $Z(p_{\beta}) \cap J$ is the singleton $\brcs{\beta}$.  Thus conditions (ii) and (iii) are automatically satisfied, and the hypothesized condition is a reformulation of (i). \qed

This example suggests more examples. Let $M$ be a compact subset of $\R^n$, preferably regular (equalling the closure of its interior), but not necessarily. Form $A = \Q[x_1, \dots, x_n]$; viewed as a ring of functions on $M$, $A$ will be norm dense in $C(M,\R)$, so equipped with the strict ordering, $A$ is a simple dimension group and also an ordered ring. For $\alpha = (\alpha_1, \dots, \alpha_n)$ in $M$, define the point evaluation at $\alpha$, $\phi_{\alpha}$; as before, this is a pure trace of $A$, and they are all of this form.

Then $\phi_{\alpha}$ will be a good trace if and only if the zero set of $\ker \phi_{\alpha}$ contains no points of $M$ other than $\alpha$ (this comes from the density criterion). If we use the notation $I_{\Q}(S)$ for the ideal of $\Q[x_1, \dots, x_n]$ consisting of those elements that vanish at a subset $S$ of $M$ (of course, if $F$ is a subfield of the complexes, then we can define $I_{F}$ similarly, as an ideal in $F[x_1, \dots, x_n]$), then $\phi_{\alpha}$ is good if and only if
$$
Z(I_{\Q} (\alpha)) \cap M = \brcs{\alpha}.
$$
(recall that for $U$ a set of polynomials, $Z(U)$  denotes  the set of complex zeros common to all elements of $U$).

Modulo the usual obnoxious boundary behaviour, we can start things off with a consequence of the implicit function theorem, again.

\Lem Proposition \fivten. Let $M$ be a compact subset of $\R^n$, and let $A = \Q[x_1, \dots,x_n]$ equipped with the strict ordering as functions on $M$. Suppose that $\alpha = (\alpha_1, \dots, \alpha_n)$ is in the interior of $M$. If the pure trace of $A$, $\phi_{\alpha}$ (evaluation at $\alpha$) is  good, then all $\alpha_i$ are algebraic over the rationals.

This requires the following  presumably well known result on generating sets for annihilator ideals.

\Lem Proposition \fivele. 
Let $K$ be a subfield of $\R$, and let $b= (\beta_1, \dots, \beta_s)$ be a point in $\R^s$ all of whose entries are algebraic over $K$. There exist polynomials $\brcs{f_i}_{i=1}^s$ of the polynomial ring $K[y_1, \dots, y_s]$ with the following properties.{\par}
\item{(i)} $f_i \in K[y_1, y_2, \dots, y_i]$ for all $1 \leq i \leq s$;
\item{(ii)} $\frac{\partial f_i}{\partial y_i}(b) \neq 0$ for all $1 \leq i \leq s$;
\item{(iii)} $I_K(b)$ is the ideal of $K[y_1, \dots, y_s]$ generated by $\brcs{f_1, \dots, f_s}$. 

\Pf
Define the following polynomials, $F_i \in K[\beta_1, \dots, \beta_{i-1}][y_i]$ (so $F_1 \in K[y_1]$), via
$F_i$ is the monic irreducible polynomial, in $y_i$, of $\beta_i$ over $K[\beta_1, \dots, \beta_{i-1}]$.  Now define $f_i \in K[y_1, \dots,y_i]$ by replacing every occurrence of $\beta_j$ in $F_i$ by $y_j$ (for every $j < i$); we may also assume (by using the algebraicity of $\beta_j$), that the degree of $y_j$ in $f_i$ is at most $K[\beta_1, \dots, \beta_j]: K[\beta_1, \dots, \beta_{j-1}]-1$.  (This last is not necessary for the argument here, but would be useful if we made it into an algorithm.) Then $f_i(\beta_1,\dots,\beta_{i-1},y_i) = F_i(y_i)$

Each $f_i$ is in   $K[y_1, \dots, y_i] \subseteq K[y_1, \dots, y_s]$, so we can evaluate them at points of $\R^s$. Let $b^{(i)}$ be the truncation to the first $i$ coordinates. Then $f_i (b) = F_i (b^{(i)}) = 0$. In particular, $f_i \in I_K (b)$.

Let $h \in K[y_1, \dots, y_s]$ be \st $h(b) = 0$. Then $h_s(y_s):=h(b^{(i-1)},y_s)$ vanishes at $y_s = \beta_s$, so that  as elements of $K[\beta_1, \dots, \beta_{s-1}][y_s]$, $F_s$ divides $h_s$. That is, there exists $G_s \in K[\beta_1, \dots,\beta_{s-1}]$ \st $h_s = F_s G_s$. We may form $g_s$ in $K[y_1, \dots, y_s]$ as before, so that $G_s (y_s) = g_s(b^{(s-1)}, y_s)$.

Now form $h^{(s-1)} = h_{s-1}(y_{s-1}, y_s)-f_{s}(b^{s-2},y_{s-1},y_{s-1}, y_s)g_{s}(b^{(s-2)},y_{s-1}, y_s) $. We can rewrite this as $\sum y_s^l \sum c_l (y_{s-1})$ where $c_l$ belongs to $K[\beta_1, \dots, \beta_{s-2}][y_{s-1}]$. Plugging in $y_{s-1} = \beta_{s-1}$, we obtain that each $c_l (\beta_{s-1})$ must be zero (since the identity is true for all values of $y_s$), hence $F_{s-1}$ divides each one of those. Thus we have $h^{(s-1)}(y_{s-1}, y_s) = F_{s-1} G_{s-1}$, where $G_{s-1}$ belongs to $K[\beta_1, \dots, \beta_{s-2}, y_{s-1}, y_{s-2}]$. As before, we can find $g_{s-1}$ with coefficients in $K[y_1, \dots, y_s]$ \st $g_{s-1} (b^{s-2},y_{s-2}, y_{s-1}) = G_{s-1}$. Then $h(b^{(s-2)}, y_{s-1},y_s)- f_{s}(b^{(s-2)},y_{s-1},y_{s-1}, y_s)g_{s}(b^{(s-2)},y_{s-1}, y_s) -
f_{s-1}(b^{(s-2)},y_{s-1},  y_s)g_{s-1}(b^{(s-2)},y_{s-1}, y_s) $ vanishes identically.

We can continue this process, obtaining successively $g_{s-j} \in K[y_1, \dots, y_s]$   ($j = 0,1,2, \dots$) \st $h(b^{(s-j-1)}, y_{s-j}, \dots, y_{s})-\sum_{k=0}^{j-1} f_{s-k}(b^{(s-k-1)}, y_{s-k}, \dots, y_{s-j}) g_{s-k}(b^{(s-k-1)}, y_{s-k}, \dots, y_s)$ vanishes identically. This continues until we run out of coordinates, and we obtain $h = \sum f_{s-j} g_{s-j}$. This yields (iii); (i) is by construction, and (ii) follows from the fact that $\frac{dF_i}{dy_i}(\beta_i) \neq 0$, as the irreducible polynomial, $F_i(y_i)$, has no repeated roots.
\qed

{\noindent \it Proof} (of Proposition \fivten) If $\brcs{\alpha_1, \dots, \alpha_n}$ is algebraically independent, then $\ker \phi_{\alpha} = \brcs{0}$, so $\phi_{\alpha}$ could not be good. Hence after relabelling, we may assume that for some $r $ with $1 \leq r <n$, the set $\brcs{\alpha_1, \dots,\alpha_r}$ is algebraically independent (over $\Q$), and each of $\beta_k:= \alpha_{r+k}$ ($k> 0$) is algebraic over the function field, $K = \Q(\alpha_1, \dots, \alpha_r)$. Let $a = (\alpha_1, \dots, \alpha_r)$. We may find the polynomials (with coefficients from $K$), $f_i$ satisfying the conditions above (where $y_i:= x_{r+i}$). We can thus write each  
$$
f_i:= \sum y^w \frac{A_{wi}}{B_{wi}}
$$
where $w$ is an $i$-tuple in $(\Z^+)^i$ (and $y^w$ is the corresponding monomial in the $y\,$s),  $A_{wi}$ and $B_{wi}$ are elements of $\Q[\alpha_1, \dots,\alpha_r]$ \st none of $B_{wi}$ are zero, and at least one of $A_{wi}$ is not zero. Because $a$ is algebraically independent, there exist unique $a_{wi}$, $b_{wi}$ in $\Q[x_1, \dots, x_r]$ \st $a_{wi}(a) = A_{wi}$ and $b_{wi} (a) = B_{wi}$; set $f_i^{0}$ to be the resulting elements of  $\Q(x_1,\dots, x_r)[y_1, \dots y_s]$. Set $q = \prod_{w,i} b_{wi}$. Since $a$ is algebraically independent, $q(a) \neq 0$. Hence there exists a neighbourhood, $N$, of $a$ in $\R^r$ \st $q$ has no zeros in $N \times \R^{n-r}$; of course, each  $qf_i^{0}$  belongs to $\Q[x_1, \dots, x_r; y_1, \dots, y_i]$. It easily follows that if $h \in \Q[x_1, \dots, x_r, y_1, \dots, y_s]$ and $h(a,b) = 0$, then there exists $p \in \Q[x_1, \dots, x_r]$, with no zeros on $N$, \st $ph$ is in the ideal of $\Q[x_1, \dots, x_r, y_1, \dots, y_s]$ generated by $\brcs{qf_i^0}$.  

Define the map $\Arrow(qf_i^0);N\times \R^s.\R^s$  by $v \mapsto (qf_1^0(v), qf_2^0(v), \dots, qf_s^0(v))$. This is analytic on $N\times R^{s}$ and  sends $(a,b)$ to $ \pmb 0$. 

Calculating the Jacobian matrix of the map,  $J = \( \smallmatrix \frac{\partial qf_i^0}{\partial y_{j}} \endsmallmatrix \)$, we see that $J$ is an upper diagonal  matrix, whose entries on the diagonal do not vanish at $(a,b)$. Hence $J(a,B)$ is invertible, so we may apply the implicit function theorem. There exist an open neighbourhood $V \subseteq N$ in $\R^r$ of $a$, an open neighbourhood $W$ of $b$ in $\R^s$, together with an analytic function $\Arrow g;V.W$ with  $g(a) = b$ and  \st the graph of $g$ is $\Set{(v,w ) \in V \times W}{(qf_i^0)(v,w) = \pmb 0}$.

Since $\alpha$ is in the interior of $M$, there exists a point in $M$ of the form $(a',b'):= (\alpha_1', \dots, \alpha_r',g(\alpha_1', \dots, \alpha_r'))$ unequal to $(a,b)$.   Hence $f_i^0(a',b') = 0$ for all $i$. If $h$ in $\Q[x_i, y_j]$ vanishes  at $(a,b)$, we saw previously that $ph$ vanishes where all the $f_i^0$ do, and since $(a',b')$ is in $N \times \R^{s}$, $p$ vanishes nowhere on this neighbourhood, and thus $h(a',b') = 0$. In particular, $\phi_{a',b'}$ kills $\ker \phi_{a,b}$, hence $\phi_{a,b}$ is not good. 
\qed

The argument proves a bit more; since there is at least a curve of zeros in the neighbourhood,  $Z(I_{\Q} (\alpha)) \cap M$ contains a continuum when not all the entries of $\alpha$ are algebraic and $\alpha$ is in the interior of $M$---so there are uncountably many pure traces that kill $\ker \phi_{\alpha}$. 

One consequence of the existence of the generating set is that if $\alpha = (\alpha_i)$ is an algebraic point, then the closure of $\ker \phi_{\alpha}$ in $C(M,\R)$ is of codimension at most $t:=\dim_{\Q} \Q[\alpha_1, \dots, \alpha_n] -1$ in the ideal $\Ann (\alpha)$ of $C(M,\R)$. Another aspect is that we need only check the generating set on at most $t$ points to decide whether $\phi_{\alpha}$ is good. 

The problem with the boundary would go away in the weird case that for every non-algebraic boundary point of $M$, there are no analytic curves whose intersection with $M$ consists locally exactly of the boundary point.

\SecT  6 Impure traces

Some  results  on goodness require $\tau$ to be pure. If $\tau$ is not pure, then it might still be worthwhile to study goodness, because it requires some relations among the pure traces in the support.

\Lem Lemma \sixtwo. Let $G$ be an approximately divisible dimension group with order unit, and let $\tau_1$ and $\tau_2$ be pure normalized traces. Let $\lambda$ be a real number in the open unit interval, and define $\tau = \lambda \tau_1 + (1-\lambda)\tau_2$. If
$\lambda \tau_1 (G) \cap (1-\lambda) \tau_2 (G) = \brcs{0}$, then  $\tau$ is not good.

\Remark If $G$ is countable, then the set of $\lambda$ for which $\lambda \tau_1 (G) \cap (1-\lambda) \tau_2 (G) = \brcs{0}$ is cocountable, so that goodness hardly ever occurs for these combinations.

\Remark As is clear from the proof, $\tau_i$ need not be pure, but simply different enough that they can be separated by elements of the affine function space.

\Pf Select positive real numbers $r$ and $s$ \st $r > \lambda s/(1-\lambda)$. There exists $h$ in the affine function space \st $h(\tau_1) = -s$ and $h(\tau_2) = r$. For all sufficiently small $\epsilon$, a perturbation of $r$ and $s$ independently by at most $\epsilon$ will guarantee that the perturbed $-s$ is negative, the perturbed $r$ is still positive, and the perturbed $r$ and $s$ will still satisfy the inequality $r' > \lambda s'/(1-\lambda)$. As $G$ is approximately divisible, we may find $g$ in $G$ \st $\|\hat g-h\| < \epsilon$. Then $|\tau_1(g) + s|, |\tau_2(g)-r| < \epsilon$ (so $\tau_1(g) < 0$ and $\tau_2 (g) > 0$) and
$$
 \tau_2 (g) > -\frac{\lambda \tau_1 (g)}{1-\lambda}.\tag *
$$

We may write $g = b-a$ where $b$ and $a$ are in $G^+$. Then $\tau_2 (b) > \tau_2 (a)$, $\tau_1 (a) > \tau_1 (b)$ (so $\tau_1 (a) > 0$, and thus $\tau(a) > 0$), and most importantly, (*) ensures that $\tau(a) < \tau(b)$.  If $0 \leq c \leq b$ and $\tau(c) = \tau(a)$, the latter implies $0 = \lambda \tau_1 (c-a) + (1-\lambda) \tau_2(c-a)$. The hypothesis on the intersection forces $\tau_i(c -a) = 0$, whence $\tau_1 (c) = \tau_1 (a) > \tau_1 (b)$, a contradiction.
\qed

Some convex combinations of pure traces are good, e.g., consider $G = \Q \oplus \Q$ with the strict ordering from the obvious inclusion in $\R^2$: a trace given as $\lambda \tau_1 + (1-\lambda) \tau_2$  with $ 0 \leq  \lambda \leq 1$, a weighted average of the two coordinates, is good, if and only if $\lambda$ is rational.

On the other hand, the following trivial result shows that a lot more convex combinations, even with rational choices for $\lambda$, are ruled out.

\Lem Lemma \sixthr. Let $G = \lim G_{\alpha}$ be a dimension group expressed as a direct limit of ordered abelian groups over a direct set, with maps $\Arrow \phi_{\alpha}; G_{\alpha}.G$. Let $\tau$ be a trace of $G$ \st for all $\alpha$, $\tau \circ \phi_{\alpha}$ is a good trace of $G_{\alpha}$. Then $\tau$ is good.

\Pf Pick $a$ and $b$ in $G^+$ \st $0 < \tau(a) < \tau(b)$. There exists (by directedness) $\alpha$ \st $\phi_{\alpha}^{-1} (a)\cap G_{\alpha}^+$ and $\phi_{\alpha}^{-1} (b)\cap G_{\alpha}^+$ are both nonempty; choose one from each, say $A$ and $B$ respectively. As $\tau\circ \phi_{\alpha}$ is good, there exists $0 \leq C \leq B$ in $G_{\alpha}^+$ \st $\tau\circ \phi_{\alpha} (C) = \tau\circ \phi_{\alpha} (A)$. Then $c:= \phi_{\alpha}(C)$ is the desired element. \qed

If $G$ is a dimension group, then of course it is a direct limit with the $G_i$ being simplicial ordered groups ($\Z^{n(i)}$ with coordinatewise ordering), so if $\tau$ is a trace on $G$, it induces a trace on each of the simplicial ones, necessarily a positive linear  combination of the coordinates. Thus if $\tau$ is not good, infinitely many of those traces must themselves be not good (in fact, it is quite easy to find lots of rational combinations of pure traces on simplicial groups that are not good).

\SecT 7 Goodness en masse and refinability

Suppose $(G,u)$ is a dimension group with order unit, and $U$ is a set of normalized traces, that is, a subset of $S(G,u)$ (in view of subsequent definitions, we can sometimes assume that $U$ is closed,  compact, and convex). We say (the set of traces) {\it $U$ is weakly
good\/} if whenever $a$ and $b$ are in $G^+$ and $\inf_{\tau \in U}(\tau(b)-\tau (a)) > 0$, then there exists $c \in [0,b]$ \st $\tau(c-a) = 0$ for all $\tau$ in $U$. We say {\it $U$ is weakly order unit good\/} if we restrict $b$ to be order units of $G$, and we drop the weakly from the definition if $a$ is allowed to vary over all of $G$, not just $G^+$.

Other definitions are possible, and possibly more interesting, e.g., if we only require $\tau(b)-\tau(a) > 0$ for all $\tau$ in not necessarily compact $U$, or even $\tau(b) \geq \tau(a) $ for all $\tau$ in $U$. Even in the simple case, these notions are apparently different.

For example, if $U = S(G,u)$, goodness is automatic (a consequence of unperforation; note that none of the individual traces need be good themselves). If $U$ is a singleton, it reverts to the definition of good (or order unit good) that we previously used. The proofs of some  results carry over to the general case.

There is a generally smaller compact convex set than $Z(\ker U)$ associated to $U$, ${\tilde U} = L(U) \cap S(G,u)$, where $L(U)$ is the closure of the affine span of $U$ (we shall discuss this construction in more detail later). Then $U$ is good implies ${\tilde U}$ is good, and in this case, $Z(\ker U) = \tilde U$ (Lemma \appbsev). 

The notions of ugly  and bad apply to subsets $U$ of $S(G,u)$ as well: $U$ is {\it ugly\/} if the image of $\ker U$ is discrete in $\Aff S(G,u)$ and $U$ is good for the group $G \otimes \Q$; $U$ is {\it bad\/} if $\ker U = \Inf (G)$, equivalently, $Z(\ker U) = S(G,u)$.

\Lem Lemma \appbtwe. (c.f., Lemma \oneone) Let $(G,u)$ be a simple dimension group and let $U$ be a subset of $S(G,u)$.
\item{(a)} If for all $a,b \in [0,u]$ \st $0 \ll \hat a|U \ll \hat b|U$ there exists $c \in [0,b]$
\st $\hat a|U = \hat c|U$, then $U$ is good.

\Pf Virtually idential to the proof of Lemma \oneone(a).\qed

Let  $U$ be a subset of $S(G,u)$ and define $\ker U = \cap_{\tau \in U} \ker \tau$. This is a subgroup of $G$ (frequently zero); define $Z \equiv Z(\ker U) = \Set{\sigma\in S(G,u)}{\sigma(\ker U) = 0}$. Then $Z$ is compact, convex, and contains $U$ (it can be bigger than the closed convex hull of $U$), and in fact, $Z = S(G,u) \cap \Cal L$ where $\Cal L$ is a closed flat (the intersection of the zero sets of all elements of $\ker U$). Note that if for all $\tau$ in $U$ we have $\tau(b)-\tau(a)\geq \delta > 0$, then the same is true for all $\tau$ in $Z$ (thus goodness of $U$ implies that of $Z$). If $Z$ happens to be a face of $S(G,u)$, we will write it as $F$.  Let $U^{\vdash} =  \Set {h \in \text{Aff\,}S(G,u)}{h|U = 0}$. This is an order ideal in $\Aff S(G,u)$ if the closed convex hull of $U$ is a face of $S(G,u)$, and in that case we write $F^{\perp}$ instead of $U^{\vdash}$. 

Let $G = \Aff K$ for some Choquet simplex $K$, equipped with the strict ordering. Then we can ask which subsets $U$ of $K$ are good. There is a connection with a classical result on faces: $U$ is good (\wrt $G = \Aff K$), equivalently $\tilde U$ is good, if and only if for all $f,g \in (\Aff K)^{++}$,
$$
f|{\tilde U} \ll g|{\tilde U} \quad \text{implies there exists $h \in (\Aff K)^{++}$ \st $h \ll g$ and $h|{\tilde U}= f|{\tilde U}$.} 
$$
A well-known result in Choquet theory (e.g., [G; 11.25]) asserts that if $F$ is a compact convex subset of $K$, then $F$ is a face if and only if for all $f,g \in (\Aff K)^+$,
$$
f|F \leq g|F \quad \text{ implies there exists $h \in (\Aff K)^+$ \st $h \leq f,g$ and $h|F = f|F$.}
$$
The second property implies the first (there exists $\epsilon > 0$ \st $g' = g - \epsilon\pmb 1 \gg f$ and $f'= f -\epsilon\pmb 1 \gg 0$, and apply the second property to the pair $g',f'$, then add $\epsilon \pmb 1$ to the resulting $h$).  The  important difference between these two properties is the extra condition that $h \leq f$.
It forces the compact convex set $F$ to be a closed face; however, all singleton sets $U$ satisfy goodness (and these are not faces unless the singleton consists of an extreme point), and we will give constructions of more complicated compact convex sets that are  good but not faces. We will also show that $U$ is good for the simple dimension group $(G,u)$ if and only if $\ker U$ has dense image in ${\tilde U}^{\vdash} := \Set{h \in \Aff S(G,u)}{h|{\tilde U} \equiv 0}$ and ${\tilde U}$ is good for $\Aff S(G,u)$ with the strict ordering.

There is an elementary technicality which can lead to confusion. Let $(H,u)$  be a dimension group. Let $L$ be a subgroup of $H$ with the following properties:
\item{(i)} $L \cap H^+ = \brcs{0}$
\item{(ii)} for $b$ in $H$ and $n$ a positive integer, if  $nb \in L$,  then $b \in L$
\item{(iii)} if $b \in H$ and $nb \in L + H^+$ for some positive integer $n$, then $b \in L + H^+$

\noindent In applications, $L = \ker U$ for $U \subset S(G,u)$.
 
Form the quotient group $H/L$. Impose on it the pre-ordering given by $a + L \geq 0$ if there exists $a' \in H^+$ \st $a-a'\in L$. Condition (i) is equivalent to this pre-ordering being a partial ordering, that is, something that is both positive and negative must be zero. Condition (ii) asserts that the quotient group is torsion free, and condition (iii) is simply a restatement of $H/L$ being unperforated \wrt this ordering (there is some redundance, since unperforation implies torsion-freeness). Note that we are {\it not\/} assuming $L$ is an order ideal, indeed this would defeat the purpose. Typically, $L = \ker U$ where $U \subseteq S(H,u)$. Then (i) and (ii) hold automatically, but (iii) need not,
 even when $H$ is simple. Appendix B discusses sufficient and necessary conditions for $G/\ker U$ to be unperforated.

If we take the rank $5$ critical group $H = \langle e_i; \sum_{i=1}^4 \alpha _ie_i\rangle$ with the strict ordering inherited from $\R^4$, with $\brcs{1, \alpha_1,\alpha_2,\alpha_3,\alpha_4 }$ linearly independent over $\Q$, then $H$ is a simple dimension group, with pure traces given by the four coordinate maps, $\tau_i$. Let $L = (1,1,-1,-1)\Z$; this is the intersection of the kernels of the four traces $\phi_1 = \tau_1 + \tau_3$, $\phi_2 = \tau_1 + \tau_4$, $\phi_3 = \tau_2 + \tau_3$, and $\phi_4 = \tau_2 + \tau_4$. In particular $H/L$ is partially ordered via the quotient ordering. However, we note that $2(1,1,0,0) = (1,1,1,1) + (1,1,-1,-1)$, so $2(1,1,0,0) +L$ is in the positive cone of $H/L$. However, for no integer (positive or negative) is $(1,1,0,0) + m(1,1,-1,-1)$ in $H^+$. So $H/L$ is not 2-unperforated. Incidentally, the pure trace space of $H/L$ is $\brcs{\phi_i}$; since there is an affine relation $\phi_1 + \phi_4 = \phi_2 + \phi_3$, the trace space is not a simplex.
In this example, $U = \brcs{\phi_i}$ and $Z(\ker U) $ agrees with $\tilde U$ and is the convex hull of $\brcs{\phi_i}$. 
We will see more examples of this sort in Appendix B.

However, if (i) holds   and both $H$ and $L$ are rational vector spaces,  then  $H/L$ is torsion-free, and also $H/L$ is unperforated, that is (iii) holds.  The importance of this lies in the following elementary result.

\Lem Lemma \appbthr. Let $(H,u)$ be an approximately divisible  dimension group, and let $L$ be a subgroup of $H$ satisfying (i--iii) above. Let $h$ be an element of $H$. If for all traces $\psi$ of $(H,u)$ \st $\psi(L) = 0$, $\psi(h) > 0$, then there exists an order unit $v$ of $H$ \st $v - h \in L$. {\par}In particular, this applies if $L$ is of the form $\cap_{\tau \in U} \ker \tau$ for some subset $U \subset S(H,u)$ \st $\(\cap_{\tau \in U} \ker \tau\) \cap H^+ = \brcs{0}$ and $H$ is a rational vector space.

\Pf Since $H$  is approximately divisible,  its image in $\Aff S(H,u)$ is dense, and thus we may impose the strict ordering on $H$, $H^{++} \cup \brcs{0}$, that is take the order units of the original ordering (those with strictly positive image in $\Aff S(H,u)$ and make them together with $0$, the new positive cone. The new ordering on $H$ is also a dimension group ordering (since $S(G,u)$ is a Choquet simplex). The trace space is unchanged, and (i,ii,iii) remain in force, since $H^{++} \subseteq H^+ $. 

Impose the quotient ordering on $H/L$ arising from the strict ordering on $H$. Now we examine the normalized traces of $(H/L,u+L)$. If $\phi$ is a trace, then since $\Arrow \pi;H .H/L$ is by construction order preserving, $\phi\circ \pi$ is a trace of $H$, and obviously $\phi\circ \pi (L) = 0$. Conversely, if $\Arrow  \psi;H.\R$ is a normalized trace of $H$ \st $\phi (L) = 0$, then we can construct $\phi \in S(H/L,u+L)$ \st $\phi \circ \pi = \psi$ in the obvious fashion: since $\psi$ kills $L$, $\psi$ induces a group homomorphism $\Arrow \phi; H/L.\R$ defined as $\phi(a+L) = \psi(a)$; this is obviously well-defined, and moreover, if $a+L \geq 0$ (in $H/L$), then there exists $a' \geq 0$ (in $H$) \st $a- a' \geq 0$. Then $\psi(a) = \psi(a')\geq 0$, and thus $\phi$ is positive, i.e., a trace, and normalization is automatic. It is now routine that the map $S(H/L) \to S(H,u) \cap L^{\vdash}$ (the latter is the set of traces killing $L$)  given by $\phi \mapsto \phi\circ \pi$ is an affine homeomorphism. 

Finally, suppose $\psi(h) > 0$ for all $\psi$ in $S(H,u) \cap L^{\vdash}$. Then $\phi(h+L)> 0$ for all $\phi \in S(H/L,u+L)$, and since $H/L$ is unperforated, $h+ L$ is a positive element of $H/L$, and obviously nonzero. By definition of the quotient ordering, there exists $v$ in $H^{++} \cup \brcs{0}$ \st $v - h \in L$. Obviously $ v\neq 0$ (evaluate at any $\psi$), so $v \in H^{++}$, that is, $v$ is an order unit \wrt the original ordering on $H$. 

If $L = \cap_{\tau \in U} \ker \tau$ then (i) occurs by hypothesis, and (ii) is obvious, and when $H$ is a rational vector space, (iii) is a consequence of (i) and (ii).\qed

From a subset $U$ of $S(G,u)$, we can first form $ \ker U = \cap_{\tau \in U} \ker \tau$, a subgroup of $G$. More general conditions under which $G/\ker U$ is (or is not) unperforated are given in Appendix B. Next, we  form $Z(\ker U) = \Set{\sigma \in S(G,u)}{\sigma(\ker U) = 0}$. Then $U \subseteq Z(\ker U)$, and the latter is of the form $S(G,u) \cap W$ where $W$ is a closed flat (the intersection of the zero sets of all the elements of $\ker U$); it is also a closed convex subset of $S(G,u)$, but as we saw earlier, need not be a simplex itself. Finally, we form a closed convex real vector subspace  of $\Aff S(G,u)$, $Z(\ker U)^{\vdash} = \Set{h \in \Aff S(G,u)}{h|Z(\ker U) \equiv 0}$. The representation $g \mapsto \hat g$ sends $\ker U = \ker (Z(\ker U))$ to $Z(\ker U)^{\vdash}$.

It follows that for the rank 5 critical group given just before Lemma \appbthr, the set $U$ consisting of the four extreme traces of $G/L$ (or its convex hull, a square in $\R^4$) is not a good subset of $S(G \otimes \Q,u)$ (for any choice of order unit) even after tensoring $G$ with the rationals. On the other hand, all the individual traces are good \wrt $G \otimes \Q$.

A measure $\mu$ on $X$ is {\it refinable\/} if for all clopen $A$, whenever $\mu(A) = \sum \alpha_i$ where $\brcs{\alpha_i}$ is a finite set of elements in the  clopen value set ($S(\mu):= \Set{\mu(B)}{B \text{ is clopen}}$), there exist pairwise disjoint clopen sets $V_i$ \st $\mu(V_i) = \alpha_i $ and $A = \cup V_i$. Good measures are refinable, but the converse fails. Refinability is related to goodness en masse, as observed in Proposition \appbsix\ and Lemma \appbsev.

This definition translates directly to  dimension groups, which however, is somewhat unwieldy. Let $(G,u)$ be a dimension group with order unit; a normalized trace $\tau$ is {\it refinable\/} \wrt $u$ if for all $b, a_1', \dots, a_n'$ in $[0,u]$ \st $\tau(b) = \sum \tau(a_i')$, there exist $a_i \in G^+$ \st $\tau(a_i ) = \tau(a_i')$ and $b = \sum a_i$. If $\tau$ is refinable \wrt all order units (that is, the $b$ and $a_i'$ are allowed to be arbitrary elements of  $G^+$), then we simply say $\tau$ is refinable. It probably is the case that refinability \wrt a single order unit implies refinability. 

Analogous to goodness, we can define refinability for a set of traces. For a dimension group with order unit $(G,u)$, and  $U \subseteq S(G,u)$, we say that $U$ is {\it refinable\/} (\wrt $u$) if whenever $b, a_i' \in [0,u]$ and $\tau(b) = \sum \tau(a_i')$ for all $\tau$ in $U$, there exist $a_i \in G^+$ \st $b = \sum a_i$. If we let $L = \ker U = \cap_{\sigma \in U} \ker \sigma$, then the hypothesis is obviously $b - \sum a_i' \in L$, and this is equivalent to $\hat b|U = \sum \hat a_i'|U$, and this in turn is equivalent to $\hat b|Z = \sum \hat a_i'|Z$ where $Z = Z(\ker U) = \Set{\sigma \in S(G,u)}{\sigma(\ker U) = 0}$. Note that refinability of $\tau$ is equivalent to refinability of $Z(\ker \tau)$. The latter is of the form $S(G,u) \cap \Cal L$ where $\Cal L$ is a closed flat. We will see below that refinability of $\tau$ forces $Z(\ker \tau)$ to be a simplex, a fairly severe restriction.

If $\ker U = \brcs{0}$, then $Z(\ker U) = S(G,u)$, which of
course is good; therefore, $U$ is refinable. This generalizes the trivial
fact that a bad trace ($\ker \tau = \brcs{0}$) is refinable, Lemma \sevone\ below.

 An element $b$ of $G^+$ is called {\it refinable \wrt $\tau$\/} if whenever $\tau(b) = \sum \tau(a_i')$ with $a_i' \in G^+$, there exist $a_i$ in $G^+$ \st $b = \sum a_i$. Obviously $\tau$ is refinable if and only if every $b$ in $G^+$ is refinable \wrt $\tau$. We say $\tau$ is {\it weakly refinable\/} if there exists a generating set $H$ of the positive cone \st every member of $H$ is refinable \wrt $\tau$. Unlike the corresponding concepts for goodness, it is not clear that weak refinability implies refinability. First however, we have a skirmish of the definitions. 

\Lem Lemma \sevele. Suppose $(G,u)$ is a dimension group, $\tau$ is a trace that is group-like \wrt $u$,  and $b$ is a  element of $[0,u]$ that is refinable \wrt $\tau$. Then $b$ is $\tau$-good.

\Pf Recall that $\tau$ is group-like \wrt $u$ asserts that $\tau([0,u]) = \tau(G) \cap [0,\tau(u)]$. Suppose $a$ is an element of $G$ \st $0 <  \tau(a) , \tau(b-a)< \tau(b) \leq \tau(u)$. Since $u$ is $\tau$-good, there exist $c,c' \in [0,u]$ \st $\tau(c) = \tau(a)$ and $\tau(c') = \tau(b-a)$. Thus $\tau(b) =\tau(c) + \tau(c')$, and now we can apply refinability of $b$ to obtain $d$ and $d'$ in $G^+$ \st $b= d + d'$, and $d-c, d'-c' \in \ker \tau$. In particular, $d$ belongs to $[0,b]$ and $\tau(d) = \tau(a)$, so $b$ is $\tau$-good.\qed

Not only are good traces refinable, so are the bad ones (by a completely trivial argument)!

\Lem Lemma \sevone. Let $(G,u)$ be a dimension group with order unit, and let $\tau$ be a  trace. Any of the following conditions is sufficient for $\tau$ to be refinable.
\item{(i)} $\tau$ is good.
\item{(ii)} $\ker \tau = \Inf(G)$ and $G$ is simple.
\item{(iii)} $\ker \tau = \brcs{0}$.

\Pf (i) Given $b, a_i' \in [0,u]$, we may assume that $\tau(b) > 0$ (or else just replace $a_1'$ by $b$ and the rest by $0$) and  relabel so that $\tau(b) > \tau(a_1')$ (if $\tau(b) = \tau(a_1')$, set $a_1 = b$ and the rest zero); by goodness, there exists $a_1 \in [0,b]$ \st $\tau(a_1) = \tau(a_1')$, and apply the same procedure to $\tau(b-a_1) = \sum_{i\geq 2} \tau(a_i')$. The process stops when (after relabelling), $\tau (\sum_{i\geq l} a_i') = 0$, in which case we set $a_l = b-\sum_{i \leq l-1}a_i'$ and $a_{l+j} = 0$ for $j > 0$.

\noindent {(ii)} Since $\tau(b) = \sum \tau(a_i')$, we have $ b-\sum a_i' := k$ is an infinitesimal; replace one of the nonzero $a_j'$ by $a_j = a_j + k$, so that $a_j$ is still an order unit, and thus positive, and set all the other $a_i' = a_i$.

\noindent{(iii)} Obviously $\tau(b) = \sum \tau(a_i')$ implies $b = \sum a_i'$. \qed

As an elementary exercise, if $G$ is a simple dimension group with exactly two pure traces and $\Inf(G) = 0$, then a refinable trace is either good or bad. This fails with three pure traces  (Example \exten, $G'$).  

\Lem Corollary \sevtwe. Let $(G,u)$ be a dimension group with order unit, and let $\tau$ be a trace of $G$. The following are equivalent.
\item{(i)} $\tau$ is good
\item{(ii)} $\tau$ is both group-like and refinable \wrt $u$
\item{(iii)} $\tau$ is both group-like and weakly refinable \wrt $u$.
 
\Pf (i) $\implies$ (ii) is an immediate consequence of \appbfou\ and the definitions, and (ii) $\implies$ (iii) is trivial. 

\noindent (iii) $\implies$ (i). Let $H$ be a generating set for $G^+$ consisting of refinable elements. We observe that $H(u):= H \cap [0,u]$ is also a generating set: if $g \in G^+$ and $g \leq nu$, there exist $u_i \in [0,u]$ \st $g = \sum u_i$. Since each $u_i$ is represented as a sum of elements of $H$, each of which must lie in $[0,u]$, $H(u)$ is a generating set. By Lemma \appbfou, every element of $H(u)$ is $\tau$-good, and so by Lemma \onethr, all elements of $G^+$ are $\tau$-good, and thus $\tau$ is good.
\qed

The following is immediate from the definitions.  

\Lem Lemma \appbfiv. Let $(G,u)$ be a dimension group, and suppose that $L$ is a subgroup satisfying (i) and (ii), so that $G/L$ with the quotient ordering is a torsion-free partially ordered abelian group; let $\Arrow \pi;G.G/L$ denote the quotient map. Set $Z = \Set{\phi \in S(G,u)}{\phi(L) =0}$. The map $S(G/L,u+L) \to Z$ given by $\psi \mapsto \psi \circ \pi$ is an affine homeomorphism.

The following does not yield the expected result, that $G/\ker U$ is a dimension group if $U$ is refinable; this would require $G/\ker U$ to be unperforated, which is not clear at the moment. It is however, simple and torsion-free. Simple torsion-free interpolation groups which are not unperforated exist; the first example was constructed by John  Lawrence [L]. 

However, if $U$ is good (\wrt $G$) and $G$ is simple, then we will show that $G/\ker U$ is  unperforated, en route to showing it is a simple dimension group.  Let $J$ be a torsion-free partially ordered abelian group (which will be set to $G/\ker U$), {\it not\/} assumed to be unperforated, and let $a$ be an element of $J$. Define a  (possibly empty) subset of $\N$ (the set of strictly positive integers), 
$P(a) = \Set{n \in \N}{na \in J^+}$. We note that $P(a)$ is closed under addition, so if $P(a)$ contains a  subset whose greatest common divisor is $1$, then $P(a)$ is cofinite. In general, $P(a)$ is either empty or  of the form $k\N \setminus T$ for some positive integer $k$ and a finite set $T$.

If $J$ is simple, then $P(a)$ is either empty or cofinite.  For nonzero $a$, if $P(a)$ is not empty, choose $n \in \N$ \st $na \in J^+$; since $J$ is torsion-free, $na$ is nonzero, hence (being in $J^+$) is an order unit. Thus there exists a positive integer $N$ \st $a \leq Nna$. This implies $Nn-1 \in P(a)$, and since $\gcd(n,Nn-1) = 1$, $P(a)$ is cofinite.
 
\Lem Proposition \appbsix. Suppose $(G,u)$ is a simple dimension group. Let $U$ be a collection of traces, and form $Z(\ker U) = \Set{\sigma \in S(G,u)}{\sigma(\ker U) = 0}$. 
\item{(a)}If $U$ or $Z(\ker U)$ is refinable, then the quotient group $G/\ker U$ satisfies Riesz interpolation, and thus $Z(\ker U)$ is a simplex.
\item{(b)} If $U$ is good, then $G/\ker U$ is unperforated.  
\item{(c)} If $U$ is good, then $U$ is refinable.
\item{(d)} If $U$ is good, then $G/\ker U$ is a simple dimension group.
\item{(e)} If $Z(\ker U)$ is 
refinable and $G/\ker U$ is unperforated, then  $\ker U$ has dense image in $Z(\ker U)^{\vdash}$.  
\item{(f)} If $Z(\ker U)$ is refinable  and $G/\ker U$ is unperforated, then $Z(\ker U)$ is good.

\Pf (a) It suffices to show $G/\ker U$ satisfies interpolation. By [G; 2.1(c)], sufficient is that if $a_i$, $b_i$ are positive elements of $G/\ker U$ \st $b_1  + b_2 = a_1 + a_2$, then we may find $a_{ij} \in (G/\ker U)^+$ \st $b_i = \sum_j a_{ij}$ and $a_j = \sum_i a_{ij}$. If any of $\brcs{a_i,b_j}$ is zero, there is nothing to  do. Since $a_i, b_j$ are in $(G/\ker U)^+$, by definition of the quotient ordering, there exist $A_i, B_j \in G^{++}$ \st $a_i = A_i + \ker U$, $b_j = B_j + \ker U$. Then $B_1 + B_2 = A_1 + A_2 - k$ for some $k \in \ker U$. By refinability applied to the $A_i$, there exist positive elements $A_i'$ \st $A_i - A_i' \in \ker U$ and $B_1 + B_2 = A_1' + A_2'$. Now interpolation applies in $G$, yielding $A_{ij}$ in $G^+$ satisfying the usual equations, and their images, $a_{ij} = A_{ij} + \ker U $ satisfy the corresponding equations,  yielding the Riesz interpolation property. By \appbfiv, the trace space of $G/\ker U$ is affinely homemorphic to $Z$ (via $\psi \mapsto \psi \circ \pi$), and it follows that $Z$ is a simplex.

\noindent (b) Let $a$ be a nonzero element of $J = G/\ker U$ with nonempty $P(a)$. Set $L = \ker U$. Suppose $m,n \in P(a)$ with $ m> n$. We show that $m-n \in P(a)$. Since $ma$ and $na$ are in the positive cone of $G/L$, from the definition of the quotient ordering, there exist $A$, $B$ in $G^+$ \st $B +L = ma$ and $A+L = na$. Since $G$ is simple, $\tau(B) $ is bounded below away from zero as $\tau$ varies over $S(G,u)$, and thus if ${\tilde \tau} \in S(G/L,u+L) $ corresponds to $\tau \in Z$, $\tilde \tau(ma) $ is bounded below away from zero. Hence $\tau(B) = \tilde\tau(ma) = m\tilde\tau(a) > n\tilde\tau(a) = \tau(A)$ and the difference is bounded below away from zero. Thus
in $G$, we have $\hat B|Z -\hat A|Z$ is bounded below away from zero, that is, $\tau(B - A)$ is bounded below as $\tau$ varies over $U$ or $Z$. By goodness of $U$, there exists $A' \in [0,B]$ \st $A - A' \in \ker Z$. Thus $B-A' \in G^+$, so its image in $G/L$, $(m-n)a$, belongs to $(G/L)^+$, so $m-n \in P(a)$. This forces $P(a) = k\N$ for some nonnegative integer $k$ (if $P(a)$ is nonempty, let $k$ be its smallest element, and use the division algorithm). Since $J$ is simple, $P(a)$ is either cofinite or empty, hence $k=1$, and thus $J$ is unperforated. 

\noindent (c) Suppose  $a_i',b \in G^+$ and $b - \sum_{i=1}^n a_i' = k \in \ker U$. Assuming $a_1 \neq 0$, we have $a_1'$ is thus an order unit; hence there exists $\delta > 0$ \st $\hat a_1 \gg \delta \pmb 1$. Hence $\hat b|U - \hat a_n|U \gg \delta \pmb 1$, whence $b - a_n' + \ker U$ is an order unit of $G/\ker U$. We may thus lift $b-a_n' + \ker U$ to an order unit $c'$ of $G^+$; since $b - c' + \ker U = a_n' + \ker U$, goodness of $U$ applies. There thus exists $c \in [0,b]$ \st $c' - c \in \ker U$. Now $b-c = \sum_1^{n-1} a_i + (k +c-a_n)$. The parenthesized term is in $\ker U$, and induction can now be applied. 

\noindent (d) By (c) and (a), $G/\ker U$ satisfies interpolation; by (b), it is unperforated. Simplicity is automatic.

\noindent (e) Select $h$ in $Z^{\vdash} = \Set {f \in \Aff K}{f| Z \equiv 0}$.  By approximate divisibility and consequent density of $G$ in $\Aff K$, given $\epsilon > 0$, there exists $g \in G$ \st $\| \hat g - h\| < \epsilon$. In particular, for all $\tau \in Z$, $|\tau(g)| < \epsilon$. By approximate divisibility again, there exists a (small) order unit, $e$, \st $\epsilon \pmb 1 < \hat e < 2 \epsilon \pmb 1$. Hence $\tau(g+e) > 0$ for all $\tau$ in $Z$, and $\| \hat g+\hat e -h\|< 3\epsilon$. Set $g_1 = g +e$.

Since $\hat g_1|Z \gg 0$, the element $g_1 + \ker U$ is an order unit of the unperforated group $G/\ker U$ ($\ker U$ automatically satisfies (i) and (ii), and by hypothesis (iii)), hence there exists an order unit $v$ of $G$ \st $\hat v|Z = \hat g_1|Z$. Let $f$ be a (small) order  unit of $G$ \st $3 \epsilon < \hat f < 4\epsilon$. Since $\hat g_1|Z < 3\epsilon$, it follows that $(\hat f - \hat v)|Z$ is an order unit of $G/\ker U$. There thus exists an order unit $w$ of $G$ \st $(\hat f - \hat v)|Z = \hat w|Z$. Thus $(\hat f - \hat w - \hat v)|Z \equiv 0$; hence there exists $k$ in $\ker U$ \st $f = w  + v +k$. Since $f$, $w$, and $v$ are all order units, hence are positive, we can apply refinability. There thus exist positive $w'$ and $v'$ \st $w' - w\in \ker U$, $v' - v \in \ker U$ and $f = w' + v'$. Thus $v' \in [0,f]$, so $\| \hat v'\| < 4\epsilon$.

Now consider $g':= g_1 - v'$. Restricted to $Z$, its value is zero (since $v'|Z = v|Z = g_1|Z$), so $g' \in \ker U$. On the other hand $\| \hat g' - h\| \leq \| \hat g_1 -h\| + \| v'\| < 3\epsilon + 4\epsilon = 7\epsilon$. So $\ker U$ has dense image in $Z(\ker U)^{\vdash}$. 

\noindent (f) The quotient being unperforated implies that order units lift (Lemma \appbthr). Suppose $a$ and $b$ are nonzero elements of $G^+$ (therefore are order units) \st $\tau(b) - \tau(a) > 0$ for all $\tau$ in $Z$. Since the set of traces of the quotient induce the elements of $Z$, $b - a + \ker U$ is an order unit of $G/\ker U$. Since order units lift, there exists an order unit $c$ in $G$ \st $c - b+a \in \ker U$. Hence $\tau(b) = \tau(a) + \tau(c)$ for all $\tau$ in $Z$. By refinability of $Z$, there exist $a', c' \geq 0$ \st $b = a' + c'$; thus $a' \in [0,b]$ and $a' - a \in \ker U$. \qed

A particular consequence of (b) is that if $Z(\ker U)$ is good, then $U$ is refinable (since $U$ is refinable if and only if $Z(\ker U)$ is refinable). The converse is almost true; it may be that refinability implies the quotient is unperforated. If $G$ is divisible, then $G/\ker U$ is automatically unperforated, and in this case, refinability of $U$ is equivalent to goodness of $Z(\ker U)$.

The difference between $U$ being good and $U$ being refinable should be clarified. We have seen that $U$ is refinable if and only if $Z(\ker U)$ is refinable, and the latter is of the form $K \cap {\Cal L}$ where $K$ is the trace space, and $\Cal L$ is a closed flat obtained as an intersection of zero sets of $\hat g$ where $g \in U$, meaning $\Cal L = \cap_{{g \in \ker U}} \hat g^{-1}(0)$. There is a similar characterization for  good $U$, involving a (generally) smaller  flat. 

For a set $U \subset K$, define $L(U)$ to be the smallest closed flat containing $U$; that is, $L(U)$ is the closure of the set $\Set{\sum \alpha_i u_i}{\alpha_i \in \R;\ \sum \alpha_i = 1; \ u_i \in U}$. Define ${\widetilde U} = K \cap L(U)$. If $U$ is a singleton, then $L(U)= {\widetilde U} = U$, and in general, ${\widetilde U} $ is a compact convex set contained in $Z(\ker U)$. Alternatively, we could define $\widetilde U$ as $K \cap \(\cap_{\Set{h \in \Aff K}{h|U \equiv 0}}  h^{-1}(0)\)$.
\comment
, so the difference between ${\widetilde U}$ and is that in the latter, the set over which we take the intersection is $U^{\vdash}$, while in the latter, it is just $ G \cap \ker U$, whose image in $\Aff K$ is typically much smaller than $U^{\vdash}$. 
\endcomment
Both are of the form $K \cap {\Cal L}$ where $\Cal L$ is a closed flat.

There is a property related to  goodness, which however, does not seem to lend itself to any particular use. We say that a subset $U$ of a  dimension group $(G,u)$ with order unit is {\it poob\/} (this is what handwritten {\it good\/} approximately looks like when viewed upside down) if for all $a$ and $b$ in $G^+$ \st $\tau(b) - \tau(a)\geq \delta >0$ for some $\delta$ and all $\tau$ in $U$, there exists $b'$ in $G^+$ \st $\tau(b) = \tau(b')$ and $a \leq b'$. 

\Lem Proposition \appbtwo. Suppose that $(G,u)$ is simple and not cyclic, and $U$ is a subset of $S(G,u)$.  If  $\ker U$ has dense image in $Z(\ker U)^{\vdash}$, then
$U$ and $Z(\ker U)$ are  poob. In particular, if $U$ is good, then it is poob. 

\Remark But not conversely! By Lemma \appbsix(a), if $U$ is good, then $Z(\ker U)$ is a simplex, but we have already seen a finite dimensional example for which $\ker U$ has dense image, but $Z(\ker U)$ is not a simplex.

\Pf Suppose $a$ and $b$ are in the positive cone of $G$ and $\tau(b) - \tau(a)> \delta> 0$ for all $\tau \in U$ and some $\delta >0$. By definition, the same inequality holds for all $\tau$ in $Z(\ker U)$. Let $H = \Aff S(G,u)$ and $L= Z(\ker U)^{\perp}$. Applying Lemma \appbthr\ to $\hat b - \hat a$ (restricted to $Z(\ker U)$, this element is strictly positive), there exists $f \in \Aff S(G,u)$ \st $f\gg 0$ and $f|Z(U) = (\hat b - \hat a)|Z(U)$, that is $e:= f - (\hat b - \hat a)$ is in $Z(U)^{\perp}$. 

Given $\epsilon > 0$, there exists (by the density hypothesis) $j \in \ker U$ \st $\|\hat j - e\| < \epsilon$, so that $\| \hat j - \hat a + \hat b - \hat f\| < \epsilon$. Hence if $\epsilon $ is chosen less than $\inf_{\sigma \in S(G,u)} f(\sigma)$ (this is of course strictly positive, since $f$ is a continuous strictly positive function on a compact space), then $j - a + b$ is an order unit of $G$. In particular, $b+j \gg a$ so the former is positive, and we can set $b' = b+j$. 

If $U$ is good, then \appbsev(a) and \appbsix(a,c) apply  to yield density.
\qed

The following is almost immediate from the definitions.

\Lem Lemma \appbsev. Let $(G,u)$ be an approximately divisible dimension group and $U$ be a subset of $S(G,u)$. {\par}
\item{(a)} If $U$ is a good subset of $S(G,u)$ then ${\widetilde U}$ is a good subset of $S(G,u)$, and in that case, $Z(\ker U) = {\widetilde U}$ and is a simplex.
\item{(b)}  $U$ is a refinable subset of $S(G,u)$ if and only if $Z(\ker U)$ is a refinable subset of $S(G,u)$.

\Pf The only assertion that is not tautological is that if $U$ is good, then $Z(\ker U) = {\widetilde U}$. Suppose $h \in \Aff K$ vanishes identically on $U$; we show $h$ is approximable by the images of elements of $\ker U$. By approximate divisibility, given $\epsilon > 0$, there exists $g_1 \in G$ \st $\| h - \hat g_1\| < \epsilon$. Then $|\tau(g_1)| < \epsilon$ for all $\tau \in U$. There exist small order units $v$ and $b$ of $G$ \st $\epsilon \pmb 1 <\hat v < 2\epsilon \pmb 1$  and $4\epsilon \pmb 1 <\hat b < 5\epsilon \pmb 1$ (again, by density of $\hat G$ in $\Aff K$, a consequence of approximate divisibility). Then $\hat v+ \hat g_1$ is strictly positive on $U$. By goodness applied to $b$ and $v+g_1$, there exists $c \in [0,b]$ \st $\hat c|U = \hat v|U + \hat g_1|U$, that is, $\tau(c) = \tau(v+g_1)$ for all $\tau \in U$. Since $0 \leq c \leq b$, it follows that $\| c \| \leq \| b\| < 5 \epsilon$. Set $g = v + g_1 -c$; this belongs to $\ker U$ and $\| \hat c - h \| \leq  \| \hat g_1 - h\| + \|\hat v \| + \| \hat c\| < 8 \epsilon$. 

Thus the image of $\ker U$ is dense in $U^{\vdash}$, and the equality $Z(\ker U) = {\widetilde U}$ follows. \qed

While Lawrence's example of a simple torsion-free interpolation group is complicated, we can give an easy example of a simple torsion-free partially ordered group with exactly two pure traces (thus the trace space is a simplex) for which unperforation fails (unfortunately, it is not clear whether interpolation holds). Here $G$ is a simple dimension group with three traces, $U$ is a subset of $S(G,u)$, $Z(U)$ is a line segment and  $G/\ker U$ is not unperforated. 

 Consider the rank 4 critical group of the form $G = \langle e_1, e_2, e_3; \sum \alpha_i e_i\rangle$ in $\R^3$ with strict ordering inherited from $\R^3$; the  pure traces are the evaluations at the coordinates, $\tau_i$. Let $U = \brcs{(2\tau_1 + 3\tau_2)/5, (3\tau_1 + 5\tau_3)/8}$ (the denominators are normalizations at the order unit $(1,1,1)$; they can be dispensed with). Then $Z(\ker U)$ is a line segment running from a point on the edge joining $\tau_1$ to $\tau_2$ to a point on the edge joining $\tau_1$ to $\tau_3$, therefore a simplex. It is easy to calculate $\ker U = k\Z$ where $k = (15,-10,-9)$.

Now consider the element $g = (2,-1,-1)$. We note that $8g -k = (1,2,1) \gg (0,0,0)$, so that $8(g + \ker U) \in (G/\ker U)^+$. However, there is no integer $n$ \st all three coordinates of $(2,-1,-1) + n(15,-10,-9)$ are positive: the first coordinate would require $n \geq 0$, but the second requires $n < 0$. So $G/\ker U$ is not $8$-unperforated, even though $Z$ is a simplex. In addition, $17g -2k$ belongs to $G^+$, and this  means that there exists $m$ \st for all $t \geq m$, $t(g +\ker U)\in (G/L)^+$ (we can take $m = 112 = (8-1) \times (17 - 1)$, probably much less). A similar but easier example is discussed in Appendix B.\qed

What is missing from the litany of implications is that if $Z = Z(\ker U)$ has dense image in $Z(\ker U)^{\vdash}$ and $Z(\ker U)$ is a simplex, then $Z$ is good for $G$. This is false as we will see in Lemma \appbele, but it can be reduced to a problem concerning $\Aff K$ ($K$ a Choquet simplex). If $Z$ happens to be a face, then it is good, as will be a consequence of an easy computation, and as we saw is a consequence of  a liftability statement about faces in Choquet simplices. There are however lots of choices of $Z$ which are not faces, which we will show to be good. 

Another question arising from this result is whether refinability of $U$ implies $G/\ker U$ is unperforated---which would imply that $U$ is refinable if and only if $Z(\ker U)$ is good.

The  sets $Z(\ker U)$ create a Zariski-like topology on $S(G,u)$; this topology depends on $G$, not just the simplex $S(G,u)$---the closed sets come from intersections of kernels of traces; that is, $Z(\ker U) $ is of the form $S(G,u) \cap {\Cal L}$ where ${\Cal L}$ is the closed flat obtained as the intersection of zero sets of elements of $\ker U$.  Good closed sets are of the form $Z(\ker U)$ where the latter is a simplex and $\ker U$ is dense in ${\tilde U}^{\vdash}$ (not $Z^{\vdash}$), but these conditions are not sufficient. In case $G = \Aff K$ itself, this topology extends the   facial topology (defined on $\partial_e K$).

Fix  $Z \subseteq K$ (both simplices); $g \in \ker Z$ is equivalent to $\hat g|Z \equiv 0$. Define
 $$
R(G,Z)= \Set{(a,b) \in G^{++}\times G^{++}}{\text{there exists $j \in \ker Z$ \st $0 \ll \hat a - \hat j \ll \hat b$}}.
$$
\Lem Lemma \appbeig. Let $(G,u)$ be a simple dimension group. Then $R(G,Z)$ is an open subset of $G^{++} \times G^{++}$ \wrt the supremum pseudo-norm; explicitly,  for all $(a,b) \in R(G,Z) $, there exists $\epsilon > 0$ \st if $\| \hat a - h \|, \|\hat b - k \| < \epsilon$, then $(h,k) \in R(\Aff S(G,u),Z)$ (as usual, $K = S(G,u)$; in addition, in this case, $\ker Z$ is replaced by $Z^{\vdash}$).

\Pf The condition for membership in $R$ is equivalent to $\hat a \gg \hat  j \gg \hat a- \hat b$ (for some $j$ in $\ker Z$). Hence if we pick $\epsilon < \min\brcs{\inf_{\sigma \in K} \sigma(a-j) ,  \inf_{\sigma \in K} \sigma(j -a +b) }/2$, we can use the image of the same $j$, $\hat j$, for $(h,k)$. A particular consequence is that if $Z(\ker U)$ is good \wrt $G$, then it is also good \wrt $\Aff S(G,u)$. The same argument also shows that if $(a',b') \in G^{++} \times G^{++}$, and $\| \hat a - \hat a'\|$ and $\| \hat b -\hat b'\|$ are sufficiently small, then $(a',b') \in R(G,Z)$.\qed

\Lem Lemma \appbnin. Let $(G,u)$ be a simple dimension group, let $U \subseteq S(G,u)$; let $Z \equiv Z(\ker U) $ denote $\Set{\sigma \in S(G,u)}{\sigma (\ker U) = 0}$. The following are equivalent:
\item{(i)} $Z$ is good as a family of traces on $\Aff S(G,u)$ (\wrt the strict ordering) and $\ker U$
 has dense range in $Z^{\vdash} = \Set{f \in \Aff S(G,u)}{f|Z \equiv 0}$
\item{(ii)} $Z$ is good as a family of traces on $G$.

\Pf (i) $\implies$ (ii) We may replace $U$ by $Z$, since $\hat c|Z$ is uniquely determined by $\hat c|U$ for all $c $ in $G$.If $a,b \in G^+$, $0 \ll \hat a|Z \ll \hat b|Z$, then from $Z$ being good as a set of traces on $\Aff S(G,u)$, there exists $h$ \st $0 \ll h \ll \hat b$. Hence $R(\Aff S(G,u),Z)$ contains $(\hat a, \hat b)$. There thus exists $k \in Z^{\vdash}$ \st $\hat a \gg k \gg \hat a - \hat b$. For $\epsilon <\min\brcs{ \min_{\sigma \in S(G,u)} (\hat a- k)(\sigma),\min_{\sigma \in S(G,u)} (k - \hat a +\hat b)(\sigma)}$, there exists $j$ in $\ker Z = \ker U \subset G$ \st $\| \hat j - k \| < \epsilon$. Then $\hat a \gg \hat j \gg \hat a - \hat b$, and since $G$ is unperforated, $j \in R(G,U) $, hence the latter is nonempty. Thus $U$ is good.

\noindent (ii) $\implies $ (i) Density follows from Proposition \appbsix(a,c).  To show $Z$ is good relative to $\Aff S(G,u)$, suppose $f,h \gg 0$ are elements of $\Aff S(G,u)$ \st $h|Z \ll f|Z$. By density of the image of $G$, given $\epsilon$, there exists $a \in G$ \st $h \gg \hat a$ and $h - \hat a \ll \epsilon \pmb 1$, and by making $\epsilon$ sufficiently small, we can ensure that $\hat a \gg 0$, so that $a$ is an order unit. Next, we may approximate $f$ by $\hat b$ with $b \in G$ \st $3 \epsilon \pmb 1 \gg f - \hat b \gg h - \hat a$. In particular, if $\epsilon $ is sufficiently small, it will follow that $\hat b \gg \hat a$. By goodness of $U$ (or of $Z$ relative to $G$), there exists $j$ in $\ker U$ \st $\hat a \gg \hat j \gg \hat a - \hat b$. However, $h \gg \hat a$ and $\hat a - \hat b \gg h - f$, yielding $h \gg \hat j \gg h-f$, and thus $f \gg h - j \gg 0$. \qed

By \appbnin\ and \appbsix, deciding whether $Z \equiv Z(\ker U)$ is good  \wrt  $G$  boils down (assuming $\ker U$ has dense image in $Z^{\vdash}$) to whether $Z$ is good for the huge group $\Aff S(G,u)$. 
For example, this is the case if $Z = \brcs{\tau}$ (a singleton) for any $\tau \in K$ (that is, any singleton is a good set of traces for $\Aff K$). If $Z$ happens to be a face, then $Z$ being good is simply a consequence of the lifting property discussed near the beginning of this section; all that is really necessary is the interpolation property of the usual ordering.

 An example wherein a face arises out of a simple dimension group occurs in our standard critical example, $G = \brcs{e_i, \sum \alpha_j e_j}$; then $\ker \tau_1$ is spanned by  elements whose images in $\Aff S(G,u) = \R^n$ are nonnegative, $\brcs{e_i}_{i \neq 1}$, and $Z(\ker \tau_1)$ is the face with vertices  $\brcs{\tau_2, \dots, \tau_n}$. Another, rather odd example is  line segment from a vertex to a point in the relative interior of a facet. This is neither a face nor a singleton; however, not many other types of $Z$ can be good, in view of the next result. 

In the following, the point $v$ need not be an extreme point; it is however, in a proper face of $K$.

\Lem Lemma \appbele. Suppose $K$ is a simplex and set $Z = K \cap {\Cal L}$ where ${\Cal L}$ is a closed flat. Suppose in addition that  there exists a closed face,  $F$, having more than one extreme point, \st $F \cap Z = \emptyset$, and  a point $v$ in $K \setminus Z$ from which there exist two distinct line segments  to $F$, both of which pass through $Z$. Then $Z$ is not good. 
{\par}
In particular, this applies if there exists a face $F$ of nonzero dimension \st $Z \cap F = \emptyset$ and $Z$ contains an interior point of $K$.

\Pf Since $F \cap Z $ is empty, $F$ is a closed face, and $Z$ is compact, there exists $b_1 \in \Aff K$ \st $ b_1 \geq 0$, $b_1|F \equiv 0$, and $b_1|Z \gg 1$. Let $x_1$ and $x_2$ be the endpoints of the line segments in $F$ (joining $F$ to $v$), and let $z_1$, $z_2$ be the intersection points of the two lines with $Z$. We can write $z_i = \lambda_i x_i + (1-\lambda_i) v$ for some $\lambda_i \in (0,1)$. By interchanging $x_1$ with $x_2$ if necessary, we may assume that 
$$
\frac{\lambda_2(1-\lambda_1)}{\lambda_1(1-\lambda_2)} \leq 1. 
$$
We may also find $a_1 \in \Aff K$ \st $a_1 \geq 0$, $a_1 (x_1) \geq \slfrac34$,  $0 < a_1(x_2) < \slfrac34 - 3\epsilon$ for some $0< \epsilon <\min\Set{b_1(k) -1}{k\in Z}$ and $a_1|Z \leq 1$. Hence $b:= b_1 + \epsilon \pmb 1$ and $a:= a_1 + (\epsilon/2)\pmb 1$ satisfy $a,b \gg 0$ and $b|Z \gg a|Z$. 

Suppose that $Z$ were good (for $\Aff K$). Then there would exist $c \in \Aff K$ \st $0 \leq c \leq b$ and $c|Z = a|Z$. Let $K'$ denote the triangle with vertices $\brcs{x_1, x_2, v}$ and let $Z'$ be $K'\cap Z$, which it is easy to see is just the line segment joining $z_1$ to $z_2$. Then $c|Z' = a|Z'$ and $0 \ll c|K' \ll b|K'$. We show this is impossible. Denote $v$ by $x_0$.

Since $z_i = \lambda_i x_i + (1-\lambda_i) x_0$ ($i=1,2$), on writing $c(x_i) = a(x_i) - \delta_i$ ($i = 0,1,2)$, we find $0 = \lambda_i \delta_i + (1-\lambda_i) \delta_0$ ($i=1,2$), whence $\delta_i = - (1-\lambda_i) \delta_0/\lambda_i$ ($i = 1,2$). The condition that $c(x_i) > 0$
 yields three inequalities for $\delta_0$, computed from $\delta_i < a(x_i)$.
$$
\max\brcs{-\frac{\lambda_1 a(x_1)}{(1-\lambda_1)}, -\frac{\lambda_2 a(x_2)}{(1-\lambda_2}}< \delta_0 < a(x_0).
$$
The conditions $c(x_i) < b(x_i)$ yield three more inequalities, 
$$
a(x_0) - b(x_0) < \delta_0 < \min\brcs{\frac{(b(x_1) -a(x_1))\lambda_1}{1- \lambda_1},\frac{(b(x_2) -a(x_2))\lambda_2}{1- \lambda_2} }.
$$
Combining all of them, we see that 
$$
\max\brcs{a(x_0) - b(x_0),-\frac{\lambda_1 a(x_1)}{1-\lambda_1}, -\frac{\lambda_2 a(x_2)}{1-\lambda_2}} < \min\brcs{a(x_0),\frac{(b(x_1) -a(x_1))\lambda_1}{1- \lambda_1},\frac{(b(x_2) -a(x_2))\lambda_2}{1- \lambda_2} }.
$$
In particular,
$$\eqalign{
b(x_1) &> a(x_1) - \frac{\lambda_2 (1-\lambda_1) }{(1-\lambda_2)\lambda_1}a(x_2) \cr
& \geq a(x_1) - a(x_2) > \frac 34 - \(\frac 34 - 3 \epsilon\) = 3 \epsilon.
}$$
This contradicts $b(x_1) = \epsilon$. 
\qed

Lemma \appbele\ applies if $K$ is finite dimensional and $Z = K \cap {\Cal L}$ is a codimension one simplex that misses at least three vertices and is not a face (for then the hyperflat must separate the one of the vertices from the other two). In particular, if the hyperflat misses all the vertices, $Z$ is not good. It also applies for most (but not all) $Z$ that are  line segments from one facet to another that missing all the vertices of  $K$ (exceptions occur when the endpoints are in faces complementary to each other).

Here is an easy construction of a good  $Z = K \cap {\Cal L}$ that is neither a face nor a singleton, corresponding to the coproduct of two good subsets (for a discussion of coproducts of simplices in the category of simplices, see [G, p\,164]).  Let $K $ be a simplex, and let $F$ and $F'$ be disjoint closed faces. Select $W$, a  subset of $F$ that is good \wrt $F$, and $W'$ a subset of $F'$ that is good \wrt $F'$. Let  $Z \equiv Z(W,F; W',F') \subset K$  be the convex hull of $W \cup W'$ inside $K$. The convex hull of $F \cup F'$ will be denoted $F \vee F'$, and is itself a  face of $K$. 

It is completely routine to verify that this construction yields good $Z$. If $a,b$ are strictly positive elements of $\Aff K$ \st $a|Z \ll b|Z$, then of course $a|W \ll b|W$ and $a|W' \ll b|W'$, so by goodness of $W,W'$ relative to $F,F'$ respectively, there exist strictly positive $d \in \Aff F$ and $d' \in \Aff F'$ \st $d|W = a|W$,  $d'|W' = a|W'$, and  $d \ll b|F$, $d' \ll b|F'$. Since $F$ and $F'$ are disjoint and closed, there exists $c \in \Aff (F \vee F')$ \st $c|F = d$ and $c|F' = d'$, and it is obvious that $0 \ll c|(F \vee F') \ll b|(F \vee F')$ and $c|Z = a|Z$. Since $F \vee F'$ is itself a closed face of $K$, we may apply the extension result for closed faces to find $c' \ll b$ \st $c'|(F \vee F') = c|(F \vee F')$. 

This generalizes the two-dimensional situation, wherein $K$ is a triangle and $Z$ is the line segment from  a vertex to a point on the opposite edge. This allows us obtain building blocks for good subsets of $K$. The building blocks are singleton subsets in the boundary and closed faces (disjoint from the face containing the singleton). The process can be iterated finitely many times (or simply deal with finitely many  disjoint faces), but it is not clear that  we can take the closure of the increasing union if infinitely many disjoint faces or singleton sets are involved. 

Thus a line segment in the triangle from a vertex is the coproduct of a zero-dimensional face (the vertex) and a point in the complementary edge. The smallest example for which we can construct a coproduct of two singletons that hits no vertices occurs in a tetrahedron. Let $K$ be such, with vertices $\brcs{x_0, x_1, x_2, x_3}$. Let $z_1$ be a point in the relative interior of the edge $F$ joining $x_0$ to   $x_1$, and choose $z_2$  in the relative interior of the edge $F'$ joinining $x_2$ to $x_3$. The construction yields a good line segment joining $z_1$ to $z_2$; it  contains no vertices of  $K$. 

On the other hand, if the two building blocks are both faces, then the construction yields only a  face. In particular, this means in any sequence of building blocks (closed faces and singleton sets), we can combine all the closed faces in one (the closure of the convex hull of a union of faces is a face), and then we are simply adding singleton sets in disjoint faces (which have to be complementary to the previous ones).

It is also routine to verify that if $W$ and $W'$ are intersections of $F$ and $F'$ respectively with flats, then the resulting $Z(W,F;W',F')$ is also of this form, that is, the closure of the affine span of $Z$ intersected with $K$ is just $Z$.

\comment
Select a point $z$ in $F'$, and set $Z $ to be the convex hull of $F \cup \brcs{z}$. The affine span consists of finite sums of the form $v= \alpha_{0,v} z + \sum \alpha_{i,v} f_i$ with $\sum \alpha_{i,v} = 1$ (and of course, some of the $\alpha$s can be negative). If $w_{j}$ is a sequence of elements of the affine span that converges  to $w \in K$, on applying the linear functional that exposes $F$, we see that $\alpha_{0,w_j}$ converges, say to $\alpha_0$. Thus $w_j - \alpha_0 z \to w - \alpha_0 z$ in the affine span of $F$. On writing $w = \lambda e + (1-\lambda )e'$ with $e \in F$ and $e' \in F'$, we see that $(\alpha_0 z - (1-\lambda)e' )/\lambda $ belongs to the affine span of $F$, which forces $e' = z$, and thus $w \in Z$. In particular, $Z$ is the intersection of a closed flat (the closure of the affine span of $Z$ itself) with $K$. 
\endcomment

\Lem  Conjecture. If  $K$ is a simplex and $Z = K \cap {\Cal L}$ where ${\Cal L}$ is a closed flat, then $Z$ is good as a set of traces on $\Aff K$ (with the strict ordering) if and only if $Z $ is obtained as the closure of an increasing  union of $Z_i = K \cap {\Cal L}_i$, where each $Z_i$ is obtained as a coproduct of $Z_{i-1}$ with a building block.

 If $K$ is finite dimensional, the conjecture is almost certainly true, and we can just forget the limit process. 
\comment 
\Lem Lemma \appbten. Suppose that $K$ is a finite dimensional simplex, and $\Cal L$ is a hyperplane  \st $Z = K \cap \Cal L$ is a simplex of codimension one in $K$. Then
\item{(a)}  $\min\brcs{|V_+|,|V_{-}|} \leq 1$.
\item{(b)} $Z$ is good for $\Aff K$ \wrt the strict ordering.

\Pf Without loss of generality, we may assume $K$ is the standard solid simplex (with vertices \brcs{\pmb 0, e_i}$) in $\R^n$. (a) For any $x$ in $V_+$ and $y \in V_{-}$, there exists unique $1 > \lambda > 0$ (depending on $x$ and $y$) \st $v_{x,y} = \lambda x + (1-\lambda)y$ belongs to $Z$. The claim is that each $v_{x,y}$ (as $(x,y)$ vary over $V_{+} \times V_{-}$) is an extreme point of $Z$. In addition the points of $V_0$ (if any) are obviously vertices of $Z$.

Suppose $v_{x,y} = \alpha_{x,y} z + (1-\alpha_{x,y} )w$ for $w, z \in K$ and $\alpha_{x,y} \in (0,1)$. We first observe that $v_{x,y}$ belongs to the edge $E$ of $K$ with endpoints $\brcs{x,y}$; since this is a face of $K$, we must have that $w, z  \in E \cap Z$.  However, $E \cap Z$ is a singleton, so that $w = z = v_{x,y}$. Hence every $v_{x,y}$ is an extreme point of $Z$.

\noindent (a) Thus   $Z$ will have at least $e= |V_0| + |V_{+}|\cdot |V_{-}|$ extreme points. If both $|V_{\pm}|> 1$, then $e > n = |V_0| + |V_{+}| + |V_{-}| -1$, contradicting $Z$ being a simplex of dimension at most $n-1$.

\noindent (b) By replacing $f$ by $-f$ if necessary, we may assume $|V_+| \leq 1$. If $V_+ = \emptyset$, then $Z$ is the face spanned affinely by $V_0$, and we already know that a closed face is good (for $\Aff K$). Otherwise, $V_{+} = \brcs{x_0}$.
Since $Z$ is an $n-1$-dimensional simplex,  the vertices  $\brcs{v_{x_0,y}} \cup V_0$ exhaust the vertices of $Z$. Index the vertices of $K$ $\brcs{x_0, x_1, \dots, x_j; x_{j+1}, \dots , x_n} $ where $V_{-} = \brcs{x_1, \dots, x_j}$ and $V_0 = \brcs{x_{j+1}, \dots , x_n}$; either (but not both) $V_{-} $ or $V_0$ may be empty---but if $V_{-}$ is empty, the situation is decidedly uninteresting---it means $n = 1$, which is a  trivial case.

Suppose $h$ and $g$ are affine functions on $K$ \st $h \gg 0$ and $h|Z \gg g|Z \gg 0$.  Set $v_t = v_{x_0, x_t}$ (for $1 \leq t \leq j$) and $v_t = x_t$ for $t > j$; these are the vertices of $Z$. Set $h(x_i) = \beta_i > 0$ (for $i \geq 0$). Then for $1 \leq t \leq j$, $v_t = \alpha_{t} x_t + (1-\alpha_t) x_0$, so $h(v_i) = \alpha_t \beta_t + (1-\alpha_t)\beta_0$, and if $t > j$, $h(v_t) = h(x_t) = \beta_t$.

Since $h|Z \gg g|Z \gg 0$, we must have $0 < g(v_i) < h(v_i)$. Hence we may replace $h$ by an affine function $h'$ \st $ h \geq h' + \epsilon\pmb 1 \geq h - \epsilon \pmb 1$ (on all of $K$) \st $h' |Z \gg  g|Z$. Set $h'(x_i) = \beta_i' < \beta_i$.

We want to find affine $g' \gg 0$ (on $K$) \st $g'|Z = g|Z$ and $h' \geq g'$. Affine functions are determined and determine uniquely their values on the vertices, so if we set $g(v_i) = y_i > 0$, we must have $y_t  = \alpha_t g'(x_t) + (1-\alpha_t)g'(x_0)$,  if $1 \leq t \leq j$, and from $y_j = g'(x_i)  < h'(v_t) $ for $t > j$. In order to obtain $g' \leq h'$, we must have $g'(x_i ) \leq h'(x_i)$.

The first set of equations requires $g'(x_t) = (y_t  - (1-\alpha_t)g'(x_0))/\alpha_t$ (in particular,  $g'(x_0) < y_t/(1-\alpha_t), h'(x_0)$ ensures positivity), and the second set (when $t > j$) $g'(x_j < h'(x_j) $. 

Now we observe that
$$
\max{\frac{y_t- \alpha_t\beta_t'}{1-\alpha_t}} \leq  \min \brcs{\frac{y_i}{1-\alpha_i}}
$$
To see this, note that the left side is just $\beta_0' = h'(x_0)$ (all the values are the same!), and so this inequality reduces to $(1-\alpha_i) \beta_0' \leq    y_i = (1-\alpha_i) \beta_0' + \alpha_i h(x_i)$. Let $\delta$ be the minimum of the right side (so $\delta > 0$).

Define $g'(x_0) = \delta$;
for $1\leq t \leq j$,  set $g'(x_t) = (y_t - (1-\alpha_t)\delta)/\alpha_t$, so that $g(v_t) = y_t =g'(y)$ for this set of $t$s,  and $g'(x_t) > 0$. Finally for $t > j$, set $g'(x_t) = y_t j= g(x_t)$. Then $g' \gg 0$, $g'|Z = g|Z$, and if $t > j$, $h(x_t)> g(x_j) = g'(x_t)$. Now for $1 \leq t \leq j$, $h'(x_t) \geq g'(x_t)$ follows from the definition of $\delta$ and the displayed inequality. Since $h \gg h'$, $h \gg g'$, and we are done. \qed

\Lem Corollary \appbele. Suppose $K$ is a finite dimensional simplex, and ${\Cal L}$ is a flat \st $Z:= K \cap {\Cal L}$ is a simplex. Suppose there exists a sequence of flats ${\Cal L} = {\Cal L_0} \subset {\Cal L}_1 \subset \dots \subset {\Cal L}_k$ \st each $Z_i: = K \cap {\Cal L}_i$ is a  simplex of one dimension less than $Z_{i+1}$, and $Z_k$ is of codimension one in $K$. Then $Z$ is good. 

\Pf We may apply the lemma to each inclusion $Z_i = {\Cal L}_i \cap Z_{i+1}$, and so extend the affine function from one simplex to the next. \qed

Since it is routine to verify that if the dimension of $K$ is less than or equal three (that is, four or fewer vertices), then the sequence of flats always exists. We thus deduce that if $(G,u)$ is a simple dimension group with at most four pure traces, then for a set of traces, $U$, $Z(\ker U)$ is good if and only if $Z(\ker U)$ is a simplex and $\ker U$ has dense range in $Z^{\vdash}$. 

If $\Cal L$ is a {\it translate\/} of the affine span of a face (not necessarily a facet) of a finite dimensional simplex $K$, then $Z = K \cap {\Cal L}$ admits such a sequence and therefore is good \wrt $\Aff K$.
Whether every (finite dimensional) $Z = K \cap {\Cal L}$ with $Z$ a simplex admits such a chain is open. To prove the existence of the chain, it would be enough to show that given $Z$, there exists a flat ${\Cal L}'$ containing $\Cal L$ \st $Z$ is strictly contained in $K \cap {\Cal L}' $, which in turn is not $K$. 
\endcomment

\Lem Corollary \sevten. Suppose that $(R,1)$ is a partially ordered ring that is a simple dimension group, with $1$ as an order unit. Then each  pure trace is refinable. 

\Remark 1 Returning to Example \exnin, $R = \Q[x]$ with the strict ordering obtained from restriction to an interval $I \subseteq [a,b]$, we saw that the pure trace $\phi_{\alpha}$ ($\alpha \in I$) is good if and only if $\alpha$ is algebraic and no other algebraic conjugates of $\alpha$ are in $I$, and $\phi_{\alpha}$ is bad if and only if $\alpha$ is transcendental. Both of these types of pure traces are obviously refinable, and the remaining pure traces ($\phi_\alpha$ where $\alpha$ is algebraic and at least one other algebraic conjugate is in $I$) are neither good nor bad (nor ugly), but are refinable by this result. 

\Remark 2 The conclusion does not hold if we weaken the hypotheses to some consequences, for example that $(G,u)$ be a simple dimension group for which $S(G,u)$ is a Bauer simplex or even finite-dimensional. The simple dimension group $G'$ of Example \exten\  (below) has a pure trace which is not refinable (nor ugly), but the pure trace space consists of three points. On the other hand, for simple dimension groups with exactly two traces, modulo the infinitesimals, every pure trace is either good,  bad (hence refinable), or ugly (hence not refinable). 

\Pf If $R = \Z$, the only trace is obviously good, hence refinable. Otherwise, $R$ has dense image in its affine representation. We show by routine arguments that $Z \equiv Z(\ker \tau)$ is a face of $S(R,1)$ and then the image of $\ker \tau$ is dense in $Z(\ker \tau)^{\perp}$. Then by  Proposition \appbtwo, $Z$ is order unit good. As $R$ is simple as  a dimension group, $Z$ is good; therefore $Z$ is refinable, and thus $\tau$ is refinable. 

We may assume $\tau$ is normalized at $1$. Since $\tau$ is pure, it is multiplicative (\throne), hence $J:= \ker \tau$ is an ideal. Let $X = \partial_e S(R,1)$; this is a compact space (since $R$ is a partially ordered ring with $1$ as order unit), and we can identify $C(X,\R)$ with $\text{Aff\,}(S(R,1))$ (and the probability measure space on $X$ with $S(R,1)$). The closure of $J$ is a closed ideal of the closure of $R$, that is, of $C(X,\R)$. Hence there exists compact $Y\subset X$ \st $\overline J = \Ann Y: = \Set{f \in C(X,R)}{f|Y \equiv 0}$. There exists a unique face $F$ of $S(R,1)$ \st $Y  = F \cap X$ (always true for Bauer simplices), and it is easy to check that $\Ann Y = F^{\perp}$. 

Now $Z(\tau)$ consists of normalized traces (that is, probability measures on $X$) that kill $\ker \tau$, so that $F \subseteq Z$. On the other hand, if $\sigma \in Z(\tau)$, then $\sigma(J) = 0$, so $\sigma \in F^{\perp}$. 
\qed

\Lem Example \exten. Simple dimension groups with three traces; pure traces which are refinable but not good, etc.

\noindent Let $\alpha$ and $\beta$ be  positive irrational numbers, and define the following four elements of $\R^3$.
$$\eqalign{
v_1 & = \(\matrix 0 & 0 & 1\\ \endmatrix\)\cr
v_1' & = \(\matrix 0 & \beta & 1\\ \endmatrix\)\cr
v_2 & = \(\matrix \alpha & 1 & 0\\ \endmatrix\)\cr
v_3 & = \(\matrix 1 & 0 & 0\\ \endmatrix\).\cr
}$$
Let $G $ be the rational vector space spanned by $\brcs{v_1, v_2, v_3}$ and let $G'$ be the rational vector space spanned by $\brcs{v_1', v_2, v_3}$. Let $\tau_i$ denote the projections on the $i$th coordinate. Each of $G$ and $G'$ is dense in $\R^3$, so equipped with the strict ordering inherited from $\R^3$, each is a simple dimension group, and $\brcs{\tau_i}$ is the set of pure traces (in both cases). 

For $G$, we consider $U = \brcs{\tau= \tau_1}$. It is a pure trace whose kernel is spanned (as a rational vector space) by $v_1$, and since $\tau_2(v_1) = 0$, it follows that $\ker \tau$ is not dense in $(\ker \tau)^{\vdash} \iso \R^2$. Hence $ \tau$ is not good; neither is it bad, since $v_1 \not\in \Inf(G) = \brcs{0}$. On the other hand, $Z(\ker U)$ is the face spanned by $\tau_1$ and $\tau_2$, and is thus good (\wrt $\Aff S(G,u)$). From \appbnin, $Z(\ker U)$ is good (\wrt $G$), therefore refinable, and thus $\tau$ is refinable; this  also follows from an (easy) construction of coproducts below. In this case, $\tilde U = \brcs{\tau_1} \neq Z(\ker U)$.

In contrast, for $G'$, again let $\tau = \tau_1$. Its kernel is $v_1'\Q$, so this is a pure trace that is neither bad nor good. However, it isn't refinable either, because $\tau_2-\beta\tau_3$ kills $v_1'$ (and thus $\ker \tau$), but there is no trace $\sigma = r\tau_2 + s\tau_3$ (where $r$ and $s$ are nonnegative real numbers) that kills $v_1'$, so that $Z(\ker \tau) = \brcs{\tau}$ (in this case, $Z = Z(\ker \tau)$ is a face), and thus $Z^{\vdash}$ is two dimensional, but $\ker \tau$ is rank one, so the image of $\ker \tau$ is not dense in $Z^{\perp}$. (Since  $G'$ is a rational vector space, there is no difference between $Z(\ker U)$ being good or refinable.)

For both groups, $\tau_2$ is refinable but neither good nor bad, and $\tau_3$ is good (by the density criterion).  \qed

An easy sufficient condition for refinability comes from the following construction, similar to the one constructing good sets out of the building blocks of faces and singletons. Just as the ordered tensor product is the product within the category of simple dimension groups with order unit (and also in the category of dimension groups with order unit), we can find a coproduct in this category (which differs from the coproduct in dimension groups, that is, the direct sum), provided we exclude the cyclic case (ordered groups isomorphic to $\Z$). Let $G$ and $H$ be noncyclic simple dimension groups. Form the direct sum, $\Cal C \equiv \Cal C(G,H)= G \oplus H$, and impose on it the strict ordering, that is, 
$$
\Cal C^+ \setminus \brcs{(0,0)} = \Set{(g,h) \in G \oplus H} {g \in G^{++}, h \in H^{++}}.
$$
It is easy to verify that $\Cal C$ is simple and that the Riesz interpolation and unperforation properties hold (the former uses the strict interpolation property that holds for simple noncyclic dimension groups); thus $\Cal C$ is a simple dimension group. The  pure trace space of $\Cal C$ is the disjoint union of the pure trace space of $G$ and of $H$ (the pure traces are precisely those of the form $\tau \oplus 0$ and $0 \oplus \sigma$ where $\tau$ varies over the pure traces of $G$ and $\sigma$ varies over the pure traces of $H$). Note however, the inclusions $G \to G \oplus \brcs{0} \subset C$ and $H \to  \brcs{0}  \oplus H \subset C$ are {\it not\/} order preserving. Instead, the projection maps $\Cal C \to G$ and $\Cal C\to H$ are. It is routine to verify that within the category of simple noncyclic dimension groups, this is the coproduct. [We have to exclude $\Z$ since, as is well known, no noncyclic simple dimension group can have a discrete trace, which is what would happen if we allow one of $G$ or $H$ to be cyclic---so the constructed group would not be a dimension group.] 

This is a special case of a pullback of simple dimension groups (in [M], Matui discusses pushouts of dimension groups; these are dual to pullbacks). In this more general setting, we have two positive onto maps $\Arrow \phi; G. K$ and $\Arrow \psi; H. K$ where $K$ is a dimension group (possibly zero), and we form ${\Cal C}=\Set{(g,h)}{\phi(g) = \psi(h)} \subseteq G \oplus H$, with the strict ordering induced by the inclusion. For the resulting subgroup to be a dimension group, we require that if $\rho$ is a pure trace of $K$, then  $\rho \circ \phi$ and $\rho \circ \psi$ are pure traces of their respective groups, or zero. The extremal trace space is then obtained by taking the union and identifying the traces that factor through the map to $K$. If $K$ is zero, we obtain ${\Cal C(G,H)}$.

\Lem Lemma \seveig. Suppose  $G$ and $H$ are noncyclic simple dimension groups, and $\tau$ is a pure refinable (good) trace of $G$. Then the extension of $\tau$ to a pure trace of $\Cal C (G,H)$ given by $\tau \oplus 0$ (that is, $(g,h) \mapsto \tau(g)$) is refinable (good). 

\Pf Fix an order unit $u$ (\wrt which $\tau$ is refinable) of $G$, and let $v$ be any order unit of $H$. Let $T $ denote $\tau \oplus 0$. Obviously $T([(0,0),(u,v)]) = \tau([0,u])$ and if nonzero $a \in [0,u]\setminus\brcs{u}$ and we pick any element $v' \in [0,v]\setminus{\brcs {0,v}}$ (such exist since $H$ is noncyclic), then $[a,v'] \in [(0,0),(u,v)]$. It easily follows that refinability (goodness) of $T$ holds. \qed

Referring back to Example \exten, we can quickly  show (again) that the trace $ \tau_1$ of $G$ is refinable. Let $G_0 = \Q v_3 + \Q v_2$; equipped with the strict ordering as a subgroup of $\R^2$  (i.e., stripping the third coordinate), $G_0$ is a simple dimension group and  the first coordinate projection, $p$, is a pure trace that is bad (since $\alpha$ is irrational), hence refinable. Now $G$ may be identified in the obvious way with $ \Cal C(G_0 ,\Q)$ (where $v_3$ yields the rational direct summand) and $p \oplus 0 = \tau_1$, so the latter is refinable, by Lemma \seveig.

\SecT Appendix A\ \  Tensor products dynamically

If $(X,S)$ is a minimal zero-dimensional dynamical systems (that is, $X$ a non-atomic separable compact zero-dimensional space---a Cantor set---and $S$ is a self-homeomorphism), then the partially ordered group K$_{0} ( C(X,\C) \times_{S} \Z) := \text{K}_0(X,S)$ has a natural pre-ordering (as all K$_0$-groups do), which makes it into a simple dimension group, alternatively described by the quotient ordering on the cokernel of $\Arrow \I-S; C(X,\Z).C(X,\Z)$; the image of the constant function $\pmb1$ is the order unit \wrt which traces are normalized, and there is a natural  affine homeomorphism between the invariant probability measures on $X$ and the normalized traces on the dimension group. Conversely, every simple dimension group arises in this fashion, even when include a specific order unit as part of the data. The simple dimension group is a complete invariant for strong orbit equivalence. This and much more  is due to [GPS1].

If we take two dimension groups with order unit, $(G,u)$ and $(H,v)$, it is a natural construction to form the ordered tensor product, as in [GH2]---on the tensor product over $\Z$ (that is, as abelian groups), $G \otimes_{\Z} H$ (the subscript $\Z$ is usually deleted), impose the positive cone generated by all pure tensors, $g \otimes h$ where $g \in G^+$ and $h \in H^+$. That this {\it is\/} a partially ordered group is not entirely obvious in general, but for dimension groups, it is even a dimension group (from [EHS]), with $u \otimes v$ as designated order unit. If both $G$ and $H$ are simple, then so is $G\otimes H$. 

The tensor product construction on dimension groups is perfectly natural, so at least in the case that the dimension groups are simple, should correspond to some construction on the minimal dynamical systems that yielded them. Explicitly, if $(X,S)$ and $(Y,T)$ are minimal zero-dimensional dynamical systems corresponding respectively to simple dimension groups with order units, $(G,u)$ and $(H,v)$, then what minimal action does $(G\otimes H, u\otimes v)$ correspond to?

We  answer this; however, the answer is not entirely satisfactory, since it only yields one of the many possible minimal $\Z$-actions on $X \times Y$ (necessarily orbit equivalent to the obvious $\Z^2$ action), and there is no obvious way to single out a natural one.  

We do this by exploiting (as is done in [GPS1]) Vershik's adic maps. 

Pick a Bratteli diagram that represents $(G,u)$ (typically, at level zero, these have a single point at the apex). Let $X_n$ be the space of paths of length $n$ beginning at the top. This is a finite space, and we may form the path space, $X = \lim_{\leftarrow} X_{n}$, where the maps $\Arrow \pi_{n}; X_{n+1} .X_n$ are given by simply removing the last vertex and edge. Two paths are {\it tail equivalent\/} if they are composed of the same edges, up to finitely many differences. 

On each $X_n$, we define an equivalence relation $R_n$ simply by declaring two paths equivalent if they end at the same vertex. On each $R_n$-equivalence class, we try to impose a total ordering, $\prec_n$, in such a way that if the finite paths, given as sequences of edges, $p:= (e_1,e_2,\dots,e_n) \prec_n q=(f_1,f_2,\dots,f_n)$, then for every admissible edge $g_{n+1}$, the new paths satisfy $(p,g_{n+1}) \prec_{n+1} (q,g_{n+1})$. If this is done, we can compare any two infinite paths that are tail equivalent to each other, inducing a partial ordering $\prec$ on $X$, and it is easy to check that if $p$ is an infinite path which is not maximal, it has an immediate successor, i.e., there exists $p' $ \st $p \prec p'$ and if $p \prec p''$ with $p \neq p''$, then $p' \prec p''$. 

Let ${\Cal M^+}$ and ${\Cal M^-}$ denote the collections of maximal and minimal paths respectively; both are compact. Vershik's adic map is given by  $p \mapsto p'$, defined on $X \setminus {\Cal M^+}$ with range in $X\setminus {\Cal M^-}$. 
As  observed in [GPS1], if $\Cal M^+$ and $\Cal M^-$  each consist of a single point, then the adic map can be extended to a minimal homeomorphism of $X$  by sending the maximal path to the minimal one. 

To construct a minimal dynamical system that yields the tensor product of the dimension groups, we deal with the orderings a bit more generally. Let $(W, \prec)$ be a partially ordered set; an equivalence relation $R$ on $W$ is said to be {\it compatible with $\prec$} if for elements $x$, $y$, and $z$, whenever $yRz$ and either $x \prec y$ or $y \prec x$, then $xRz$. For example, if $R$ is the trivial equivalence relation ($xRy$ for all $x$ and $y$), this obviously holds, but it also holds if $R$ is the transitive closure of the relation given by $x \prec y$ or $y \prec x$. In fact, the latter is the minimal equivalence relation compatible with $\prec$. (It may turn out to be trivial, e.g., if $W$ has a  minimal element less than or equal to every element of $W$.) Even better, the  relation $x \prec y$ or $y \prec x$ could turn out to be an equivalence relation, and this is what happens in the motivating example. 

If $(W,\prec,R)$ and $(Z,\prec',S)$ are partially ordered sets with respective compatible equivalence relations, then we can construct a new partially ordered set (with a compatible equivalence relation), which we call the {\it $R$-lexicographic product,} defined as follows. The underlying set is $W\times Z$, the equivalence relation is $R \times S$, and the ordering is given by 
$$
(w,z) \prec'' (w_1,z_1) \text{ when } \cases  z \prec' z_1& \text{if $z \neq z_1 \  \&\   wRw_1$, or}\\
w \prec w_1 & \text{if $z = z_1$}\\
\endcases
$$
It is routine to verify that $W \times Z$ becomes a poset (this does not require the equivalence relation $S$), and moreover, $R \times S$ is compatible \wrt $\prec''$. We abbreviate this new object, $W\times_{R\text{-lex}} Z$. If $R$ is the trivial equivalence relation, then this is the usual lexicographic product. 

In our situation, arising from Bratteli diagrams, we simply take $X_n \times_{R_n\text{-lex}} Y_n$, and then the new ordered compactum is $\lim_{\leftarrow} X_n \times_{R_n\text{-lex}} Y_n$, with the products of the usual maps. Since inverse limits commute with finite direct products, the underlying toplogical space is just $X \times Y$; however, the ordering does not commute with inverse limits.

Now it is routine to check that the transition matrices (that is, obtained from 
$X_{n+1} \times_{R_{n+1}\text{-lex}} Y_{n+1} \to X_n \times_{R_n\text{-lex}} Y_n $ are exactly the matrix tensor products of the transition matrices; this is, if $\pi_{n}$ is induced by $M_n$ (a rectangular matrix with entries from the positive integers) and $N_n$ for $X_{n+1}\to X_n$ and $Y_{n+1}\to Y_n$ respectively, then the transition matrix on the new construction is simply $M_n \otimes N_n$ (the Kronecker product, in old-fashioned terminology). 

Again it is routine to verify that the unique successor property holds, and so the adic map is defined on an open subset. Furthermore, the collection of maximal paths of the new thing is just the cartesian product, $\Cal M_X^+ \times \Cal M_Y^+$, and the same applies to the minimal paths. Hence if each ordered Bratteli diagram just has one maximal and one minimal path, then so does the new object, and thus we create  a minimal homeomorphism, call it $S\otimes T$, on an ordered Bratteli diagram, whose dimension group is the tensor product of the original dimension groups. Unfortunately, this is not canonical, in that we make a lot of choices on how to order the product spaces; we picked the $R_n$-lexicographical order because it seemed reasonable. 

The effect of the adic map in this construction is to sometimes choose $S\times 1$ and other times choose $1 \times T$, acting on a path. Hence to each element of $X \times Y$ we can associate an element of $2^{\N}$; its coordinate in the $n$th position will be a $0$  if the $n$th iterate of the path under the $\Z$ action is given by $S\times 1$, otherwise the entry will be one. This gives a factor map from $S\otimes T $ to a subshift of the  two-shift.  

It is known that any product action of a finite number of minimal homeomorphisms is hyperfinite, that is it is orbit equivalent to a single homeomorphism. The following observation
exhibits this minimal homeomorphism explicitly.
Let $B$ and
$B'$ two Bratteli diagrams corresponding to simple dimension groups $G$ and $G'$ respectively.  Let $S $ and $T$ be minimal Vershik homeomorphisms
defined on the path spaces of $B$ and $B'$. Then the sets of ergodic (that is, extremal) invariant measures, $M(X, S)$
and $M(X', T)$ are naturally in bijection wtih  $\partial_e S(G, u)$ and $\partial_e S(G', u')$ respectively.
Construct the Bratteli diagram $B \otimes B'$ and the Vershik map $S \otimes T$ acting on its path
space as above. It is well known that the product
action of $\Z^2$  on $X \times X'$ generated by $S \times  \text{I}$ and $\text{I} \times  T$ is orbit equivalent to a $\Z$-action
by a single minimal homeomorphism, say $V$. Then $V$ and $S \otimes T$ have the same sets
of invariant ergodic probability measures, namely $M(X, S) \times M(X', T) \iso \partial_e S(G \otimes  G', u \otimes  u')$.
Hence $S \otimes T$ is orbit equivalent to $V$, and therefore $S \otimes T$ is orbit equivalent to the
$\Z^2$-action.
This also applies to product $\Z^n$-actions, by induction.


The following is presumably well known; its proof is standard in this context.

\Lem Lemma \appone. Let $(X,S)$ and $(Y,T)$ be compact zero-dimensional dynamical systems, and let $\nu$ be a probability measure on $X \times Y$  extremal among measures invariant under the action of $\Z^2$ generated by $S \times 1$ and $1 \times T$. Then there exist an extremal $S$-invariant probability measure $\mu_1$ on $X$ and an  extremal $T$-invariant probability measure $\mu_2 $ on $Y$ \st $\nu = \mu_1 \times \mu_2$ (the product measure).

\Pf Define $\mu_1$ and $\mu_2$, measures on $X$ and $Y$ respectively, via $\mu_1 (A) = \nu (A \times Y) $ and $\mu_2 (B) = \nu(X \times B)$. These are invariant under $S$ and $T$ respectively, and their extremality (within the set of invariant measures) follows from the corresponding property of $\nu$.

For each clopen set $B$ of $Y$ \st $q:=\nu(X \times B) \neq 0$ (that is, $\mu_2 (B) \neq 0$), define a measure on $X$, $\mu_B$, via $\mu_B (A) = \nu (A \times B)/\mu_2 (B)$ (for $A$ an arbitrary clopen of $X$). It is routine  that this is an $S$-invariant measure, and $\mu_B (A) \leq \frac 1q \mu_1 (A)$ for all clopen sets $A$, and thus for all elements in the $\sigma$-algebra they generate. By extremality of $\mu_1$, we must have $\mu_B  = \mu_1$. Hence $\nu(A \times B)/\mu_2 (B) = \mu_1 (A)$, that is $\nu(A \times B) = \mu_1 (A) \mu_2 (B)$ for all clopen sets $B$ \st $\mu_2 (B) \neq 0$.

 If $\mu_2 (B) = 0$, then $\nu (A \times B) \leq \nu(X \times B) = 0$;  thus for all clopen sets $B$ in $Y$, we have $\nu (A \times B) = \mu_1 (A) \mu_2 (B)$. Let $\nu' $ be the product measure on $X \times Y$ obtained from $\mu_1$ and $\mu_2$. Then $\nu$ agrees with $ \nu'$ on rectangles; since the rectangles generate the full $\sigma$-algebra generated by the clopen sets in the product, we must have $\nu = \nu' = \mu_1 \times \mu_2$.
\qed

\Lem Corollary \apptwo. Let $(X,S)$ and $(Y,T)$ be zero dimensional compact dynamical systems, and form the space $(X \times Y, \langle S\times1, 1 \times  T\rangle)$ (the natural action of $\Z^2$---not of $\Z$ (via $S\times T$)---on $X \times Y$). If $\alpha$ and $\beta$ are  pure invariant measures  on $X \times Y$ whose corresponding traces agree on K$_0\(C(X \times Y)\times_{S,T} \Z^2\)$, then $\alpha = \beta$.

\Pf Since $\alpha$ and $\beta$ are pure invariant measures, by the preceding, there exist $\mu_i$, $\nu_i$, pure invariant measures on the relevant sets \st $\alpha = \mu_1 \times \nu_1$ and $\beta = \mu_2 \times \nu_2$. If $\mu_1 \neq \mu_2$, then there exists a clopen $A$ of $X$ \st $\mu_1 (A) \neq \mu_2 (A)$, so that $\mu_1 (\chi_A) \neq \mu_2 (\chi_A)$, and thus the corresponding traces on K$_0 (C(X )\times_S \Z)$ have different values at the equivalence class $[\chi_A]$. The natural map $\text{K}_0 (C(X )\times_S \Z)  \otimes \text{K}_0 (C(Y )\times_T \Z) \to \text{K}_0 (C(X\times Y )\times_{\langle S,T\rangle } \Z^2) $ is order preserving (but not onto---there is a free rank one contribution from K${}_1 \otimes \text{K}_1$), and $\alpha(A \times Y) \neq \beta (A \times Y)$ entails the corresponding traces on $\text{K}_0$ of the $\Z^2$ crossed product disagree on $[\chi_{A \times Y}]$. Hence $\mu_1 = \mu_2$, and similarly, $\nu_1 = \nu_2$, so that $\alpha = \beta$. \qed

\SecT Appendix B Unperforation of quotients

Unperforation of quotients by convex subgroups (of the form $\ker U$) seems to be a significant contributor to the difference between goodness (of sets of traces) and refinability, as in  Proposition \appbsix(b \& f), so it seemed worthwhile to investigate this in more detail. First, we give  sufficient conditions for unperforation, generalizing from the divisible case. Then we show that many quotients $G/\ker U$ are very badly not unperforated. 

Here is a sufficient condition for $G/L$ to be unperforated, where $L$ is a convex subgroup ofÊ $G$ \st $G/L$ is torsion free. It is likely necessary as well in this context (simple dimension groups). 

\Lem Proposition \sevthr. Let $(G,u)$ be a simple dimension group, and let $L$ be a convex subgroup of $G$ \st $G/L$ is torsion free (as an abelian group). Sufficient for $G/L$ to be unperforated \wrt the quotient ordering is the following:
{\par} \item{(*)} For all $k \in L$, for all $\epsilon > 0$, there exist $k_1, k_2, k_3 \in L$ \stÊ 
$$
\hat k = 2\hat k_1 + \hat k_2 -\hat k_3ÊÊÊ \qquad\text{and }\qquadÊ \hat k_2, \hat k_3 \geq -\epsilon \pmb 1.
$$

\Remark The coefficient $2$ ofÊ $k_1$ can obviously be replaced by any integer exceeding $1$.

\Remark The ostensibly weaker condition, 
\item{($\frac *2$)} For all $k \in L$, for all $\epsilon > 0$, there exist $k_1, k_2, k_3 \in L$ \st 
$$
\left\|\hat k - (2\hat k_1 + \hat k_2 -\hat k_3) \right\|Ê < \frac {\epsilon}2Ê \qquad\text{and }\qquadÊ \hat k_2, \hat k_3 \geq -\frac{\epsilon}2 \pmb 1
$$

\noindent 
is equivalent to (*), since we can incorporate the difference, $k - (2k_1 + k_2 - k_3)$, into either $k_2$ or $k_3$.

\Pf Since $G$ is simple, so is $G/L$ (with the quotient ordering); thus for $a$ in $G/L$, $P(a)$ (defined as before as $\Set { n\in \N}{na \in (G/L)^{++}}$) is either empty, all ofÊ $\N$ (so $a \geq 0$), or cofinite. If the last holds, there exists a positive integer $l$ \st $2^l a \geq 0$. Hence it suffices to show that ifÊ $a$ is in $G/L$ andÊ $2a \geq 0$, then $a \geq 0$.

To this end, suppose $a = g +L$ and $2a \geq 0$. Then there exist $k \in L$ and $p \in G^{++}$ \st $2g + k = p$. Select $\epsilon < \inf_{\sigma \in S(G,u)} \sigma(p)/2$, and apply (*) to $k$. Then $2g + 2k_1 + k_2 = p + k_3$; define $h = g + k_1 + k_2$ so thatÊ $2\hat h = 2(\hat g + \hat k_1 + \hat k_2) = \hat p + \hat k_2 + \hat k_3$. Hence $2\hat h \geq \hat p - 2\epsilon \pmb 1 \gg 0$. Thus $\hat h \gg 0$, so that $h$ is an order unit, and $h + L = g + L= a$. So $a \in (G/L)^{++}$.\qed

If $J$ is a closed subgroup of a Banach space, then the divisible subgroup of $J$ (the maximal rational subspace of $J$), $D(J)$, is a closed real vector space and the maximal real vector subspace of $J$. If $J$ is the closure of an abelian group of finite rank or more generally, if the vector subspace generated by $J$ (denoted $J\R$) is finite dimensional, then $D(J)$ is  also open, but this does not extend to the infinite rank case. For $J$ a closed subgroup of $\Aff K$, let $\Cal P(J)$ be the subgroup of $J$ generated by $\Set{j \in J}{j \geq 0}$. An elementary condition for $J$ to be a real vector space (that is, $ J = D(J)$) when $J $ arises as  the closure of the image of $\ker U$ (for $U $ a subset of $S(G,u)$),  is that for all $k \in \ker U$, the element $\hat k/2 \in \Aff S(G,u)$ can be uniformly approximated by   images of elements of $\ker U$. 

If $J= \Cal P(J)$, then $Z(\ker J)$ is a closed face of $K$. If we take critical simple dimension groups of the form $G = \langle e_i; \sum \alpha_i e_i \rangle$ (where $\brcs{1, \alpha_1, \alpha_2,\dots }$ is linearly independent over the rationals), and set $U = \brcs{\tau_1}$ (consisting of the the trace obtained from the projection on the first coordinate ), then $\ker U = \langle e_i \rangle_{i\neq 1}$. Each $\hat e_i = e_i$ (a consequence of the notation), so is nonnegative as a function, and thus $J = \Cal P(J)$ and it follows that the quotient group is $\Z + \alpha_1 \Z$ with the ordering inherited from the reals.

On the other hand, if we alter this somewhat, there is a drastic change. Take $G = \langle e_1,e_2; (\alpha,\beta)\rangle$ with $\brcs{1,\alpha,\beta}$ linearly independent over the rationals, so that $G$ is dense in $\R^2$. Set $\tau = \tau_1 + \tau_2:(a,b)\mapsto a+b$; then $\ker \tau = (1,-1)\Z$. It turns out that $H = G/\ker \tau$
 is very holey. 

For example, we can construct  for every positive integer $m$ an element $a_m$ in $H$ \st $P(a_m) = \N \setminus \brcs{1,2,3,\dots, m}$ (that is, $n a_m \in G^+$ if and only if $n \geq m+1$). To see this, by density of $G$ in $\R^2$, we may find $A_m = (r,s) \in G$ \st $1/(m+1) < r < 1/m$ and $-1/(2m+1) < s < -1/(2m+2)$. Then $(m+1)A_m + (1,-1) \gg 0$, but for no $j < m+1$ does there exist an integer $n$ (positive or negative)
\st $jA_m + (n,-n)\gg 0$. Hence if $a_m = A_m + \ker \tau$, then $P(a_m) \cap \brcs{1,2,\dots, m+1} = \brcs{m+1}$, and   if $l > m$, then $l \in P(a_m)$.

 It is  easy to describe the ordering on the quotient: $(r,s)$ has positive image if and only if both $r+s > 0$ and the interval $[-s,r]$ contains an integer.  In particular, if $r+s> 1$, then $(r,s)$ has positive image in $G/\ker \tau$. It is plausible that as $a$ varies over the elements of $H$, the corresponding $P(a)$ run over all cofinite subsemigroups of $\N$. It would be interesting to decide whether interpolation holds for this example. 

\Lem Corollary \sevfiv. Let $L$ be a subgroup of the simple dimension group $(G,u)$ \st $L \cap G^+ = \brcs{0}$ and $G/L$ is torsion-free. Let $J$ denote the closure ofÊ the image of $L$ in $ \Aff S(G,u)$. Either of the following is sufficient for $G/L$ to be unperforated \wrt the quotient ordering.
{\par} \item{(a)} $J$ is a real vector space 
\item{(b)} $D(J) + \Cal P(J)$ is dense in $J$.

\Pf In both cases, ($\frac *2$) applies; in the first case,Ê $k_2 = k_3 = 0$.\qed

Unbelievably, there is an example wherein  Proposition \sevthr\ applies but Corollary \sevfiv\ does not. This  was originally constructed as a potential example of a singleton set $U$ that was refinable and not good with  $U = Z(\ker U)$ (no such examples are known); however, it turned out not to be refinable, as $G/\ker U$ is unperforated but $\ker U$ does not have dense range in $Z(\ker U)^{\vdash}$.

\Lem Example \exthi. A simple dimension group with four pure traces together with a singleton set $U =\brcs{\tau}$ \st $G/\ker U$Ê is unperforated, but does not satisfy either sufficient condition of Corollary \sevfiv.

\noindent Let $\brcs{a,b,c,r,s}$ be an algebraically independent set of real numbers (this is likelyÊ overkill). Define the following six elements ofÊ $\R^4$.
$$
\matrix 
Êv_1 & = & (0 & 0 & 1 & -1) \\
Ê v_2 & = & (0 & 0 & s & -s) \\ 
Êv_3 & = & (1 & -1 & a & 0) \\ 
Êv_4 & = & (r & -r & 0 & b) \\ 
Êv_5 & = & (1 & 0 & 1 & 1) \\ 
Êv_6 & = & (c &1 & 0 & 0) \\ 
\endmatrix
$$
Let $G$ be the subgroup ofÊ $\R^4$ generated by $\brcs{v_i}$. Let $\tau_j$ ($ j = 1,2,3,4$) be the coordinate projections, and set $\tau = \tau_1 + \tau_2$. Then $\tau^{\vdash} = \Set{(x,-x,y,z)}{x,y,z \in \R}$. 
We first observe that $\ker \tau = \langle v_1, v_2, v_3, v_4 \rangle$ (an easy consequence of linear independence ofÊ $\brcs{1,r,c}$ over the rationals). Now let $J$ be the closure of $\ker \tau$. We immediately see that $v_1\R \subset J$ (from linear indepence ofÊ $\brcs{1,s}$ over $\Q$). Thus $J$ is the closure of the span of $\brcs{v_1 \R, v_3, v_4':= (r,-r,-b,0)}$. We note that $\brcs{v_1, v_3, v_4}$ is a real basis forÊ $\tau^{\vdash}$. Hence $J = v_1\R + v_3\Z + v_4'\Z$, expressed as a direct sum of a vector space and the discrete abelian group $\langle v_3, v_4 \rangle$. In particular, $J \iso \R \oplus \Z^2$. 

Now we show $G$ is dense in $\R^4$. Since $J$ is in the closure of $G$, we can add elements of $v_1 \R$ to the generators of $G$ without affecting its closure. This shows that  the elements $v_4'$ and $ v_5' = (1, 0, 2, 0)$  are in $J$. Now consider the subgroup of $\R^3 \times \brcs{0}$ generated by $v_3, v_4', v_5', v_6$; delete the terminal zero of each, and calculate all four $3 \times 3$ determinants arising from this set of four vectors. The result is a linearly independent set over $\Q$ (a consequence of algebraic independence of $\brcs{a,b,c,r}$). Hence this subgroup is dense in $\R^3 \times \brcs{0}$. Since $v_1\R + (\R^3 \times \brcs{0}) = \R^4$, $G$ is dense in $\R^4$. 

Hence $G$, equipped with the strict ordering, is a simple dimension group. Let $U = \brcs{\tau}$, so $\ker U = \ker \tau$. Now $\Cal P(J) = \brcs{0}$: this is obvious from the decomposition $J = v_1 \R \oplus v_3\Z \oplus v_4'\Z$, and since $D(J) = v_1 \R$, the criterion in the corollary cannot be applied. We show  that the hypothesis of Proposition \sevthr\ does apply.

Select $k \in \kerÊ \tau$. IfÊ $k$ has aÊ nonzero component in $\langle v_1, v_2\rangle$, we can immediately reduce to the zero component there. Hence we can write $k = tv_3 + uv_4 = (ta + ur, -(ta + ur), tu, ub)$ for integers $t,u$. Since $\brcs{1,r}$ is linearly independent over $\Q$, there exist integers $t', u'$ \st $|(t+ ur)/2 - (t' + u'r)| < \epsilon/8$.Ê Set $k' = t'v_3 + u'v_4$,Ê so that $k - 2k' = (2\delta, -2\delta, (t - 2t')a, (u - 2u')b)$ where $|\delta| < \epsilon/8$. Next, we may find integers $m,n$ \st $|m + ns - (t -2t')a/2| < \epsilon/8$. Set $k'' = mv_1 + nv_2$, so that 
$k - 2k' - 2k'' = (2\delta, -2\delta, 2 \eta, (u-2u')b -2(m-ns)):= k'''$ where $|\eta| < \epsilon/8$. 

Let $f$ denote the fourth coordinate of $k'''$. Set $k_1 = k' + k''$. If $fÊ \geq 0$, then $k''' \geq \min\brcs{-2|\delta|, -2|\eta|}$, so we can set $k_2 = k - 2k' - 2k''$ and $k_3 = 0$; if $f' < 0$, set $k_3 = -k'''$ and $k_2 = 0$. 

Hence $G/\ker \tau$ is unperforated. It is easy to see that $Z(\ker \tau) = \brcs{\tau}$, so $G/\ker \tau$ has unique trace, that induced by $\tau$. Being simple, unperforated, and with unique trace, $G/\ker \tau$ is order isomorphic to the subgroup $\tau(G) $ of the reals (with the relative ordering), in particular, $G/\ker \tau \iso \Z + c\Z \subset \R$.

In this example, and since $J \neq \tau^{\vdash}$ (as $J$ is not a real vector space), $\tau$ is not good; since $G/\ker \tau$ is unperforated, $Z(\ker \tau) $ is not refinable, and therefore $\tau$ is not   refinable. Since $\ker \tau$ contains a real basis for $\tau^{\vdash}$, $\tau \otimes 1$ is a good trace ofÊ $G \otimes \Q$.

Now we analyze properties of the individual pure traces and their finite subsets. Among the cases considered, only one of the quotients $G/\ker U$ is holey (not unperforated). This is probably atypical.

 \itemitem{$U = \brcs{\tau_1}$} We note that  $\ker \tau_1 = \langle v_1, v_2, v_5 - v_3\rangle$; the closure is $v_1 \R \oplus (v_5 - v_3) \Z$, and it is easy to check that regardless of the value of $a$, $J = \Cal P(J)$, and thus by Corollary \sevfiv(b), $G/\ker \tau_1$ is unperforated. However, $\ker \tau_1$ does not contain a real basis for $\tau_1^{\perp}$, so even $\tau_1 \otimes 1$ on $G\otimes \Q$ is not refinable.

 \itemitem{$U = \brcs{\tau_2}$} We have $\ker \tau_2 = \langle v_1, v_2, v_3 + v_6, v_5\rangle$; here $J \iso \R \oplus \Z^2$ and $ J = \Cal P (J)$, so Corollary \sevfiv(b) applies again; thus $G/\ker \tau_2 $ is unperforated. As in the case with $\tau$, $\tau_2$ is not refinable. 

\itemitem{$U = \brcs{\tau_3}$}  or  $\brcs{\tau_4}$ We note $\ker \tau_3 = \langle v_4, v_6, v_5 -v_1\rangle$ (since $\brcs{1,a,s}$ is rationally independent). Here $\ker \tau_3$ is discrete and contains a basis, so $\tau_3$ is ugly. It is routine to check that once again, $J = \Cal P(J)$, so the quotient is unperforated, and thus $\tau_3$ is not refinable. Similar computations reveal that $\tau_4$ has the same properties.
\itemitem{$U = \brcs{\tau_1, \tau_2}$} Then  $J = v_1 \R$, so $J = \Cal D(J)$ and thus $G/\ker U$ is unperforated. In this case, $Z(\ker U)$ is the face spanned by $U$, and since $\ker U$ is not dense in $U^{\perp} \iso \R^2$, $U$ is not good (nor refinable). 

\itemitem{$U = \brcs{\tau_2, \tau_3}$} Then  $\ker U = \langle v_5 - v_1 \rangle $, which is discrete of rank one, hence $U$ is not refinable; however, $J = \ker U = \Cal P(J)$, so $G/\ker U$ is unperforated.

\itemitem{$U = \brcs{\tau_3, \tau_4}$} Then $\ker U = \langle v_6 \rangle$, again discrete. Here $Z(\ker U)$ is again the affine span of $U$, and it follows from the rank of $\ker U$ that $U$ is not refinable. However, there are now two cases to consider.

\itemitem{} If $c > 0$, then $J = \Cal P(J)$, so $G/\ker U$ is unperforated. On the other hand, if $c < 0$, there is a trace in the interior of the trace space that kills $\ker U$, and by Proposition  \sevfou\ below, $G/\ker U$ is {\it not\/} unperforated.

\itemitem{$U = \brcs{\tau_2, \tau_4}$} We see $\ker U = \langle v_6 + v_3, v_5 + v_1\rangle$ which is discrete and equal to $P(J)$. Hence $\ker U$ is unperforated, and thus $U$ is not refinable; here $U$ is ugly.

\itemitem{$U = \brcs{\tau_1, \tau_3}$} or  $\brcs{\tau_1, \tau_4}$ We observe that  $\ker U = \brcs{0}$, and thus $Z(\ker U) = K$, hence $U$ is refinable.
\qed

\comment
The proof of the following does not require the interpolation property. The condition that the closure of $\ker U$ be  a real vector space can be restated simply:
\item{$\bullet$} for all $k \in \ker U$, the element $\hat k/2 \in \Aff S(G,u)$ can be uniformly approximated by the images of elements of $\ker U$. 
 
\Lem Proposition \sevthr. Let $(G,u)$ be an unperforated  simple patially ordered abelian group. If $U$ is a subset of $K:= S(G,u)$ \st the closure of the image of $\ker U$  in $Z(\ker U)^{\vdash} \subset \Aff S(G,u)$ is a real vector space, then $G/\ker U$ is unperforated.

\Pf Suppose $m(g + \ker U) \geq 0$ for some positive integer $m$ and element $g$ of $G \setminus \ker U$. This means there exists
$k \in \ker U$ and $p \in G^+$ \st $mg = p + k$. The element $\hat k/m$ of $\Aff K$ is in the rational, hence real, vector space generated by the image of $\ker U$, and since the closure of the latter is a real vector space containing the image, it follows, given $\epsilon$, there exists $k' \in \ker U$ \st $\| \hat k/m - \hat k' \| < \epsilon$. Set $g' = g - k'$, so that $\| m \hat g' - \hat p\| < m\epsilon$. Since $G$ is simple, $p$ is an order unit, so that $\delta:= \inf_{\sigma \in K}(\sigma (p)) > 0$. If we select $\epsilon < \delta/m$, then $m \hat g' \gg \delta - m\epsilon > 0$, that is $m\hat g'$ is strictly positive, and therefore $\hat g' \gg 0$. As $G$ is unperforated, $g' $ is an order unit. Hence $g + \ker U = g' + \ker U$ is in $(G/\ker U)^+$.
\qed

Together with  Proposition \appbsix, it follows that if $\ker U$ has dense image in $Z(\ker U)^{\vdash}$ (or even with the weaker hypothesis, the closure of the image of $\ker U$ is a real vector space) and $U$ is refinable, then $G/\ker U$ is unperforated, and thus $Z(\ker U)$ is good. This means that if $U$ is refinable but $Z(\ker U)$ is not good, then not only is $G/\ker U$ not unperforated, but the image of $\ker U$ is far from being dense in $Z(\ker U)^{\vdash}$ (and see the next two results).

There are however other situations wherein the quotient can be unperforated, and the closure of the image is not a real vector space. Here is an easy-to-prove case; when it applies, 
$Z(\ker U)$ is a face of $S(G,u)$. This situation really does occur in some simple dimension groups, despite its simplicity.

\Lem Lemma \sevtwo. Let $(G,u)$ be a simple dimension group, and let $U$ be a subset of $S(G,u)$ with the property that the image of $\ker U$ in $\Aff S(G,u)$ is generated as an abelian group by elements of the form $\Set{\hat g \in \Aff S(G,u)}{\hat g|U = 0 \ \&\  \hat g \geq 0}$. Then $G/\ker U$ is unperforated and $Z(\ker U)$ is a closed face of $S(G,u)$. 

\Remark Generation by nonnegative elements really means that in the image of $G$, equipped with the archimedean ordering, that is, as functions on $S(G,u)$, the image of $\ker U$ is an order ideal. Of course, neither $\ker U$ nor its image is an order ideal \wrt the simple orderings. 

\Pf Suppose $mg + \ker U \in (G/\ker U)^+$ for $g \in G\setminus \ker U$ and $m$ a positive integer. This entails there exists $p$ in  $G^+ \setminus\brcs{0}$ and $k$ in $\ker U$ \st $mg + k = p$. Thus $m\hat g + \hat k = \hat p \gg 0$. From the hypothesis, there exist $k_1, k_2$ in $K$ \st $\hat k_i \geq 0$ and $\hat k = \hat k_1 - \hat k_2$. Thus $m\hat g + \hat k_1 = \hat p  + \hat k_2$, so that $m\hat g + \hat k_1 \geq \hat p \gg 0$. Hence $m(\hat g + \hat k_1) \gg 0$ (since $\hat k_1 \geq 0$), and thus $\hat g + \hat k_1 \gg 0$. Since $G$ is unperforated, this forces $g + k_1$ to be an order unit in $G$, and thus $g + \ker U = g+k_1 + \ker U$ is in $(G/\ker U)^+$. That $Z(\ker U)$ is a face is immediate.\qed

For example, 
The next result  includes the two preceding ones, but the combined criteria lead to an awkward hypothesis. If $J$ is a closed subgroup of a Banach space, then the divisible subgroup of $J$ (the maximal rational subspace of $J$), $D(J)$, is a closed real vector space and the maximal real vector subspace of $J$. If $J$ is the closure of an abelian group of of finite rank or more generally, if $J\R$ is finite dimensional, then $D(J)$ also open, but this does not extend to the infinite rank case. For $J$ a closed subgroup of $\Aff K$, let $\Cal P(J)$ be the subgroup of $J$ generated by $\Set{j \in J}{j \geq 0}$. We now consider the case that the closure of $\ker U$ in $\Aff K$, $J$, contain $D(J) + \Cal P(J)$ as a dense subgroup. The previous two results considered the cases where $J = D(J)$ and a special case of $J = \Cal P(J) $ respectively.

 \Lem Proposition \sevfiv. Let $(G,u)$ be a simple dimension group, let $U$ be a subset of $K = S(G,u)$, and denote the closure of the image  of $\ker U$ in $\Aff K$ by $J$. If $ D(J) + \Cal P(J)$ is dense in $J$, then $G/\ker U$ is unperforated. 

\Pf Obviously $G/\ker U$ is torsion-free. Suppose $mg + \ker U $ is nonzero and in the positive cone of $G/\ker U$ for some $g \in G$ and positive integer $m$; then there exist nonzero $p$ in $G^+$ and $k \in \ker U$ \st $mg = p -k$. By hypothesis, given  $\epsilon < \inf_{\sigma\in K}\sigma(p)/2$, there exist $k_1, k_2, k_3$ in $\ker U$ \st $\hat k_2, \hat k_3 \geq -\epsilon \pmb 1/2m$ and $\| \hat k - (m\hat k_1 - \hat k_2 + \hat k_3) \| < \epsilon$. Thus for any $\sigma \in K$, $\sigma( p - k  -( mk_1 - k_2 + k_3)) \geq \sigma (p) - \sigma (k -( mk_1 + k_2 - k_3) )  > \sigma(p) - \epsilon > \frac 12 \sigma(p) > 0$. Hence $p' = p - k  -( mk_1 - k_2 + k_3)$ is an order unit of $G$. 

Now $mg = p' +k_3 - k_2 + mk_1$. In particular, $m(g -k_1) +k_2 = p' + k_3$. 
We observe, for any $\sigma \in K$,
$$\eqalign{
m\sigma(g - k_1 + k_2) & = (m-1) \sigma(k_2) + \sigma(p') + \sigma(k_3) \cr
& >  \frac{-(m-1)\epsilon }{2m} + \frac {\sigma(p)}{2} + \frac{-\epsilon}{2m} = \frac{\sigma(p) - \epsilon}2  > 0\quad \text{and thus,}\cr
\sigma(g - k_1 + k_2) & > 0.
}$$
 By unperforation of $G$, $g - k_1 + k_2$ is an order unit of  $G$, so that its image in $G/\ker U$, $g+\ker U$, is positive.
\qed

It is likely that the converse is true as well; the difficulty lies in the fact that $\Cal P(J)$ need not contain any nonzero elements of the image of $\ker U$. It leads to a partial confirmation (Corollary \sevsix) of the conjecture that $U$ is refinable implies $Z(\ker U)$ is good; unfortunately, the requirement that $Z(\ker U)$ be a face excludes some of the other constructions.
\endcomment

In contrast, here is a reasonably general result on the failure of unperforation for   quotients by convex subgroups.
It does not require simplicity, merely approximate divisibility. Compare with Lemma \appbele.

\Lem Proposition \sevfou. Suppose $(G,u)$ is an approximately divisible dimension group, with normalized trace space $K = S(G,u)$. Let $U$ be a subset of $K$ for which $Z = Z(\ker U)$ satisfies $Z \cap \text{Int}(K) \neq \emptyset $. If there exists a trace $\tau$ \st $\tau(\ker U) $ is nonzero and cyclic, then $G/\ker U$ is not unperforated.

\Remark The condition that $Z$ meet the interior of $K$ obviously rules out faces; however, there are good $Z$ (\wrt $\Aff K$) that contain an interior point. 

\Remark The proof actually exhibits a great deal of holeyness among the potential order units. 

\Pf There exists $z \in Z $ and in the interior of $K$; hence the line segment joining $\tau$ to $z$ extends outside $Z$ but still within $K$. The affine span of this segment hits the boundary of $K$  at $\tau_2 \in \partial K$ (neither $\tau_2$ nor $\tau$ need be pure). We thus have a decomposition $z = \lambda \tau  + (1-\lambda)\tau_2$ for some $0 < \lambda < 1$. There exists $\delta \neq 0$ \st $\tau(\ker \tau) = \delta \Z$; without loss of generality, we may assume $\delta > 0$. 
For any $v \in \ker U$, $z(v) = 0$, so that $\tau_2 (v) = -\lambda \tau(v)/(1-\lambda)$; in particular, $\tau_2$ restricted to $\ker U$ has discrete (and nonzero) range. 

Fix an element $k \in \ker U$ \st $\tau(k) = \delta$. Since $\tau_2 (k) < 0$, $k$ is not positive. Now we find $g$ in $G$ \st $g+k \gg 0$ but $\tau(g) < 0$. By density of the range of $G$ in $\Aff K$, there exists $j \gg 0$ \st 
$0 \ll \hat j \ll (\delta/2)\pmb 1$, and set $g = j -k$. Then $\tau(g) = \tau(j) - \delta < - \delta/2$. 

As $\tau_2 (k ) < 0$ and $g + k$ is an order unit, it follows that $\tau_2(g) > 0$. If $m$ is a sufficiently large positive integer, then 
$$
\frac{|\tau(g)|}m < \frac \delta3 \quad \text{ and }\quad \frac{\tau_2(g)}m < \frac {\lambda\delta}{3(1-\lambda)}.
  $$
Choose positive real $\epsilon$ small enough that
$$
\epsilon < \min \brcs{\frac{m\delta}2,\frac{m\lambda\delta}{2(1-\lambda)},\frac{|\tau(g)|}2, \frac{\tau_2(g)}2,\min_{\sigma \in K} \sigma (g+k)}.
$$
Again using density of $G$ in $\Aff K$, we may find $g_1 \in G$ \st $\|\hat g_1 - \hat g/m\| < \epsilon/m$. Then $\| m\hat g_1 + \hat k - (\hat g + \hat k)\| < \epsilon < \min_{\sigma \in K} \sigma (g+k)$. Hence $mg_1 + k \gg 0$, so $mg_1 + k$ is an order unit in $G$, and thus  $mg_1 + \ker U$ is an order unit (in particular, it is positive) in $G/\ker U$. 

Now we show $g_1 + \ker U$ is not in the positive cone of $G/\ker U$. Suppose to the contrary that it is. If  $b:= g_1 + k_1 $ is a positive element of $G$ but not an order unit for some $k_1 \in \ker U$, we obtain an immediate contradiction: for no positive integer $l$ can $ u \leq l b$---but this contradicts $mb$ being an order unit. Hence  $g_1 + k_1 \gg 0$. There exists an integer $n$ (possibly negative) \st $\tau(k_1) = n\delta$; then $\tau_2 (k_1) = -\lambda n \delta/(1-\lambda)$. Applying $\tau$ and $\tau_2$ to $g_1 + k_1$, we obtain
$$
\tau(g_1) + n\delta > 0 \quad \text{ and }\quad \tau_2 (g_1) - \frac {n\lambda\delta}{(1-\lambda)}> 0.
$$
Since $0 < \lambda < 1$ and $\delta > 0$, this yields
$$
\frac{\tau_2 (g_1)(1-\lambda)}{\lambda \delta} > n > -\frac{\tau(g_1)}\delta.
$$
Since $\| \hat g_1 - \hat g/m \| < \epsilon/m$, it follows that $|\tau_2 (g_1) - \tau_2 (g)/m |< \epsilon/m$ and $|\tau (g_1) - \tau (g)/m |< \epsilon/m$. The latter yields $\tau(g_1) < \tau(g)/m + \epsilon/m < \tau(g)(1 -1/2) /m < 0$, and the former yields
$$
\tau_2(g_1) < \frac{\epsilon} m + \frac{\tau_2 (g)}m < \frac{\epsilon}m + \frac{\lambda \delta}{3(1-\lambda)} <  \frac 56 \frac{\lambda \delta}{1-\lambda}.
$$
Therefore $n < 5/6$, so $ n \leq 0$. However, $-\tau(g_1) > 0$, so $n > 0$, a contradiction. 
\qed

Altering some of the numbers can  increase the set of integers $r$ \st $rg_1 + \ker U$ is not positive. 

By similar methods, we can at least obtain a necessary condition for refinable sets, which together with the previous results is tantalizingly close to the conjectural result that refinability of $U$ implies $G/\ker U$ is unperforated (which would imply that $Z(\ker U)$ is good).

\Lem Proposition \sevnin. Suppose that $(G,u)$ is an approximately divisible dimension group, and $U$ is a refinable subset of $S(G,u)$. If $\sigma$ is a trace of $G$, then either $\sigma(\ker U)$ is zero or dense in $\R$. 

\Pf If not, there exists $\tau \in S(G,u)$ \st $\tau (\ker U)= \delta \Z$ for some positive real number $\delta$. Select $k \in \ker U$ \st $\tau(k) = \delta$. There exist (lots of) order units $v$ \st $v\gg k$ (any sufficiently large integer multiple of any order unit will satisfy this); pick one of them, call it $v$, and suppose $\tau(v)/\delta = r > 0$. Select an integer $N > r$. Now choose positive $\epsilon$ \st 
$$
\epsilon < \min \brcs{(N-r)\delta, \inf_{\sigma \in S(G,u)} \sigma(v), \inf_{\sigma \in S(G,u)} \sigma(v-k)}.
$$
By density of the image of $G$ in $\Aff S(G,u)$, we may find $v'$ in $G$ \st 
$\|\hat v' - \hat v/N\| < \epsilon/N$. For any trace $\sigma$, $\sigma (v') > (\sigma (v) - \epsilon)/N > 0$. Hence $v'$ is an order unit. In addition, $\sigma(Nv' - k) > \sigma (v -k) - \epsilon > 0$, so that $Nv' -k$ is also an order unit. Moreover, $\tau(v') < (\tau(v) + \epsilon)/N = (r\delta + \epsilon)/N < \delta$. Define $a_i = v'$ for $i=1,2,\dots, N$ and set $b = Nv'-k$. Then $b = \sum a_i - k$ where $b,a_i$ are all order units, and $k \in \ker U$. By refinability, there exist $c_i \in G^+$ \st $c_i - v' = k_i \in \ker U$ and $b = \sum c_i$. Applying $\tau$, we have the following (with integer $n(i)$),
$$\eqalign{
\tau(c_i ) &= \tau(v') + \tau(k_i) = \tau(v') + n(i)\delta \cr
N\tau(v') - \delta & = \tau(b) = \sum \tau(c_i)\cr
}$$
This yields $\sum n(i) = -1$. Hence at least one of the $n(i)$, call it $n(j)$, is less than or equal $-1$. However, $\tau(c_j) = \tau (v') + n(j)\delta \leq \tau (v') - \delta < 0$, contradicting positivity of $c_j$.
\qed 

\Lem Corollary \sevsix. Let $(G,u)$ be a simple dimension group with finitely many pure traces, and let $U$ be a refinable subset of $S(G,u)$ \st any one of the following conditions holds: {\par} \item{(i)}  $Z \equiv Z(\ker U)$ is a face;
\item{(ii)} $Z$ contains an interior point of $K = S(G,u)$;
\item{(iii)} $Z(\ker U)$ is contained in a facet, but contains a relative interior point of that facet.{\par } \noindent Then $Z(\ker U)$ is good.

\Pf  Since $U$ is refinable for $G$, so is it refinable for $G' = G\otimes \Q$. In this case, $G'/\ker_{G'} U$ is unperforated (since it is a rational vector space), and since  $Z(\ker U)$ is refinable for $G'$, from Proposition \appbsix(e), the image of $\ker_{G'} U = (\ker U) \otimes \Q$ is dense in $Z^{\vdash}$.

Let $J$ denote the closure of the image of $\ker U$. Since the pure trace space is finite dimensional, $J$ can be written as $D(J) \oplus W$
where $W$ is a discrete abelian group (therefore, free). Since $J\R = Z^{\vdash}$, any $\Z$-basis for $W$ together with a real basis for the real vector space $D(J)$ will constitute a real basis for $Z^{\perp}$. By Proposition \appbsix(f), it suffices to show that $G/\ker U$ is unperforated. 

\noindent (i) Now  suppose $Z$ is a face. We claim that $W \subseteq P(J)$. The set of elements of $Z^{\perp}$ that are nonnegative (as functions on $S(G,u)$) forms a closed cone in $Z^{\perp}$ that contains an open ball (true for any closed face). By taking a ball of sufficiently large radius inside this cone, we can guarantee that it contains a fundamental parallelotope of $W$, that is, a $\Z$-basis of $W$. Thus $W \subseteq P(J)$. 

Hence $J = D(J) + P(J)$, so by Proposition \sevfiv, $G/\ker U$ is unperforated, and this case is done.

\noindent (ii) Now assume that $Z$ contains an interior point, $z$, of $K = S(G,u)$. If $W$ is not zero, there exists a bounded linear functional on $\Aff K$, $\rho$, \st $\rho(\ker \tau) = \delta\Z $ for some $\delta > 0$. Since $z$ is an interior point of $K$, there exists a positive integer $N$ \st $\tau:= \rho + Nz$ is a trace (a positive linear combination of the pure traces). Clearly $\tau (\ker U) = \rho(\ker \tau)$, so $\tau$ has nonzero cyclic image. By Proposition \sevnin, this contradicts refinability of $U$. Hence $W = \brcs{0}$ and thus $J $ is a real vector space. By Proposition \sevfiv, $G/\ker U$ is unperforated.

\noindent (iii) Let the facet be given by vertices $\brcs{\tau_0, \dots, \tau_{n-1}}$, and let $\tau_n$ denote the complementary vertex. We can write $\sigma = \sum_{i= 0}^nÊ r_i \tau_i$. The relative interior point of the facet belonging to $Z(\ker U)$ is written in the form $\phi= \sum_{i=0}^{n-1} R_i \tau_i$ where $R_i > 0$. 

By multiplication by $-1$ if necessary, we may assume $r_n \geq 0$. Then there exists $N$ \st all the coefficients of $\tau:= \sigma + N \phi$ are nonnegative, and thus $\tau$ is a trace. Since $\phi(\ker U) = 0$, it follows that $\tau (\ker U) =\Z \delta$. Since $\tau$ is a trace, Proposition \sevnin\ applies, contradicting refinability. Hence the discrete part of $J$ is zero, so $J$ is divisible, thus (as it contains a real basis for $(\ker U^{\vdash}$ and is complete), $J = (\ker U)^{\vdash}$ and so $Z(\ker U)$ is good.
\qed

Combining this with the earlier results, we see that when $G$ has finitely many pure traces and $Z(\ker U)$ is a face or contains an interior point, then $U$ is refinable if and only if $\ker U$ has dense image in $Z(\ker U)^{\perp}$. In particular, if $\tau$ is a trace for which we could verify that  $Z(\ker \tau)$ is a face, then $\tau$ is refinable if and only if $\ker \tau$ is dense $Z(\ker \tau)^{\vdash}$.
This excludes all the pure traces of our now standard example, $G=\langle e_i; \sum \alpha_i e_i \rangle$, from being refinable. The first case for which the conjecture that refinability of $U$
 implies goodness of $Z(\ker U)$ is unknown occurs when $U$ is a singleton in the relative interior of an edge of a tetrahedron. Along these lines is the following supplementary result. 

\Lem Lemma \sevthi. Suppose $(G,u)$ is a simple dimension group with finite dimensional trace space $K$, $F$ is a proper face ofÊ $K$, and $U$ is a refinable subset for $G$. If $F$ is the smallest face containing $U$ and $F$ is good (for $G$), then $Z(\ker U)$ is good.

\Remark The hypothesis that $F$ be a good subset for $G$ is not entirely satisfactory---what would be preferable is merely that $F$ be refinable (which would imply that $Z(\ker F)$---which can strictly contain $F$ even though the latter is a face---is good). However, we were not able to prove this. Likely sufficient is that $F$ (the face generated by $U$) be of codimension one in a good face.

\Pf This is by induction on the dimension of $K$.
\comment
ÊList the extreme points of $F$ as $\brcs{\tau_i}_{i=1}^s$, and define the map $\Arrow p; G. \R^s$ via $p = (\tau_i)$. Since $G$ has dense range in $\Aff K$, it easily follows that $H = p(G)$ is a dense subgroup of $\R^s \iso \Aff F$. Hence with the strict ordering, $H $ is a simple dimension group, and its trace space is $F$.
Let $U'$ be the image of $U$ in the trace space ofÊ $H$. If $H/\ker$
Now $U \subseteq F$. We show that it is refinable \wrt $H$. Suppose $b,a_i' \in H^+$, andÊ $b - a_i \in \ker_H U$. Since $\ker F \subseteq \ker U$, we can apply ref inability of $F$, to obtain $$
\endcomment

Since $F$ is good, $G_1:= G/\ker F$ is itself a simple dimension group, and of course, its trace space is naturally homeomorphic to $F$. Now
$\ker F \subseteq \ker U$; let $V$ be the image of $\ker U$ in $G_1$ (it is $\ker U$ relative to $G_1$). Then $G/\ker U \iso G_1/V$. If $a_i,b \in G_1^+$, and $b = \sum a_i + k$ with $k \in V$, we can lift (since $\ker F$ is good) each toÊ $A_i, B \in G^+$ \st $A_i + V = a_i$ and $B + V = b$, and then $BÊÊ -\( \sum A_i \) \in \ker U$. By ref inability of $U$, there exist $A_i' \in G^+$ \st $A_i' - A_i \in \ker U$ and $B = \sum A_i'$. Setting $a_i'Ê = A_i' + V$, we see that $U$ is refinable for the dimension group $G/\ker F$. 

We may replace $U$ by $Z(\ker U) \cap K$; the latter is compact convex, and since $F$ is the smallest face generated by $U$, $L(U)$ contains a point in the relative interior of $F$. Since $U$ is refinable relative to $G_1$, now Lemma \sevsix(ii) applies, so that $Z(\ker U)$ is good \wrt $G_1$. Thus $G_1/V$ is unperforated, and because of the order isomorphism with $G/\ker Z(\ker U)$, the latter is unperforated. Finally, by \appbsix(f), $U$ refinable
and $G/\ker U$ unperforated imply $Z(\ker U)$ is good. \qed


\long\def\Rf[#1] #2, #3. #4\par%
{\vskip 2pt \itemitem{[#1]} #2, {\it #3,} #4\par\vskip2pt}
\def\Ibid{---\!---\!---}

\SecT References

\Rf [Ak1] E Akin, Measures on Cantor space. Topology Proc,  {24} (1999)  1--34.

\Rf [Ak2] \Ibid, Good measures on Cantor space. Trans Amer Math Soc,  {357} (2005)  2681--2722.

\Rf [ADMY] {E Akin, R Dougherty,  RD Mauldin, and A Yingst}, Which Bernoulli measures are good measures?. Colloq  Math,  {110} (2008)  243--291.

\Rf [Al] EM Alfsen, 	
Compact convex sets and boundary integrals. Springer--Verlag 1971, ix + 210\, pp.

\Rf [AE] L Asimow \& AJ Ellis, Convexity theory and its applications in functional analysis. London Mathematical Society Monographs 16 (1980).

\Rf [Au] TD Austin, A pair of non--homeomorphic product measures on the Cantor set. Math Proc Cam Phil Soc,  {142} (2007) 103--110.

\Rf [BK] {S  Bezuglyi and O Karpel}, Homeomorphic measures on stationary Bratteli diagrams. J Funct Analysis, 261 (2011) 3519--3548.

\Rf [Bk2] {S Bezuglyi and J Kwiatkowski}, Topological full group of a Cantor minimal system is dense in the full group. Topological Methods in Nonlinear Analysis 16 (2000) 2, 371--397.

\Rf [BH1] Mike Boyle \&  David Handelman, Entropy versus orbit equivalence for minimal homeomorphisms.  Pacific J Math  164  (1994),  no. 1, 1--13.

\Rf [BH2] \Ibid, Ordered equivalence, flow equivalence and ordered cohomology. Israel J Math 95 (1996) 169--210.

\Rf [BH3] \Ibid, unpublished drafts. 1993--2011. 

\Rf [DMY] {R Dougherty, R Daniel Mauldin, and A Yingst}, On homeomorphic Bernoulli measures on the Cantor space. Trans Amer Math Soc 359 (2007), 6155--6166

\Rf [E] EG Effros, Dimensions and C*-algebras. CBMS Regional Conference Series in Mathematics 46, Conference Board of the mathematical sciences, Washington (1981).

\Rf[EHS] {EG Effros,   David Handelman,   \&  Chao Liang Shen}, Dimension groups and their affine representations.  Amer  J  Math  102  (1980), no. 2, 385--407.

\Rf [GPS1] {T Giordano,  IF Putnam, \& CF Skau},  Topological orbit equivalence and $C^*$--crossed products.  J Reine Angew Math  469  (1995), 51--111.

\Rf [GPS2] \Ibid, Full groups of Cantor minimal systems.  Israel J Math.  111  (1999), 285--320.

\Rf [GW] {E  Glasner and B Weiss}, Weak orbital equivalence of minimal Cantor systems. Internat J Math, 6 (1995) 559--79.

\Rf [G] KR Goodearl, Partially ordered abelian groups with interpolation. Mathematical Surveys and Monographs, 20, American Mathematical Society, Providence, RI, 1986, xxii + 336\, pp.

\Rf [GH1] KR Goodearl \&  David Handelman,  Metric completions of partially ordered abelian groups.  Indiana Univ Math J  29  (1980), no\. 6, 861--895.

\Rf [GH2] \Ibid,   Tensor products of dimension groups and $K_0$ of unit--regular rings.  Canad. J. Math.  38  (1986),  no\. 3, 633--658.

\Rf [GH3] \Ibid,   Stenosis in dimension groups and AF C*-algebras.  Crelle's J   332  (1982),  1--98.

\Rf [H] David Handelman, Extensions for AF C*-algebras and dimension groups. Trans Amer Math Soc  (1982) 537--573.

\Rf [H1] \Ibid, Free rank $n+1$ dense subgroups of $\text{\/\bf R}^{n}$ and their endomorphisms.  J Funct Anal  46  (1982), no\. 1, 1--27.

\Rf [H2] \Ibid, Imitation product type actions on UHF algebras. J Algebra 99 (1986) 1--21.

\Rf [H3] \Ibid,   Positive polynomials, convex integral polytopes, and a random walk problem. Lecture Notes in Mathematics, 1282, Springer--Verlag, Berlin, 1987, xii+136 pp.

\Rf [H4] \Ibid, Representing polynomials by positive linear functions on compact convex polyhedra.  Pacific J Math  132  (1988),  no\. 1, 35--62.

\Rf [H5]  \Ibid, Iterated multiplication of characters of compact connected Lie groups.  J Algebra  173  (1995),  no\. 1, 67--96.

\Rf [H6] \Ibid, Simple archimedean dimension groups. Proc Amer Math Soc, to appear.

\Rf [HPS] {RH Herman, IF Putnam, \& CF Skau}, Ordered Bratteli diagrams, dimension groups and topological dynamics.  Internat J Math  3  (1992),  no\. 6, 827--864.

\Rf [K] O Karpel,  Infinite measures on Cantor spaces. J Difference Equ\. Appl\.,
(to appear).

\Rf [KRW] {KH Kim, FW Rousch, SG Williams}, Duality and its consequences for ordered cohomology of finite type shifts. Combinatorial and computational mathematics, Pohang (Korea), 2001.

\Rf [L] John Lawrence, A simple torsion-free perforated interpolation group. unpublished, {\it ca\,}1981.

\Rf [L2] \Ibid, Countable abelian groups with a discrete norm are free. Proc Amer Math Soc 90 (1984) 352--354. 

\Rf [M] Hiroki Matui, Dimension groups  topological joinings and non--coalescence of Cantor minimal systems.  Pacific J Math  204  (2002),  no\. 1, 163--176.

\Rf [Me] K Medynets, Cantor aperiodic systems and {Bratteli} diagrams.
CR Math, Acad Sci Paris, 342 (2006), no\. 1, 43--46.

\Rf [P] YT Poon, A K-theoretic invariant for dynamical systems. Trans Amer Math Soc 311 (1989) 2 515--533.

\Rf [V] AM Vershik, Uniform algebraic approximation of shift and multiplication operators.
Dokl Akad Nauk SSSR (259) 526--529, 1981 (Russian).

\Rf [Pu] IF Putnam, The C* algebras  associated with minimal  homeomorphisms of the Cantor set.
Pacific J Math 136 (1989) 329--353.

\Rf [Y] Andrew Q Yingst, A characterization of homeomorphic Bernoulli trial measures. Trans  Amer Math  Soc 360 (2008), 1103--1131.

{}

\vskip 10pt
\noindent Authors' affiliations \vskip 2pt

SB: Institute for Low Temperature Physics, 47 Lenin Avenue, 61103 Kharkov, Ukraine; bezuglyi\@ilt.kharkov.ua

DH: Mathematics Department, University of Ottawa, Ottawa ON  K1N 6N5, Canada; dehsg\@uottawa.ca

\end